\documentclass[twoside]{article} 
\usepackage[accepted]{aistats2017}
\pdfoutput=1 %

\usepackage[round]{natbib}
\usepackage{multibib} %
\newcites{sup}{Supplementary References}

\usepackage{amsthm}
\usepackage[leqno]{amsmath}
\usepackage{caption}
\usepackage{subcaption}
\usepackage{bm,color}
\usepackage{mathdots}
\usepackage{booktabs}       %
\usepackage{amsfonts}       %
\usepackage{nicefrac}       %
\usepackage{microtype}      %
\usepackage{graphicx}
\usepackage{hyperref}       %
\usepackage{algcompatible} %
\usepackage{algorithm}
\usepackage{pict2e}
\usepackage{blindtext}

\hypersetup{ %
    pdftitle={Frank-Wolfe Algorithms for Saddle Point Problems},
    pdfkeywords={},
    pdfborder=0 0 0,
	pdfpagemode=UseNone,
    colorlinks=true,
    linkcolor=blue, %
    citecolor=blue, %
    filecolor=blue, %
    urlcolor=blue, %
    pdfview=FitH,
    pdfauthor={Gauthier Gidel, Tony Jebara, Simon Lacoste-Julien},
}

\newcommand{\N}{\mathbb{N}}
\newcommand{\R}{\mathbb{R}}
\renewcommand{\L}{\mathcal{L}}
\newcommand{\X}{\mathcal{X}}
\newcommand{\Y}{\mathcal{Y}}
\newcommand{\Z}{\mathcal{Z}}
\newcommand{\M}{\X \times \Y}
\newcommand{\prodscal}[2]{\left\langle#1,#2\right\rangle}
\newcommand{\x}{\bm{x}}
\newcommand{\y}{\bm{y}}
\newcommand{\z}{\bm{z}}
\renewcommand{\u}{\bm{u}}
\newcommand{\s}{\bm{s}}
\newcommand{\xt}{\bm{x}^{(t)}}
\newcommand{\xtt}{\bm {x}^{(t+1)}}
\newcommand{\zt}{\bm{z}^{(t)}}
\newcommand{\ztt}{\bm {z}^{(t+1)}}
\newcommand{\st}{\bm{s}^{(t)}}
\newcommand{\stt}{\bm{s}^{(t+1)}}
\newcommand{\vt}{\bm{v}^{(t)}}
\newcommand{\dt}{\bm{d}^{(t)}}
\newcommand{\rt}{\r^{(t)}}
\newcommand{\yt}{\bm{y}^{(t)}}
\newcommand{\ytt}{\bm{y}^{(t+1)}}
\newcommand{\xtm}{\widehat{\bm{x}}^{(t)}}
\newcommand{\ytm}{\widehat{\bm{y}}^{(t)}}
\newcommand{\wt}{w_t}
\newcommand{\wtt}{w_{t+1}}

\newcommand{\gnx}{\nabla_x\L(\xt,\yt)}
\newcommand{\gny}{\nabla_y\L(\xt,\yt)}

\renewcommand{\SS}{ \bm{S}}

\newcommand{\gap}{g_{t}}

\newcommand{\CondNumb}{\nu}

\newcommand{\FF}{\r} %

\newtheorem*{rep@theorem}{\rep@title}
\newcommand{\newreptheorem}[2]{%
\newenvironment{rep#1}[1]{%
 \def\rep@title{#2 \ref{##1}}%
 \begin{rep@theorem}}%
 {\end{rep@theorem}}}

\newreptheorem{lemma}{Lemma'}
\newreptheorem{proposition}{Proposition'}
\newreptheorem{theorem}{Theorem'}
\newreptheorem{remark}{Remark'}

\makeatletter
\newcommand{\leqnomode}{\tagsleft@true\let\veqno\@@leqno}
\newcommand{\reqnomode}{\tagsleft@false\let\veqno\@@eqno}
\makeatother

\newtheorem{definition}{Definition}
\newtheorem{theorem}[definition]{Theorem}
\newtheorem{proposition}[definition]{Proposition}
\newtheorem{lemma}[definition]{Lemma}
\newtheorem{corollary}[definition]{Corollary}

\DeclareMathOperator*{\conv}{conv}

\DeclareMathOperator*{\argmin}{\arg\min}
\DeclareMathOperator*{\argmax}{\arg\max}

\providecommand{\norm}[1]{\left\lVert#1\right\rVert}
\providecommand{\normEucl}[1]{\left\lVert#1\right\rVert_2}

\newcommand{\domain}{\mathcal{\X}} %
\newcommand{\stepsize}{\gamma}
\newcommand{\stepmax}{\stepsize_{\textnormal{\scriptsize max}}} %
\newcommand{\FW}{{\hspace{0.05em}\textnormal{FW}}}
\newcommand{\PW}{{\hspace{0.05em}\textnormal{PFW}}}
\newcommand{\away}{{\hspace{0.05em}\textnormal{A}}}
\newcommand{\Cf}{C_{\hspace{-0.08em}f}}

\newcommand{\strongConvAFW}{\mu_{\hspace{-0.08em}f}^\away}

\newcommand{\dirW}{\mathop{dirW}}
\newcommand{\dd}{\bm{d}}

\newcommand{\vv}{\bm{v}} %

\DeclareMathOperator*{\lmo}{LMO_{\!\Vertices}}
\newcommand{\A}{\mathcal{A}}
\newcommand{\B}{\mathcal{B}}

\newcommand{\Vertices}{\mathcal{A}\times \mathcal{B}} %
\newcommand{\Coreset}{\mathcal{S}}
\renewcommand{\SS}{\mathcal{S}}
\renewcommand{\aa}{\bm{\alpha}}
\renewcommand{\r}{\bm{r}}
\newcommand{\PdirW}{\mathop{PdirW}}
\newcommand{\PWidth}{\mathop{PW\!idth}}
\newcommand{\innerProd}[2]{\left\langle #1 , #2 \right\rangle}
\newcommand{\innerProdCompressed}[2]{\langle #1 , #2 \rangle}

\newcommand{\Kface}{\mathcal{K}}

\newcommand{\xc}{\x_c} %
\newcommand{\yc}{\y_c} %
\newcommand{\muIntL}{\mu^{\text{int}}_\L} %

\newcommand{\0}{\mathbf{0}} %

\newcommand{\ignore}[1]{}%

\newcommand{\remove}[1]{} %

\begin{document}
\reqnomode %

\twocolumn[

\aistatstitle{Frank-Wolfe Algorithms for Saddle Point Problems}

\aistatsauthor{Gauthier Gidel
\And
Tony Jebara
\And
Simon Lacoste-Julien}

\aistatsaddress{ INRIA - Sierra Project-team\\
\'Ecole normale sup\'erieure, Paris 
\And Department of Computer Science\\
Columbia U. \& Netflix Inc., NYC
\And
Department of CS \& OR (DIRO) \\
Universit\'e de Montr\'eal, Montr\'eal} ]

\begin{abstract}
We extend the Frank-Wolfe (FW) optimization algorithm to solve constrained smooth convex-concave saddle point (SP) problems.
Remarkably, the method only requires access to linear minimization oracles. 
Leveraging recent advances in FW optimization, we provide the first proof of convergence of a FW-type saddle point solver over polytopes, thereby partially answering a 30 year-old conjecture.
We also survey other convergence results and highlight gaps in the theoretical underpinnings of FW-style algorithms. 
Motivating applications without known efficient alternatives are explored through structured prediction with combinatorial penalties as well as games over matching polytopes involving an exponential number of constraints. 
\end{abstract}

\section{Introduction}
The Frank-Wolfe (FW) optimization algorithm~\citep{Frank:1956vp}, also known as the conditional gradient method~\citep{demyanov1970approximate}, is a first-order method for smooth constrained optimization over a compact set. 
It has recently enjoyed a surge in popularity thanks to its ability to cheaply exploit the structured constraint sets appearing in machine learning applications~\citep{jaggi2013revisiting,lacoste2015global}. 
A known forte of FW is that it only requires access to a \emph{linear minimization oracle} (LMO) over the constraint set, i.e., the ability to minimize linear functions over the set, in contrast to projected gradient methods which require the minimization of \emph{quadratic} functions or other nonlinear functions.
In this paper, we extend the applicability of the FW algorithm to solve the following convex-concave saddle point problems:
\vspace{-2mm}  
\begin{equation} 
\label{eq:convex_concave}
   \min_{x \in \X} \max_{y \in \Y} \L(\x,\y), \vspace{-1mm}  
\end{equation}
\begin{equation*} 
  \text{with only access to }\text{LMO}(\mathbf{r}) \in \argmin_{\s \in \X \times \Y} \innerProd{\s}{\mathbf{r}}, \
\end{equation*}
where $\L$ is a smooth (with $L$-Lipschitz continuous gradient) \emph{convex-concave function}, i.e., $\L(\cdot,\y)$ is convex for all $\y \in \Y$ and $\L(\x,\cdot)$ is concave for all $\x \in \X$.
We also assume that $\X\times\Y$ is a convex compact set such that its LMO is cheap to compute. 
A \emph{saddle point solution} to~\eqref{eq:convex_concave} is a pair $(\x^*, \y^*) \in \X \times \Y$~\citep[VII.4]{hiriart2013convexI} such that: $\forall \x \in \X,\, \forall \y \in \Y$,
\begin{equation}\label{eq:saddle_point}
   \L(\x^*,\y) \leq \L(\x^*, \y^*) \leq \L(\x, \y^*).  
\vspace{-1mm}
\end{equation}
\paragraph{Examples of saddle point problems.} 
\citet{taskar2006dualExtrag} cast the maximum-margin estimation of structured output models as a bilinear saddle point problem $\L(\x,\y) = \x^\top M\y$, where $\X$ is the regularized set of parameters and $\Y$ is an encoding of the set of possible structured outputs. 
They considered settings where the projection on $\X$ and $\Y$ was efficient, but one can imagine many situations where only LMO's are efficient. 
For example, we could use a structured sparsity inducing norm~\citep{martins2011structuredSparsity} for the parameter~$\x$, such as the overlapping group lasso for which the projection is expensive~\citep{bach2012structuredSparsity}, while $\Y$ could be a combinatorial object such as a the ground state of a planar Ising model (without external field) which admits an efficient oracle~\citep{barahona1982} but has potentially intractable projection.

Similarly, two-player games~\citep{von2007theory} can often be solved as bilinear minimax problems. 
When a strategy space involves a polynomial number of constraints, the equilibria of such games can be solved efficiently~\citep{koller1994}. 
However, in situations such as the Colonel Blotto game or the Matching Duel~\citep{ahmadinejad2016}, the strategy space is intractably large and defined by an exponential number of linear constraints. 
Fortunately, despite this apparent prohibitive structure, some linear minimization oracles such as the blossom algorithm~\citep{edmonds1965} can efficiently optimize over the matching polytopes.

Robust learning is also often cast as a saddle point minimax problem~\citep{KimMagBoy05}. 
Once again, a FW implementation could leverage fast linear oracles while projection methods would be plagued by slower or intractable sub-problems. 
For instance, if the LMO is max-flow, it could have almost linear runtime while the corresponding projection would require cubic runtime quadratic programming~\citep{Kelner2014}.
Finally, note that the popular generative adversarial networks~\citep{goodfellow14GAN} are formulated as a (non-convex) saddle point optimization problem.

\vspace{-2mm}
\paragraph{Related work.}
The standard approaches to solve smooth constrained saddle point problems are projection-type methods (surveyed in~\citet{xiu2003surveyVIP}), with in particular variations of Korpelevich's extragradient method~\citep{korpelevich1976extragradient}, such as~\citep{nesterov2007dualExtrag} which was used to solve the structured prediction problem~\citep{taskar2006dualExtrag} mentioned above.
There is surprisingly little work on FW-type methods for saddle point problems, although they were briefly considered for the more general \emph{variational inequality} problem (VIP):
\vspace{-0.5mm}
\begin{equation} \label{eq:VIP}
\text{find } \,\, \z^* \in \Z \,\, \text{ s.t.} \,\, \innerProd{\FF(\z^*)}{\z - \z^*} \geq 0, \;\; \forall \z \in \Z,
\end{equation}
where $\FF$ is a Lipschitz mapping from $\R^p$ to itself and $\Z \subseteq \R^p$. By using $\Z = \X \times \Y$ and $\FF(\z) = (\nabla_x \L(\z),-\nabla_y \L(\z))$, the VIP~\eqref{eq:VIP} reduces to the equivalent optimality conditions for the saddle point problem~\eqref{eq:convex_concave}. 
\citet{hammond1984solving} showed that a FW algorithm with a step size of $O(1/t)$ converges for the VIP~\eqref{eq:VIP} when the set $\Z$ is strongly convex,
while FW with a generalized line-search on a saddle point problem is sometimes non-convergent when $\Z$ is a polytope (see also~\citep[\S~3.1.1]{patriksson1999nonlinear}). 
She conjectured though that using a step size of $O(1/t)$ was also convergent when $\Z$ is a polytope -- a problem left open up to this point. 
More recently, \citet{juditsky2016VIP} (see also~\citet{cox2015decompositionLMO}) proposed a method to transform a VIP on $\Z$ where one has only access to a LMO, to a ``dual'' VIP on which they can use a projection-type method. 
\citet{lan2013FW} proposes to solve the saddle point problem~\eqref{eq:convex_concave} by running FW on $\X$ on the \emph{smoothed} version of the problem $\max_{\y \in \Y} \L(\x,\y)$, thus requiring a projection oracle on $\Y$. 
In contrast, in this paper we study simple approaches that do not require any transformations of the problem~\eqref{eq:convex_concave} nor any projection oracle on $\X$ or $\Y$. Finally, \citet{he2015semi} introduced an interesting extragradient-type method to solve~\eqref{eq:VIP} by approximating the projections using linear oracles. In contrast to our proposal, their work does not cover the geometric convergence for the strongly convex case.

\vspace{-1.5mm}
\paragraph{Contributions.}
In \S~\ref{sec:algorithms}, we extend several variants of the FW algorithm to solve the saddle point problem~\eqref{eq:convex_concave} that we think could be of interest to the machine learning community. 
In \S~\ref{sec:SP-FW_strong_convexity}, we give a first proof of (geometric) convergence for these methods over polytope domains under the assumptions of sufficient strong convex-concavity of~$\L$, giving a partial answer to the conjecture from~\citet{hammond1984solving}. 
In \S~\ref{sec:strongly_set}, we extend and refine the previous convergence results when $\X$ and $\Y$ are strongly convex sets and the gradient of $\L$ is non-zero over $\X \times \Y$, while we survey the pure bilinear case in \S~\ref{sec:bilinear}.
We finally present illustrative experiments for our theory in \S~\ref{sec:experiments}, noticing that the convergence theory is still incomplete for these methods.

\begingroup
\setlength{\intextsep}{-1.5mm}
\setlength{\floatsep}{-1.5mm}
\setlength{\textfloatsep}{-1.5mm}
\setlength{\intextsep}{-1.5mm}
\begin{figure*}[t]
\begin{minipage}[t]{.45\textwidth}
  \begin{algorithm}[H]
    \caption{Frank-Wolfe algorithm}\label{FW}
    \begin{algorithmic}[1]
      \STATE Let $\x^{(0)} \in \X$
      \FOR{$t=0 \ldots T$} 
      \vspace*{.4mm}
      \STATE Compute $\rt = \nabla f(\x^{(t)})$
      \vspace*{.85mm}
      \STATE Compute $\st := \underset{\s \in \X}{\text{ argmin }} \innerProd{\s}{\rt}$
      \vspace*{.45mm}
      \STATE Compute $g_t := \left\langle \xt-\st, \rt\right\rangle$
      \vspace*{.85mm}
      \STATE \textbf{if} $ \gap \leq \epsilon$ \textbf{then} \textbf{return} $\x^{(t)}$ 
      \vspace*{.85mm}
      \STATE Let $\gamma = \frac{2}{2+t}$ (or do line-search)
      \vspace*{.85mm}
      \STATE Update $\xtt := (1- \gamma) \xt + \gamma \st$
      \ENDFOR
    \end{algorithmic}
  \end{algorithm}
\end{minipage}
\hspace{-2mm}
\begin{minipage}[t]{.55\textwidth}
\begin{algorithm}[H]
  \caption{Saddle point Frank-Wolfe algorithm: \textbf{SP-FW}}\label{alg:SP-FW}
  \begin{algorithmic}[1]
      \STATE Let $\z^{(0)} =(\x^{(0)},\y^{(0)}) \in \M$
      \FOR{$t=0 \ldots T$}
      \STATE Compute $ \rt := \begin{pmatrix} 
                  \nabla_x \L(\xt,\yt) \\
                  -\nabla_y \L(\xt,\yt)\end{pmatrix}$ \label{algLine:rk}
      \STATE Compute $\st := \underset{\z \in \X\times \Y}{\text{ argmin }} \innerProd{\z}{\rt}$
      \STATE Compute $g_t := \left\langle \zt-\st, \rt\right\rangle $
      \STATE \textbf{if} $ g_t\leq \epsilon$ \textbf{then} \textbf{return} $\z^{(t)}$
      \STATE Let $\gamma =  \min\left(1, \frac\CondNumb{2C}\gap\right) $ or $\gamma = \frac{2}{2+t}$ 
      \hfill \raisebox{0pt}[0pt][0pt]{
      		\begin{minipage}{2.5cm}
      			\emph{\small($\CondNumb$ and $C$ set as \\ case~(I)~in~Thm.~\ref{thm:conv})}
      		\end{minipage}
      		}
      \STATE Update $\ztt := (1- \gamma) \zt + \gamma \st$
      \ENDFOR
\end{algorithmic}
\end{algorithm}
\end{minipage}
\vspace{-3mm}
\end{figure*}
\begin{figure*}[ttt]
\begin{minipage}[t]{\textwidth}
\begin{algorithm}[H]
    \caption{Saddle point away-step Frank-Wolfe algorithm: \textbf{SP-AFW}$(\z^{(0)}, \Vertices, \epsilon)$}\label{alg:AFW}
  \begin{algorithmic}[1]
  \STATE Let $\z^{(0)}=(\x^{(0)},\y^{(0)}) \in \Vertices$, $\Coreset_x^{(0)} := \{\x^{(0)}\}$ and $\Coreset_y^{(0)} := \{\y^{(0)}\}$
  \FOR{$t=0\dots T$}
    \STATE Let $\st := \lmo \!\left(\rt\right)$ %
       and $\dt_\FW := \st - \z^{(t)}$ \hspace{27mm}~~ \emph{\small($\rt$ as defined in L\ref{algLine:rk} in Algorithm~\ref{alg:SP-FW})}
    \STATE Let $\vt \in \displaystyle\argmax_{\vv \in \Coreset_x^{(t)} \times \Coreset_y^{(t)} } \textstyle\left\langle \rt, \vv \right\rangle$ and $\dt_\away := \z^{(t)} - \vt$ \hspace{51mm} \emph{\small(the away direction)}  \label{algLine:awayCorner}
    \STATE \textbf{if} $g_t^\FW  := \left\langle -\rt, \dt_\FW\right\rangle  \leq \epsilon$ \textbf{then} \textbf{return} $\z^{(t)}$ \hspace{38mm}\emph{\small(FW gap is small enough, so return)} \label{algLine:gapFW}
      \IF{$\left\langle -\rt, \dt_\FW\right\rangle  \geq \left\langle -\rt, \dt_\away\right\rangle$ } \label{algLine:begAFW}
      \STATE $\dt :=  \dt_\FW$, and $\stepmax := 1$  
           \hspace{77mm}\emph{\small(choose the FW direction)}
      \ELSE
      \STATE $\dt :=  \dt_\away$, and $\stepmax := \min\left\{\tfrac{\alpha_{\vv_x^{(t)}}}{1- \alpha_{\vv_x^{(t)}}},\tfrac{\alpha_{\vv_y^{(t)}}}{1- \alpha_{\vv_y^{(t)}}}\right\}$
        \emph{\small(maximum feasible step size; a \emph{drop step} is when $\stepsize_t = \stepmax$)} \label{algLine:drop_step}
       \vspace{-2mm}	
      \ENDIF \label{algLine:endAFW}
      \STATE Let $\gap^{\PW} = \innerProd{-\rt}{\dt_{\FW} + \dt_{\away}}$ \textbf{and} $\gamma_t = \min \left\{\gamma_{\max} , \frac{\CondNumb^{\PW}}{2 C}\gap^{\PW} \right\}$ \hfill \emph{\small($\CondNumb$ and $C$ set as case (P) in Thm.~\ref{thm:conv})} \label{algLine:gapPFW}
      \STATE Update $\z^{(t+1)} := \z^{(t)} + \stepsize_t \dt$  
        \hspace{4mm}\emph{\small(and accordingly for the weights $\aa^{(t+1)}$, see~\citet{lacoste2015global})} 
      \STATE Update $\Coreset_x^{(t+1)} := \{\vv_x \in \A \,\: \mathrm{ s.t. } \,\: \alpha^{(t+1)}_{\vv_x} > 0\}\;$ and $\;\Coreset_y^{(t+1)}:=\{\vv_y \in \B \,\: \mathrm{ s.t. } \,\: \alpha^{(t+1)}_{\vv_y} >0\}$ \label{algLine:activeSet}
  \ENDFOR
\end{algorithmic}
\end{algorithm}
~\vspace{-1.5mm}
\begin{algorithm}[H]
  \caption{Saddle point pairwise Frank-Wolfe algorithm: \textbf{SP-PFW}$(\z^{(0)}, \Vertices, \epsilon)$}\label{alg:PFW}
  \begin{algorithmic}[1]
  \STATE In Alg.~\ref{alg:AFW}, replace L\ref{algLine:begAFW} to~\ref{algLine:endAFW} by:  $\,\, \dd^{(t)} := \dd_{\PW}^{(t)} \!:= \s^{(t)} - \vv^{(t)}$, and $\stepmax := \min\left\{ \alpha_{\vv_x^{(t)}},\alpha_{\vv_x^{(t)}}\right\}$.
  \end{algorithmic}
\end{algorithm}
\end{minipage}
\vspace{-4mm}
\end{figure*}
\endgroup

\vspace{-0.3mm}
\section{Saddle point Frank-Wolfe (SP-FW)} \label{sec:algorithms}
\vspace{-2.4mm}
\paragraph{The algorithms.} %
\label{par:the_algorithms}
This article will explore three SP extensions of the classical {\em Frank-Wolfe} (FW) algorithm (Alg.~\ref{FW}) which are summarized in Alg.~\ref{alg:SP-FW}, \ref{alg:AFW} and~\ref{alg:PFW}.\footnote{Alg.~\ref{alg:SP-FW} was already proposed by~\citet{hammond1984solving} for VIPs, while our step sizes and Alg.~\ref{alg:AFW} \&~\ref{alg:PFW} are novel.}
We denote by $\zt := (\xt,\yt)$ the iterate computed after $t$ steps.
We first obtain the {\em saddle point FW} (SP-FW) algorithm (Alg.~\ref{alg:SP-FW}) by simultaneously doing a FW update on both convex functions $\L(\cdot,\yt)$ and  $-\L(\xt,\cdot)$ with a properly chosen step size. As in standard FW, the point $\zt$ has a sparse representation as a convex combination of the points previously given by the FW oracle, that is,
\begin{equation}  
\x^{(t)} = \sum_{\vv_x \in \Coreset_x^{(t)}} \alpha_{\vv_x} \vv_x \;\; \text{and}\;\; \y^{(t)} = \sum_{\vv_y \in \Coreset_y^{(t)}} \alpha_{\vv_y} \vv_y.
\end{equation}
These two sets $\Coreset_x^{(t)}$, $\Coreset_y^{(t)}$ of points are called the \emph{active sets}, and we can maintain them separately (thanks to the product structure of $\X \times \Y$) to run the other two FW variants that we describe below (see L\ref{algLine:activeSet} of Alg.~\ref{alg:AFW}).

If we assume that $\X$ and $\Y$ are the convex hulls of two finite sets of points $\A$ and $\B$, we can also extend the {\em away-step Frank-Wolfe} (AFW) algorithm~\citep{Guelat1986AFW,lacoste2015global} to saddle point problems.
As for AFW, this new algorithm can choose an \emph{away} direction~$\dd_\A$ to remove mass from ``bad" atoms in the active set, i.e. to reduce~$\alpha_{\vv}$ for some~$\vv$ (see L\ref{algLine:drop_step} of Alg.~\ref{alg:AFW}),
thereby avoiding the zig-zagging problem that slows down standard FW~\citep{lacoste2015global}. Note that because of the special product structure of the domain, we consider more away directions than proposed in~\citep{lacoste2015global} for AFW (see Appendix~\ref{par:active_set} for more details). Finally, a straightforward saddle point generalization for the {\em pairwise Frank-Wolfe} (PFW) algorithm~\citep{lacoste2015global} is given in Alg.~\ref{alg:PFW}.
The proposed algorithms all preserve several nice properties of previous FW methods (in addition to only requiring LMO's): simplicity of implementation, affine invariance~\citep{jaggi2013revisiting}, gap certificates computed for free, sparse representation of the iterates {and the possibility to have adaptive step sizes using the gap computation.
We next analyze the convergence of these algorithms.}

\paragraph{The suboptimality  error and the gap.} %
\label{par:the_subopti}
To establish convergence, we first define several quantities of interest. 
In classical convex optimization, the suboptimality error $h_t$ is well defined as $h_t := f(\xt)- \min_{\x \in \X} f(\x)$. 
This quantity is clearly non-negative and proving that $h_t$ goes to 0 is enough to establish convergence. 
Unfortunately, in the saddle point setting the quantity $\L(\xt,\yt) - \L^*$ is no longer non-negative and can be equal to zero for an infinite number of points $(\x,\y)$ while $(\x,\y) \notin (\X^*,\Y^*)$. 
For instance, if $\L(\x,\y)=\x \cdot \y$ with $\X = \Y = [-1,1]$, then $\L^* = 0$ and $(\X^*,\Y^*)=\{(0,0)\}$. 
But for all $\x \in \X$ and $\y \in \Y$, $\; \x \cdot 0 = 0 \cdot \y = \L^*$.
The saddle point literature thus considers a non-negative gap function (also known as a merit function \citep{larsson1994class,zhu1998convergence} and \cite[Sec 4.4.1]{patriksson1999nonlinear}) which is zero only for optimal points, in order to quantify progress towards the saddle point. 
We can define the following \emph{suboptimality error} $h_t$ for our saddle point problem:
\vspace{-2mm}
\begin{equation}
\begin{aligned}\label{def:subopt} 
   h_t := \L(\xt,\ytm) - \L(\xtm,\yt), \\
   \text{where} \quad \xtm := \argmin_{\x \in \X} \L(\x,\yt), \\[-1mm]
   \text{and }  \quad \ytm := \argmax_{\y \in \Y} \L(\xt,\y).
\end{aligned}
\vspace{-1mm}
\end{equation}
This is an example of \emph{primal-dual} gap function by noticing that 
\begin{align}
h_t &= \L(\xt,\ytm) -\L^* +\L^* - \L(\xtm,\yt) \notag\\ 
    &= p(\xt) - p(\x^*) + g(\y^*) - g(\yt),  
\end{align}
where $p(\x) := \max_{\y \in \Y} \L(\x,\y)$ is the convex primal function and $g(\y) := \min_{\x \in \X} \L(\x,\y)$ is the concave dual function.
By convex-concavity, $h_t$ can be upper-bounded by the following FW linearization gap 
\citep{jaggi2011sparse,jaggi2013revisiting,larsson1994class,zhu1998convergence}:
\begin{equation}
\begin{aligned}\label{eq:gap}
  \hspace{-2mm} g_t^{\FW}
  & \!:= \left. \max_{\s_x \in \X} \innerProd{\xt- \s_x}{\nabla_x \L(\xt,\yt)} \right\} {\small \!:= g^{(x)}_t}  \\
  & \!+ \left. \max_{\s_y \in \Y} \innerProd{\yt - \s_y}{-\nabla_y \L(\xt,\yt)} \right\} {\small \!:=g^{(y)}_t} \!.
\end{aligned} \hspace{-4mm}
\end{equation}
This gap is easy to compute and gives a stopping criterion since $g_t^{\FW} \geq h_t$.

\paragraph{Compensation phenomenon and difficulty for SP.} %
\label{par:the_compensation_phenomenon}
Even when equipped with a suboptimality error and a gap function (as in the convex case), we still cannot apply the standard FW convergence analysis. 
The usual FW proof sketch uses the fact that the gradient of $f$ is Lipschitz continuous to get
\begin{equation} \label{eq:FWstandard}
  h_{t+1} \leq h_t - \gamma_t g_t^{\FW} + \gamma_t^2 \frac{L\|\dt\|^2}{2}
\end{equation}
which then provides a rate of convergence. 
Roughly, since $g_t \geq h_t$ by convexity, if $\gamma_t$ is small enough then $(h_t)$ will decrease and converge. 
For simplicity, in the main paper, $\|\cdot\|$ will refer to the $\ell_2$ norm of $\R^d$. 
The partial Lipschitz constants and the diameters of the sets are defined with respect to this norm (see~\eqref{def:lipL} in Appendix~\ref{sub:the_lipschitz_constants} for more general norms).

Using the $L$-Lipschitz continuity of $\L$ and letting $\L_t := \L(\xt, \yt)$ as a shorthand, we get
\begin{equation}
\begin{aligned}
  \L_{t+1}
  & \leq \L_t + \gamma_t \innerProd{\dt_x}{\nabla_x \L_t}  + \gamma_t \innerProd{\dt_y}{\nabla_y \L_t}   \\
  & \quad+ \gamma_t^2\frac{L\|\dt\|^2}{2}
\end{aligned}
\end{equation}
where $\dt_x = \st_x - \xt$ and $\dt_y=\st_y-\yt$. 
Then 
\begingroup
\setlength{\thinmuskip}{0mu}
\setlength{\medmuskip}{0mu}
\setlength{\thickmuskip}{0mu}
\begin{equation}\label{eq:oscilation}
  \hspace{-1mm} \L_{t+1} - \L^*
   \hspace{1mm} \leq \hspace{1mm} \L_t - \L^*
    - \gamma_t \left( g^{(x)}_t -  g_t^{(y)}\right) 
    + \gamma_t^2 \frac{L \|\dt\|^2}{2}. \hspace{-1mm}
\end{equation}
\endgroup
Unfortunately, the quantity $\gap^{\FW}$ does \emph{not} appear above and we therefore cannot control the oscillation of the sequence (the quantity $g^{(x)}_t -  g_t^{(y)}$ can make the sequence increase or decrease). 
Instead, we must focus on more specific SP optimization settings and introduce other quantities of interest in order to establish convergence.

\paragraph{The asymmetry of the SP.} %
\label{par:the_asymmetry_of_the_sp}
\citet[p.~165]{hammond1984solving} showed the divergence of the SP-FW algorithm with an extended line-search step-size on some bilinear objectives. She mentioned that the difficulty for SP optimization is contained in this bilinear coupling between $\x$ and $\y$.  
More generally, most of the examples of SP functions cited in the introduction can be written in the form:
\begin{equation}
\L(\x,\y) = f(\x) + \x^\top \!M \y - g(\y), \; \text{$f$ and $g$ convex.}
\end{equation}
In this setting, the bilinear part $M$ is the only term preventing us to apply theorems on standard FW.  
\citet[p.~175]{hammond1984solving} also conjectured that the SP-FW algorithm with $\stepsize_t=\nicefrac{1}{(t+1)}$ performed on a uniformly strongly convex-concave objective function (see~\eqref{def:uniform}) over a polytope should converge. We give a partial answer to this conjecture in the following section.

\section{SP-FW for strongly convex functions}
\label{sec:SP-FW_strong_convexity}
\paragraph{Uniform strong convex-concavity.} %
\label{par:strong_convex_concavity}
In this section, we will assume that $\L$ is uniformly $(\mu_\X,\mu_\Y)$-strongly convex-concave, which means that the following function is convex-concave:
\begin{equation}\label{def:uniform}
 (\x,\y) \mapsto \L(\x,\y) - \frac{\mu_\X}{2} \|\x\|^2 + \frac{\mu_\Y}{2}\|\y\|^2.
\end{equation}
\paragraph{A new merit function.} %
\label{par:the_new_quantity}
To prove our theorem, we use a different quantity $w_t$ which is smaller than $h_t$ but still a valid merit function in the case of \emph{strongly convex-concave} SPs (where $(\x^*,\y^*)$ is thus unique); see~\eqref{eq:relate_h_w} below. 
For $(\x^*,\y^*)$ a solution of~\eqref{eq:convex_concave}, we define the non-negative quantity $w_t$:
\begin{equation}
 w_t
   := \underbrace{\L(\xt,\y^*) -\L^*}_{:=  w_t^{(x)}} + \underbrace{\L^* - \L(\x^*,\yt)}_{:=w_t^{(y)}}.
\end{equation}
Notice that $ w_t^{(x)}$ and $w_t^{(y)}$ are non-negative, and that $w_t \leq h_t$ since:
\[
  \L(\xt,\ytm) - \L(\xtm,\yt) \geq \L(\xt,\y^*) - \L(\x^*,\yt).
\]
In general, $w_t$ can be zero even if we have not reached a solution.
For example, with $\L(\x,\y) = \x \cdot \y$ and $\X = \Y = [-1,1]$, then $\x^* = \y^* = \0$, implying $w_t = 0$ for any $(\x^{(t)},\y^{(t)})$. 
But for a uniformly strongly convex-concave~$\L$, this cannot happen and we can prove that~$w_t$ has the following nice property (akin to $\|\x - \x^*\| \leq \sqrt{\nicefrac{2}{\mu} (f(\x) - f(\x^*))}$ for a $\mu$-strongly convex function~$f$; see Proposition~\ref{prop:relation_primal} in Appendix~\ref{sub:relation_between_the_primal_errors}):
\begin{equation} \label{eq:relate_h_w}
  h_t \leq \sqrt{2} P_\L \sqrt{w_t} \, , \vspace{-1mm} 
\end{equation}
where 
\begin{equation}
  \!\!\! P_\L \leq \sqrt{2} \underset{\z \in \M}{\sup} \left\{\frac{\|\nabla_x \L(\z)\|_{\X^*}}{\sqrt{\mu_\X}} ,\frac{\|\nabla_y \L(\z)\|_{\Y^*}}{\sqrt{\mu_\Y}}   \right\} .
\end{equation}
\paragraph{Pyramidal width and distance to the border.} %
\label{par:pyramidal_width}
We now provide a theorem that establishes convergence in two situations: \eqref{enu:situation1} when the SP belongs to the interior of $\M$; \eqref{enu:situation2} when the set is a polytope, i.e. when there exist two finite sets such that $\X = \conv(\mathcal A)$ and $\Y = \conv(\mathcal B)$).
Our convergence result holds when (roughly) the strong convex-concavity of~$\L$ is big enough in comparison to the cross Lipschitz constants $L_{XY}$, $L_{YX}$ of $\nabla \L$ (defined in~\eqref{eq:cross-lip} below) multiplied by geometric ``condition numbers'' of each set. The condition number of $\X$ (and similarly for $\Y$) is defined as the ratio of its \emph{diameter} $D_\X : = \sup_{\x, \x' \in \X} \|\x -\x'\|$ over the following appropriate notions of ``width'':
\begin{align} 
\hspace{-7mm}\text{\small border distance: }\, \delta_\X &:= \min_{\s \in \partial  \X} \|\x^* -\s\|  \!\!\! &&\text{for \eqref{enu:situation1},} \hspace{-2mm} \\
\hspace{-7mm}\text{\small pyramidal width: }\, \delta_{\mathcal A} &:=\PWidth(\mathcal A)  &&\text{for \eqref{enu:situation2}.} & \label{eq:PwidthL} \hspace{-2mm}
\end{align}
The pyramidal width~\eqref{eq:PwidthL} is formally defined in~Eq.~9 of~\citet{lacoste2015global} and in  Appendix~\ref{sub:affine_invariant_measures_of_strong_convexity}.
Given the above constants, we can state below a non-affine invariant version of our convergence theorem (for simplicity). The affine invariant versions of this theorem are given in Thm.~\ref{thm:conv_affine_invariant} and~\ref{thm:conv_sublin} in Appendix~\ref{subsec:proof_main_thm} (with proofs).

\begin{theorem}\label{thm:conv}
  Let $\L$ be a convex-concave function and $\M$ a convex and compact set. 
  Assume that the gradient of  $\L$ is $L$-Lipschitz continuous, that $\L$ is $(\mu_\X,\mu_\Y)$-strongly convex-concave, and that we are in one of the two following situations:
  \leqnomode %
  \begin{align*}
  &\parbox{7cm}{The SP belongs to the interior of $\M$. In this case, set $\gap = \gap^{\FW}$ (as in L\ref{algLine:gapFW} of Alg.~\ref{alg:AFW}), $\delta_\mu:= \sqrt{\min(\mu_\X\delta_\X^2,\mu_\Y\delta_\Y^2)}$ and $a :=  1 $. ``Algorithm'' then refers to SP-FW.
    }
    \label{enu:situation1} 
    \tag{I}
      \\[2mm]
    &\parbox{7cm}{The sets $\X$ and $\Y$ are polytopes. In this case, set $\gap = \gap^{\PW}$ (as in L\ref{algLine:gapPFW} of Alg.~\ref{alg:AFW}), $\delta_\mu := \sqrt{\min(\mu_\X \delta_{\mathcal A}^2, \mu_\Y \delta_{\mathcal B}^2)}$ and $a:=  \frac{1}{2}$. ``Algorithm'' then refers to SP-AFW. Here $\delta_\mu$ needs to use the Euclidean norm for its defining constants.
    }
    \label{enu:situation2} 
    \tag{P}
  \end{align*}
  \reqnomode
  In both cases, if $\CondNumb:= a -  \tfrac{\sqrt{2}}{\delta_\mu}\max\left\{ \tfrac{D_\X L_{XY}}{\sqrt{\mu_\Y}}, \tfrac{D_\Y L_{YX}}{\sqrt{\mu_\X}}\right\}$ is positive, then the errors $h_t$~\eqref{def:subopt} of the iterates of the algorithm with step size
  $\gamma_t = \min\{\stepmax, \frac\CondNumb{2C}\gap\}$ decrease
  geometrically as 
  \begin{equation*} 
    h_t = O \left( (1- \rho)^{\frac{k(t)}{2}} \right) \: \text{ and }
 \: \min_{s\leq t} g_s^\FW = O \left( (1- \rho)^{\frac{k(t)}{2}} \right)
  \end{equation*} 
  where
  $\rho := \CondNumb^2\frac{\delta_\mu^2}{2C}$, $C := \tfrac{L D_\X^2+ L D_\Y^2}{2}$ and $k(t)$ is the number of non-drop step after $t$ steps (see L\ref{algLine:drop_step} in Alg.~\ref{alg:AFW}). In case~\eqref{enu:situation1} we have $k(t)=t$ and in case~\eqref{enu:situation2} we have $k(t)\geq t/3$.
  For both algorithms, if $\delta_\mu >2\max\left\{ \tfrac{D_\X L_{XY}}{\mu_\X}, \tfrac{D_\Y L_{YX}}{\mu_\Y}\right\}$, we also obtain a sublinear rate with the universal choice $\gamma_t =\min\{\stepmax, \frac{2}{2+k(t)} \}$. This yields the rates:
  \vspace{-3mm}
  \begin{equation}
         \min_{s\leq t} h_s \leq  \min_{s\leq t} g_s^\FW = O\left( \frac{1}{t} \right).
    \end{equation}
\end{theorem}

Clearly, the sublinear rate seems less interesting than the linear one but has the added convenience that the step size can be set without knowledge of various constants that characterize $\L$.  Moreover, it provides a partial answer to the conjecture from~\citet{hammond1984solving}.

\paragraph{Proof sketch.} %
\label{par:proof_sketch}

Strong convexity is an essential assumption in our proof; it allows us to relate~$w_t$ to how close we are to the optimum. 
Actually, by $\mu_\Y$-strong concavity of $\L(\x^*,\cdot)$, we have
\begingroup
\setlength{\thinmuskip}{0.5mu}
\setlength{\medmuskip}{0mu}
\setlength{\thickmuskip}{0.5mu}
\begin{equation}\label{eq:close_y_w}
   \| \yt - \y^*\| \leq \sqrt{\frac{2}{\mu_\Y} \left( \L^*  - \L(\x^*,\yt)\right)} = \sqrt{\frac{2}{\mu_\Y} w_t^{(y)}}.
 \end{equation}
\endgroup
Now, recall that we assumed that $\nabla \L$ is Lipschitz continuous. In the following, we will call $L$ the {\em Lipschitz continuity constant} of $\nabla \L$ and $L_{XY}$ and $L_{YX}$ its (cross) {\em  partial Lipschitz constants}. For all $\x, \, \x' \in \X, \;\y, \, \y'  \in \Y$, these constants satisfy
\begin{equation}
\begin{aligned}\label{eq:cross-lip}
   \!\! \|\nabla_x\L(\x,\y) - \nabla_x\L(\x,\y')\|_{\X^*} \leq L_{XY} \|\y- \y'\|_\Y,\\
   \!\! \|\nabla_y\L(\x,\y) - \nabla_y\L(\x',\y)\|_{\Y^*} \leq L_{YX} \|\x- \x'\|_\X.
\end{aligned}
\end{equation}
Note that $L_{XY},L_{YX}\leq L$ if $\|(\x,\y)\| := \|\x\|_\X + \|\y\|_\Y$.
Then, using Lipschitz continuity of the gradient,
\begin{align}
 \L(\xtt,\y^*)  & \leq  \L(\xt,\y^*)  +  \stepsize \innerProdCompressed{\dt_x}{\nabla_x \L(\xt,\y^*)} \notag \\
                  & \quad+ \stepsize^2 \frac{L \|\dt_x\|^2}{2}.
\end{align}
Furthermore, setting $(\x,\y) = (\xt, \y^*)$ and $\y'= \yt$ in Equation~\eqref{eq:cross-lip}, we have
\begin{equation}
\begin{aligned}\label{eq:intermediary_scheme}
 w_{t+1}^{(x)} & \leq w_t^{(x)} 
                  - \stepsize g_t^{(x)} 
                  + \stepsize D_\X L_{XY} \|\yt - \y^*\| \\
                 & \quad + \stepsize^2 \frac{L D_\X^2}{2} \, .
\end{aligned}
\end{equation}
Finally, combining \eqref{eq:intermediary_scheme} and \eqref{eq:close_y_w}, we get
\begin{equation}
 \begin{aligned}\label{eq:intermediary_scheme_2}
   w_{t+1}^{(x)} & \leq w_t^{(x)} 
                  - \stepsize g_t^{(x)} 
                  + \stepsize D_\X L_{XY}\sqrt{\frac{2}{\mu_\Y}} \sqrt{w_t^{(y)}} \\
                  & \quad + \stepsize^2 \frac{L D_\X^2}{2}.
 \end{aligned}
\end{equation}
A similar argument on $-\L(\x^*, \ytt)$ gives a bound on $w_t^{(y)}$ much like~\eqref{eq:intermediary_scheme_2}. Summing both yields:
\begin{align}
    w_{t+1} &\leq w_t - \stepsize g_t + 2\stepsize \max\left\{\tfrac{D_\X L_{XY}}{\sqrt{\mu_\Y}},\tfrac{D_\Y L_{YX}}{\sqrt{\mu_\X}}\right\} \sqrt{w_t}  \notag \\
            & \quad + \stepsize^2 \frac{L D_\X^2+ L D_\Y^2}{2}.
\end{align}
We now apply recent developments in the convergence theory of FW methods for strongly convex objectives. \citet{lacoste2015global} crucially upper bound the square root of the suboptimality error on a convex function with the FW gap if the optimum is in the interior, or with the PFW gap if the set is a polytope (Lemma~\ref{lemme:lingap} in Appendix~\ref{app:gapInequalities}). We continue our proof sketch for case~\eqref{enu:situation1} only:\footnote{The idea is similar for case~\eqref{enu:situation2}, but with the additional complication of possible drop steps.}
\begin{equation}
\begin{aligned}
 \label{eq:def_delta_x}
  & {2 \mu_\X \delta_\X^2}\left(\L(\xt,\yt)-\L(\x^*,\yt)\right) \leq \left(\gap^{(x)}\right)^2  \\
  &\qquad \text{where} \quad \delta_\X :=  \min_{\s \in \partial  \X} \|\x^* -\s\|.
 \end{aligned}
\end{equation}
  We can also get the respective equation on $\y$ with $\delta_\Y :=  \min_{\y \in \partial  \Y} \|\y^* -\y\|$ and sum it with the previous one \eqref{eq:def_delta_x} to get:
  \vspace{-1mm}
\begin{equation}\label{eq:def_delta}
     \delta_\mu \sqrt{2w_t} \leq \gap
     \;\, \text{where} \;\,
      \delta_\mu := \sqrt{\min(\mu_\X\delta_\X^2,\mu_\Y\delta_\Y^2)}.
\end{equation}
Plugging this last equation into \eqref{eq:intermediary_scheme_2} gives us 
\begingroup
\setlength{\thinmuskip}{2mu}
\setlength{\medmuskip}{2mu}
\setlength{\thickmuskip}{2mu}
\begin{equation}
\begin{aligned}\label{eq:cst}
    w_{t+1} & \leq w_t 
                  - \CondNumb\stepsize g_t 
                  + \stepsize^2 {C}
                  \;\;\;\; \text{where} \;\;\;\;  C := \tfrac{ L D_\X^2+ L D_\Y^2}{2}
                  \\
    &\text{and}  \;\;\;\; \CondNumb  := 1-  \tfrac{\sqrt{2}}{\delta_\mu}\max\left\{ \tfrac{D_\X L_{XY}}{\sqrt{\mu_\Y}}, \tfrac{D_\Y L_{YX}}{\sqrt{\mu_\X}}\right\}.
\end{aligned} 
\end{equation}
 \endgroup
The recurrence~\eqref{eq:cst} is typical in the FW literature. We can re-apply standard techniques on the sequence~$w_t$ to get a sublinear rate with $\stepsize_t = \frac{2}{2+t}$, or a linear rate with $\stepsize_t= \min\left\{ \stepmax, \frac{\CondNumb g_t}{2C} \right\}$ (which minimizes the RHS of~\eqref{eq:cst} and actually guarantees that $w_t$ will be \emph{decreasing}). Finally, thanks to strong convexity, a rate on $w_t$ gives us a rate on $h_t$ (by~\eqref{eq:relate_h_w}). \qed

\section{SP-FW with strongly convex sets}
\label{sec:strongly_set}
  \paragraph{Strongly convex set.} %
  \label{par:strongly_convex_set}
  One can (roughly) define strongly convex sets as sublevel sets of strongly convex functions~\citep[Prop.~4.14]{vial1983strong}. In this section, we replace the strong convex-concavity assumption on~$\L$ with the assumption that $\X$ and $\Y$ are $\beta$-strongly convex \emph{sets}.
  \begin{definition}[\citet{vial1983strong,polyak1966existence}]\label{def:strong_set}
    A convex set $\X$ is said to be $\beta$-strongly convex with respect to~$\|.\|$ 
    if for any $\x,\y \in \X$ and any $\gamma \in [0,1]$, $B_\beta(\gamma,\x,\y) \subset \X$ where $B_\beta(\gamma,\x,\y)$ is the $\|.\|$-ball of radius ${\gamma(1 - \gamma) \frac{\beta}{2}\|\x-\y\|^2}$ centered at $\gamma \x + (1- \gamma)\y$.
  \end{definition}

  Frank-Wolfe for convex optimization over strongly convex sets has been studied by
  \citet{levitin1966constrained,demyanov1970approximate} and \citet{dunn1979rates}, amongst others. 
  They all obtained a linear rate for the FW algorithm if the norm of the gradient is lower bounded by a constant. 
  More recently, \citet{garber2014faster} proved a sublinear rate $O(1/t^2)$ by replacing the lower bound on the gradient by a strong convexity assumption on the function.
  In the VIP setting \eqref{eq:VIP}, the linear convergence has been proved if the optimization is done under a strongly convex set but this assumption does \emph{not} extend to $\X \times \Y$ which \emph{cannot} be strongly convex if $\X$ or $\Y$ is not reduced to a single element. 
  In order to prove the convergence, we first prove the Lipschitz continuity of the \emph{FW-corner} function $\s(\cdot)$ defined below. 
  A proof of this theorem is given in Appendix \ref{sec:strong_conv_proof}.
  \begin{theorem}\label{thm:s_lip}
   Let $\X$ and $\Y$ be $\beta$-strongly convex sets. 
   If $\min(\|\nabla_{\!x} L(\z)\|_{\X^*}, \|\nabla_{\!y} L(\z)\|_{\Y^*}) \geq\delta>0$ for all $\z \in \M$, then the oracle function $\z \mapsto \s(\z) := \arg\min_{\s \in \M}\innerProd{\s}{\FF(\z)}$ is well defined and is $\frac{4L}{\delta \beta }$-Lipschitz continuous (using the norm $\|(\x,\y)\|_{\X \times \Y} := \|\x\|_\X + \|\y\|_\Y$), where $\FF(\z) := \left(
     \nabla_x \L(\z),
     -\nabla_y \L(\z)
     \right)$.
  \end{theorem}

  \paragraph{Convergence rate.} %
  \label{par:convergence_rate}
  When the FW-corner function~$\s(\cdot)$ is Lipschitz continuous (by Theorem~\ref{thm:s_lip}), we can actually show that the FW gap is decreasing in the FW direction and get a similar inequality as the standard FW one~\eqref{eq:FWstandard}, but, in this case, on the \emph{gaps}: $g_{t+1} \leq \gap (1- \gamma_t) + \gamma_t^2 \|\st - \zt \|^2 C_\delta$. 
  Moreover, one can show that the FW gap on a strongly convex set~$\X$ can be lower-bounded by $\|\s_x^{(t)} - \x^{(t)} \|^2$ (Lemma~\ref{lemma:lower_gap} in Appendix~\ref{sec:strong_conv_proof}), by using the fact that~$\X$ contains a ball of sufficient radius around the midpoint between~$\s_x^{(t)}$ and~$ \x^{(t)}$. 
  From these two facts, we can prove the following linear rate of convergence (\emph{not} requiring any \emph{strong} convex-concavity of~$\L$). 
    \begin{theorem}\label{thm:conv_strong}
     Let $\L$ be a convex-concave function and $\X$ and $\Y$ two compact $\beta$-strongly convex sets. 
     Assume that the gradient of  $\L$ is $L$-Lipschitz continuous and that there exists $\delta>0$ such that $\min(\|\nabla_{\!x} L(\z)\|_*, \|\nabla_{\!y} L(\z)\|_*) \geq \delta \;\, \forall \z \in \M$. Set $C_\delta := 2L + \frac{8L^2}{\beta \delta}$. 
     Then the gap $\gap^{\FW}$~\eqref{eq:gap} of the SP-FW algorithm with step size $\gamma_t = \tfrac{\gap^{\FW}}{\|\st-\zt\|^2 C_\delta}$
     converges linearly as $\gap^{\FW} \leq g_0 \left( 1- \rho \right)^{t}$,
    where $\rho :=  \tfrac{\beta \delta}{16 C_\delta}$.
    \end{theorem}

\section{SP-FW in the bilinear setting}\label{sec:bilinear}

\paragraph{Fictitious play.} %
\label{par:the_fictitious_play_}
  In her thesis, \citet[\S~4.3.1]{hammond1984solving} pointed out that for the bilinear setting:
\begin{equation}\label{eq:pb_bilin}
  \min_{\x \in \Delta_p} \max_{\y \in \Delta_q} \x^\top M \y 
\end{equation}
where $\Delta_p$ is the probability simplex on $p$ elements,
the SP-FW algorithm with step size $\gamma_t = 1/ \left( 1+t \right)$ is equivalent to the fictitious play (FP) algorithm introduced by \citet{brown1951iterative}. The FP algorithm has been widely studied in the game literature. Its convergence has been proved by \citet{robinson1951iterative}, while \citet{shapiro1958note} showed that one can deduce from Robinson's proof a $O(t^{-1/(p+q-2)})$ rate.
Around the same time, \citet{karlin1960mathematical} conjectured that the FP algorithm converged at the better rate of $O(t^{-1/2})$, though this conjecture is still open and Shapiro's rate is the only one we are aware of.
Interestingly, \citet{daskalakis2014counter} recently showed that Shapiro's rate is also a lower bound if the tie breaking rule gets the worst pick an infinite number of times.
Nevertheless, this kind of adversarial tie breaking rule does not seems realistic since this rule is a priori defined by the programmer.
In practical cases (by setting a fixed prior order for ties or picking randomly for example), Karlin's Conjecture~\citep{karlin1960mathematical} is still open. Moreover, we always observed an empirical rate of at least $O(t^{-1/2})$ during our experiments, we thus believe the conjecture to be true for realistic tie breaking rules.

 \paragraph{Rate for SP-FW.} %
 \label{par:Rate for SP-FW.}
 Via the affine invariance of the FW algorithm and the fact that every polytope with $p$ vertices is the affine transformation of a probability simplex of dimension $p$, any rate for the fictitious play algorithm implies a rate for SP-FW.
\begin{corollary}\label{cor:Bilin}
For polytopes $\X$ and $\Y$ with $p$ and $q$ vertices respectively and $\L(\x,\y) = \x^\top M \y$, the SP-FW algorithm with step size $\gamma_t = \frac{1}{t+1}$ converges at the rate $ h_t = O \left(t^{-\frac{1}{p+q-2}} \right).$
\end{corollary}
 This (very slow) convergence rate is mainly of theoretical interest, providing a safety check that the algorithm actually converges.
 Moreover, if Karlin's strong conjecture is true, we can get a $O (1 /\sqrt{t})$ worst case rate which is confirmed by our experiments.

\section{Experiments}\label{sec:experiments}

\begin{figure*}
    \centering
    \begin{subfigure}[t]{0.32\linewidth}
        \centering
        \includegraphics[height = .8 \linewidth, width = 1 \linewidth]{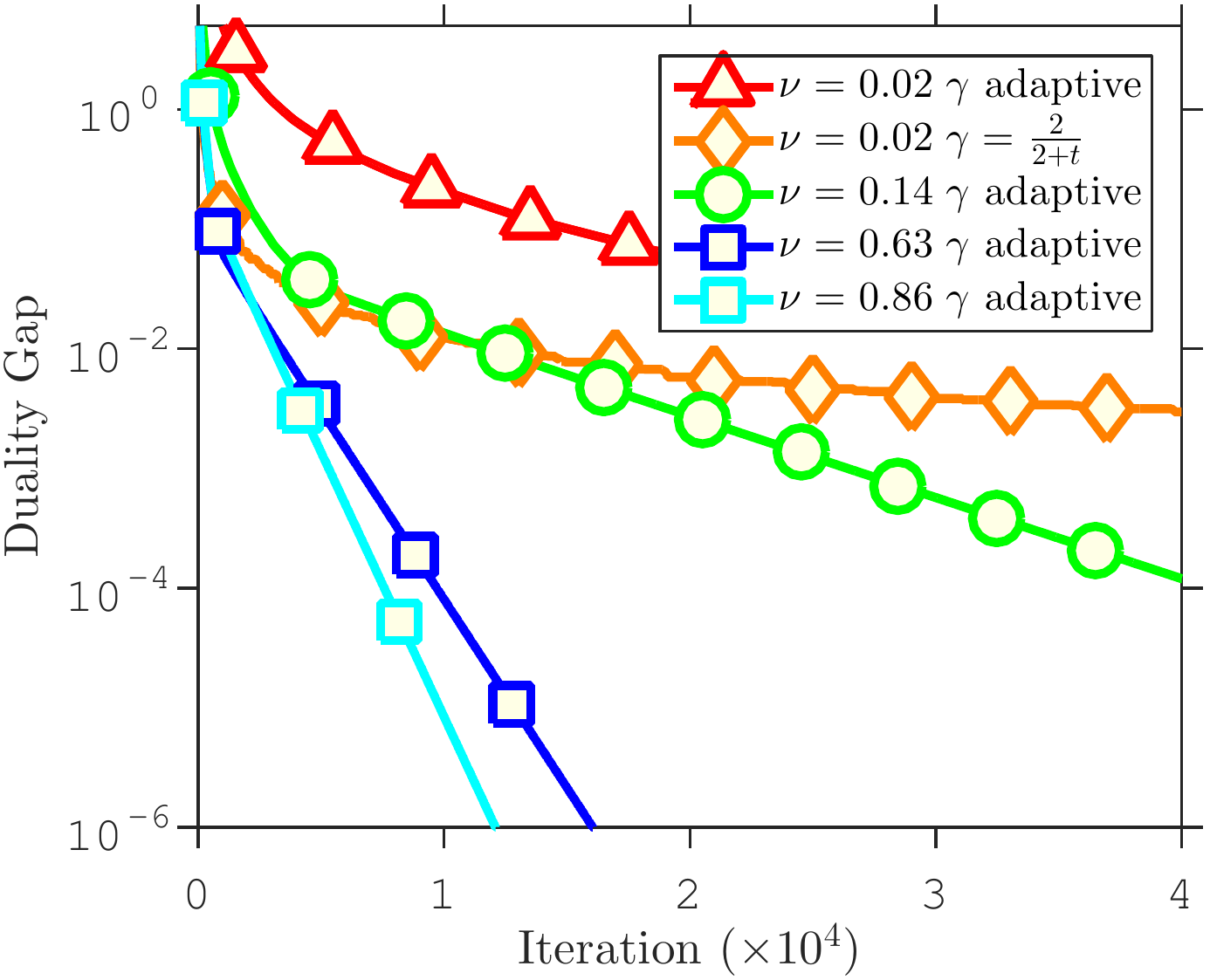}
        \caption{SP in the interior, $d=30$}
        \label{fig:conv_int}
    \end{subfigure}
    \;
    \begin{subfigure}[t]{0.32\linewidth}
        \centering
        \includegraphics[width = 1 \linewidth]{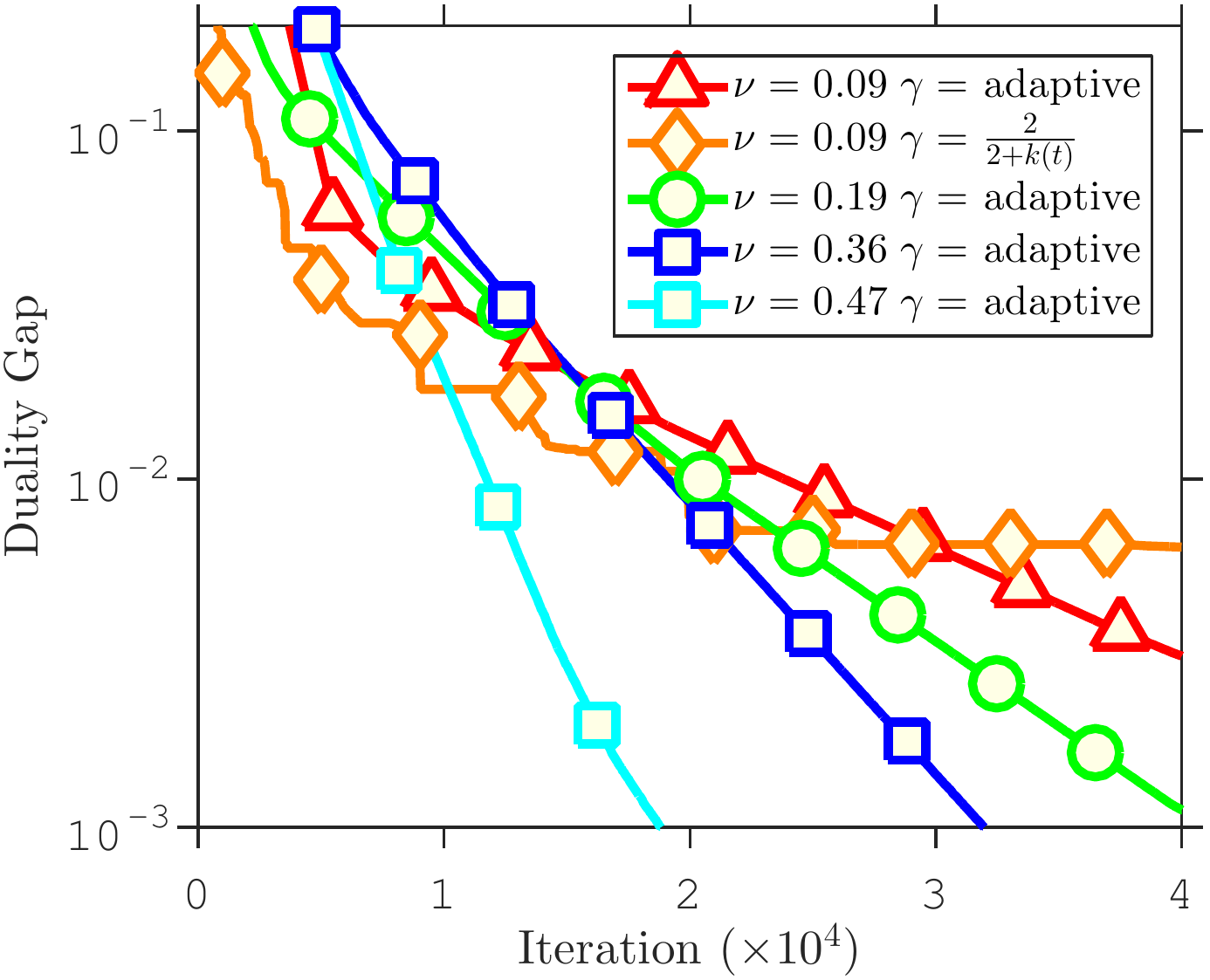}
        \caption{$\M$ is a polytope, $d=30$}
        \label{fig:conv_geom}
    \end{subfigure}
    \;
    \begin{subfigure}[t]{0.32\linewidth}
        \centering
        \includegraphics[height = .8 \linewidth,width = 1 \linewidth]{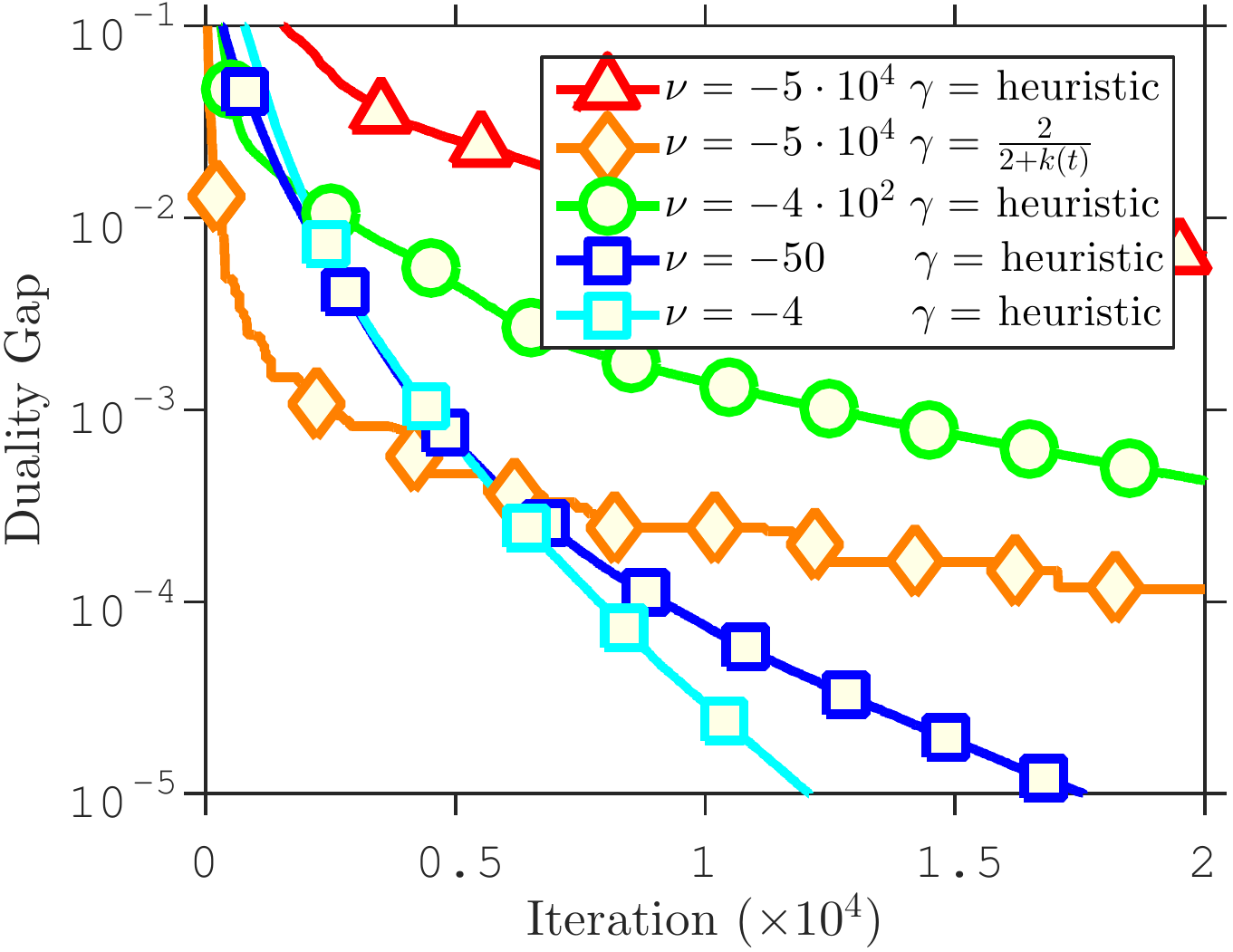}
        \caption{$\M$ polytope, $d=30$, $\CondNumb<0$}
        \label{fig:better_step_size}
    \end{subfigure}
    \\
    \begin{subfigure}[t]{0.32\linewidth}
        \centering
        \includegraphics[width = 1 \linewidth]{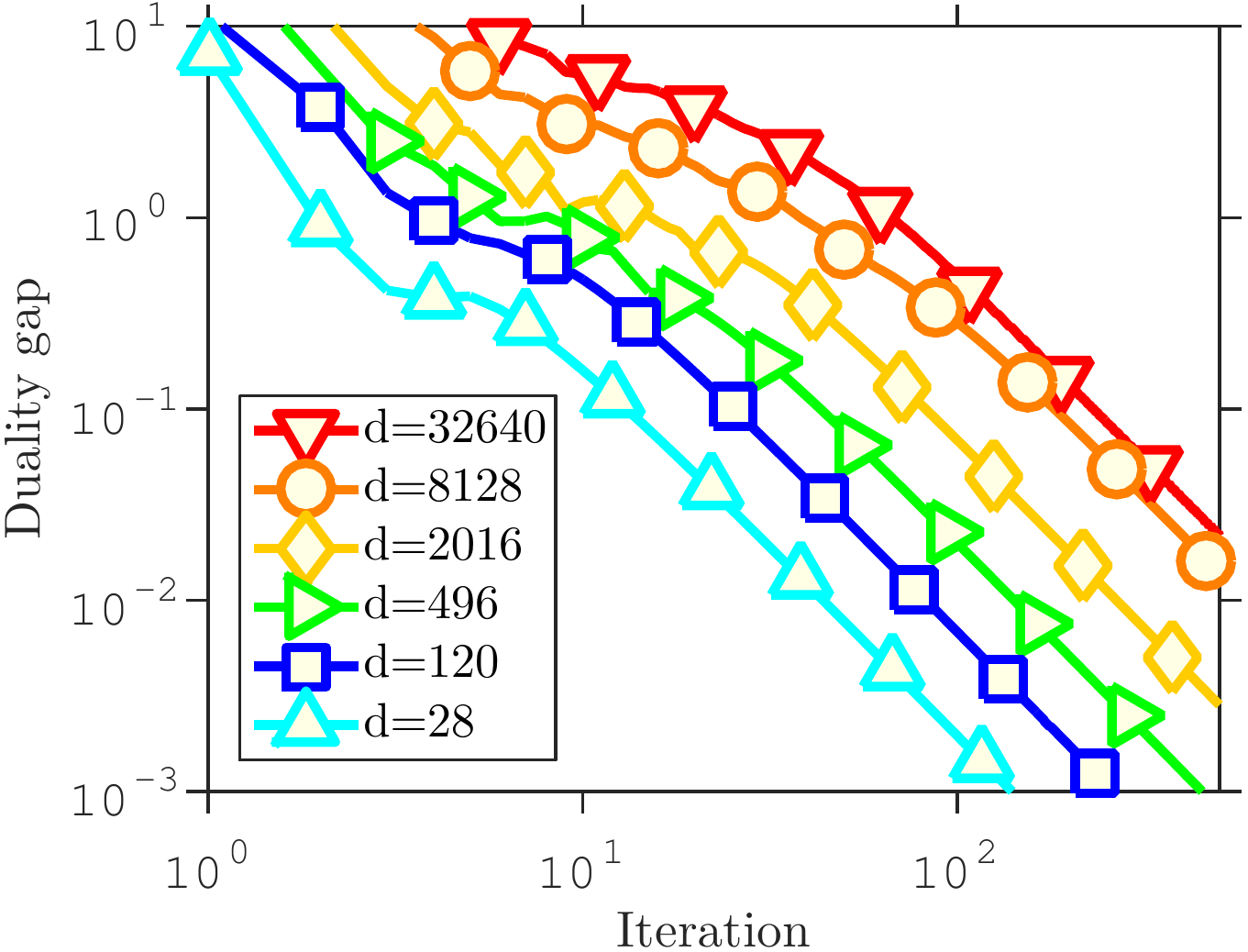}
        \caption{Graphical games}
        \label{fig:graphical}
    \end{subfigure} %
    \;
    \begin{subfigure}[t]{0.32\linewidth}
        \centering
        \includegraphics[height = .75 \linewidth,width = 1 \linewidth]{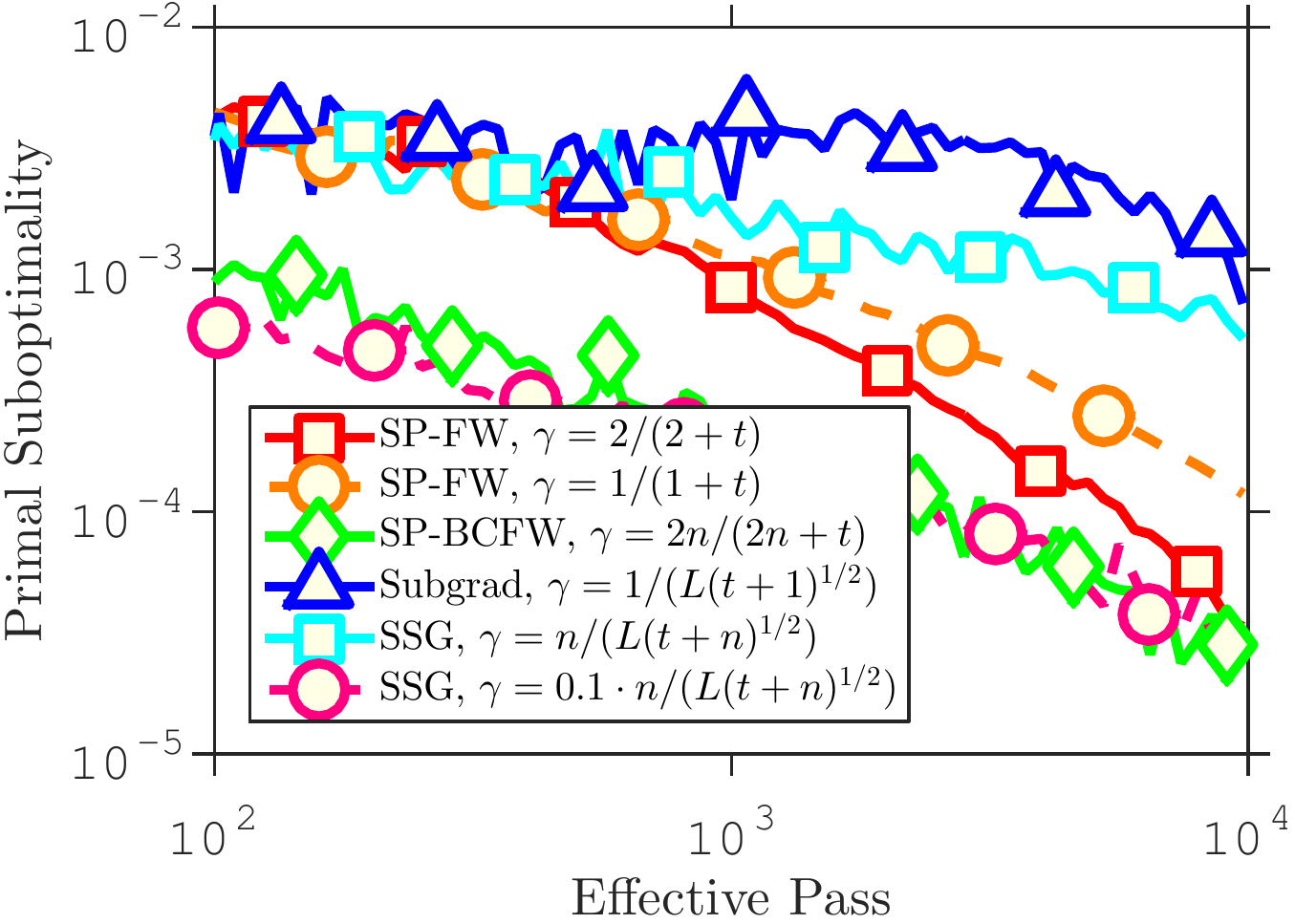}
        \caption{OCR dataset, $R=0.01$.}
        \label{fig:beta_small}
    \end{subfigure}
    \;
    \begin{subfigure}[t]{0.32\linewidth}
        \centering
        \includegraphics[height = .75 \linewidth, width = 1 \linewidth]{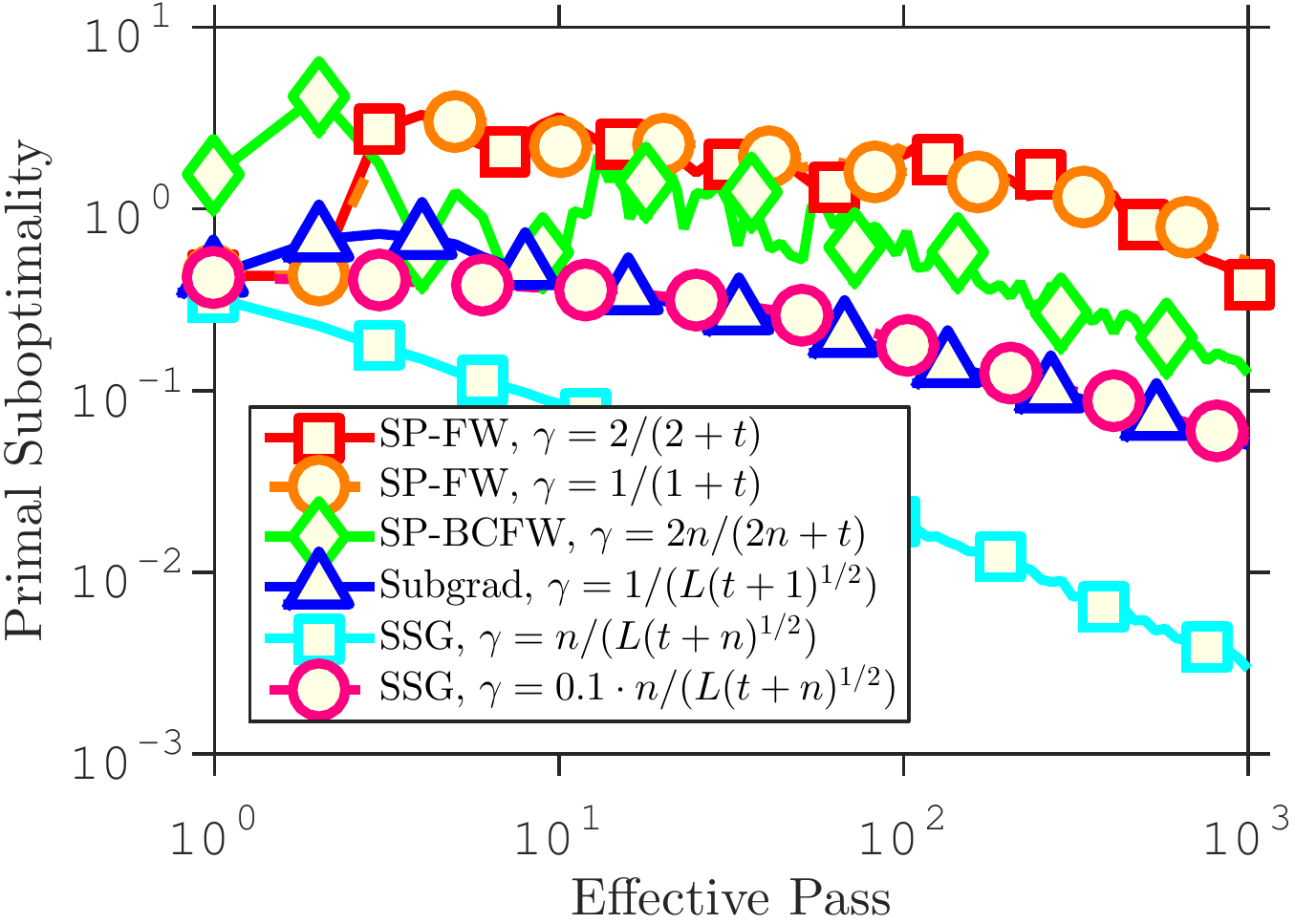}
        \caption{OCR dataset, $R=5$.}
        \label{fig:beta_big}
    \end{subfigure}
    \caption{\small On Figures~\ref{fig:conv_int}, \ref{fig:conv_geom} and~\ref{fig:better_step_size}, we plot on a semilog scale the best gap observed $\min_{s\leq t} g^\FW_s$ as a function of $t$. For experiments~\ref{fig:graphical}, \ref{fig:beta_small} and~\ref{fig:beta_big}, the objective function is bilinear and the convergence is sublinear. An effective pass is one iteration for SP-FW or the subgradient method and $n$ iterations for SP-BCFW or SSG.
    We give more details about these experiments in Appendix~\ref{sec:details_on_the_experiments}.
    }
    \vspace{-2mm}
\end{figure*}
  \paragraph{Toy experiments.} %
  \label{par:experiments_on_toy_examples}
First, we test the empirical convergence of our algorithms on a simple saddle point problem over the unit cube in dimension~$d$ (whose pyramidal width has the explicit value $1/\sqrt{d}$ by Lemma~4 from \citet{lacoste2015global}). 
Thus $\X = \Y := [0,1]^d$ and the linear minimization oracle is simply $\text{LMO}(\cdot)\! = \!-0.5 \cdot(\text{sign}(\cdot)-\bm{1})$. 
We consider the following objective function:
  \begin{equation}
    \frac{\mu}{2}\|\x- \x^*\|_2^2 + (\x-\x^*)^\top M(\y-\y^*) - \frac{\mu}{2} \|\y-\y^*\|_2^2
  \end{equation}
for which we can control the location of the saddle point $(\x^*,\y^*) \in \M$. We generate a matrix $M$ randomly as $M \sim \mathcal U([-0.1,0.1]^{d\times d})$ and keep it fixed for all experiments. For the interior point setup~\eqref{enu:situation1}, we set~$(\x^*,\y^*) \sim \mathcal U ([0.25,0.75]^{2d})$, while we set~$\x^*$ and $\y^*$ to some fixed random vertex of the unit cube for the setup~\eqref{enu:situation2}. With all these parameters fixed, the constant $\CondNumb$ is a function of $\mu$ only. We thus vary the strong convexity parameter $\mu$ to test various $\CondNumb$'s.

We verify the linear convergence expected for the SP-FW algorithm for case~\eqref{enu:situation1} in Figure~\ref{fig:conv_int}, and for the 
SP-AFW algorithm for case~\eqref{enu:situation2} in Figure~\ref{fig:conv_geom}. 
As the adaptive step size (and rate) depends linearly on~$\CondNumb$, the linear rate becomes quite slow for small $\CondNumb$. In this regime (in red), the step size $2/(2+k(t))$ (in orange) can actually perform better, despite its theoretical sublinear rate.

Finally, figure \ref{fig:better_step_size} shows that we can observe a linear convergence of SP-AFW even if $\CondNumb$ is negative by using a different step size. In this case, we use the heuristic adaptive step size $\gamma_t := \gap / \tilde{C}$ where 
$\tilde{C} := L D_\X^2+LD_\Y^2+ L_{XY} L_{YX}\left( D_\X^2/\mu_\X + D_\Y^2/\mu_\Y \right) $. Here $\tilde{C}$ takes into account the coupling between the concave and the convex variable and is motivated from a different proof of convergence that we were not able to complete. The empirical linear convergence in this case is not yet supported by a complete analysis, highlighting the need for more sophisticated arguments.

\paragraph{Graphical games.} We now consider a bilinear objective $\L(\x,\y) = \x^\top M\y$ where exact projections on the sets is intractable, but we have a tractable LMO. The problem is motivated from the following setup. We consider a game between 
two universities ($A$ and $B$) that are admitting $s$ students and have to assign pairs of students into dorms.
If students are unhappy with their dorm assignments, they will go to the other university. The game has a payoff matrix $M$ belonging to $\R^{(s(s-1)/2)^2}$ where $M_{ij,kl}$ is the expected tuition that $B$ gets (or $A$ gives up) if $A$ pairs student $i$ with $j$ and $B$ pairs student $k$ with $l$.
Here the actions $\x$ and $\y$ are both in the marginal polytope of all perfect unipartite matchings. Assume that we are given a graph $G=(V,E)$ with vertices $V$ and edges $E$. For a subset of nodes $S \subseteq V$, let the induced subgraph $G(S)=(S,E(S))$. \citet{edmonds1965} showed that any subgraph forming a triangle can contain at most one edge of any perfect matching. This forms an exponential set of linear equalities which define the matching polytope ${\cal P}(G) \subset \mathbb{R}^E $ as
\begingroup
\setlength{\thinmuskip}{1.3mu}
\setlength{\medmuskip}{1.3mu}
\setlength{\thickmuskip}{1.3mu}
\begin{equation} 
\{ \x \,| \, \x_e \geq 0, \hspace{-2.1mm} %
\sum_{e \in E(S)}\hspace{-1.7mm} \x_e \leq k, \, \forall S \subseteq V,\, { |S|=2 k+1 },\forall e \in E \}.
\end{equation}
\endgroup
While this strategy space seems daunting, the LMO can be solved in ${\cal O}(s^3)$ time using the blossom algorithm~\citep{edmonds1965}. We run the SP-FW algorithm with $\stepsize_t = \nicefrac{2}{(t+2)}$ on this problem with $s = 2^j$ students for $j=3, \ldots, 8$ with results given in Figure~\ref{fig:graphical} ($d=s(s-1)/2$ in the legend represents the dimensionality of the $\x$ and $\y$ variables). The order of the complexity of the LMO is then $O(d^{3/2})$.
In Figure~\ref{fig:graphical}, the observed empirical rate of the SP-FW algorithm (using $\stepsize_t = \nicefrac{2}{(t+2)}$) is $O(1/t^2)$. Empirically, faster rates seem to arise if the solution is at a corner (a pure equilibrium, to be expected for random payoff matrices in light of~\citep{barany2007nash}).

 \paragraph{Sparse structured SVM.} %
  \label{par:structured_svm}
  We finally consider a challenging optimization problem arising from structured prediction. We consider the saddle point formulation~\citep{taskar2006dualExtrag} for a $\ell_1$-regularized structured SVM objective that minimizes the primal cost function
$  p(\bm{w}) :=
    \frac{1}{n} \sum_{i=1}^n \tilde H_i(\bm{w}) 
$,
where $\tilde H_i(\bm{w})= \max_{\y \in \mathcal Y_i} L_i(\y) - \innerProd{\bm{w}}{\bm{\psi}_i(\y)}$ is the structured hinge loss (using the notation from~\citet{lacoste2013block}). We only assume access to the linear oracle computing $\tilde H_i(\bm{w})$. Let $M_i$ have $\big(\bm{\psi}_i(\y)\big)_{\y \in \mathcal Y_i}$ as columns.
  We can rewrite the minimization problem as a bilinear saddle point problem:
  \vspace{-2mm}
  \begin{equation}
  \begin{aligned}
   & \qquad \min_{\|\bm{w}\|_1  \leq R} \frac{1}{n} \sum_i \Big( \max_{\y_i \in \mathcal Y_i} \; \bm{L}_i^\top \y_i - \bm{w}^\top M_i \y_i \Big)  \\
   & = \min_{\|\bm{w}\|_1 \leq R} \frac{1}{n} \sum_i \Big( \max_{\bm{\alpha}_i \in \Delta(|\mathcal Y_i|)} \!\!\bm{L}_i^\top \bm{\alpha}_i - \bm{w}^\top M_i \bm{\alpha}_i \Big). \label{eq:structuredSVM}
  \end{aligned}
  \end{equation}
  Projecting onto $\Delta(|\Y_i|)$ is normally intractable as the size of $|\Y_i|$ is exponential, but the linear oracle is tractable by assumption.
  We performed experiments with 100 examples from the OCR
  dataset ($d_\omega \!= \!4028$)~\citep{taskarmax}. We encoded the structure $\Y_i$ of the $i^{th}$ word with a Markov model: its $k^{th}$ character $\Y_i^{(k)}$ only depends on $\Y_i^{k-1}$ and $\Y_i^{k+1}$. In this case, the oracle function is simply the Viterbi algorithm~\citet{viterbi1967error}. 
 The average length of a word is approximately 8, hence the dimension of $\Y_i$ is $d_{\Y_i} \approx 26^2 \cdot 8 = 5408$ leading to a large dimension for $\Y$, $d_\Y := \sum_{i=1}^n d_{\Y_i} \approx  5 \cdot 10^5$.  
 We run the SP-FW algorithm with step size $\gamma_t = 1/(1+t)$ for which we have a convergence proof (Corollary~\ref{cor:Bilin}), and with $\gamma_t = 2/(2+t)$, which normally gives better results for FW optimization. We compare with the projected subgradient method (projecting on the $\ell_1$-ball is tractable here) with step size $O(1/\sqrt{t})$ (the subgradient of $\tilde H_i(\bm{w})$ is $-\bm{\psi}_i(\y_i^*)$). Following~\citet{lacoste2013block}, we also implement a block-coordinate (SP-BCFW) version of SP-FW and compare it with the stochastic projected subgradient method (SSG). As some of the algorithms only work on the primal and to make our result comparable to~\citet{lacoste2013block}, we choose to plot the primal suboptimality error $p(\bm{w}_t)- p^*$ for the different algorithms in Figure~\ref{fig:beta_small} and~\ref{fig:beta_big} (the $\bm{\alpha}_t$ iterates for the SP approaches are thus ignored in this error).
The performance of SP-BCFW is similar to SSG when we regularize the learning problem heavily (Figure~\ref{fig:beta_small}). 
However, under lower regularization (Figure~\ref{fig:beta_big}), SSG (with the correct step size scaling) is faster. This is consistent with the fact that $\bm{\alpha_t} \neq \bm{\alpha}^*$ implies larger errors on the primal suboptimality for the SP methods, but we note that an advantage of the SP-FW approach is that the scale of the step size is automatically chosen.

\paragraph{Conclusion.}\label{sec:conclusion} 
We proposed FW-style algorithms for saddle-point optimization with the same attractive properties as FW, in particular only requiring access to a LMO. We gave the first convergence result for a FW-style algorithm towards a saddle point over polytopes by building on the recent developments on the linear convergence analysis of AFW.
However, our experiments let us believe that the condition $\CondNumb > 0$ is not required for the convergence of FW-style algorithms.
We thus conjecture that a refined analysis could yield a linear rate for the general uniformly strongly convex-concave functions in both cases~\eqref{enu:situation1} and~\eqref{enu:situation2}, paving the way for further theoretical work.

\clearpage
\subsubsection*{Acknowledgments} %
\label{par:acknowledgments} 
Thanks to N. Ruozzi and
A. Benchaouine for helpful discussions. Work supported in
part by DARPA N66001-15-2-4026, N66001-15-C-4032 and NSF
III-1526914, IIS-1451500, CCF-1302269.
  \bibliographystyle{abbrvnat} 
  {
  \bibliography{bib}
  }

  \clearpage

  \appendix
  \onecolumn
  \fontsize{11}{13}
  \selectfont

  {\huge  Appendix}
  \paragraph{Outline.} %
  \label{par:outline}
  Appendix~\ref{sec:away_step_frank_wolfe} provides more details about the saddle point away-step Frank-Wolfe (SP-AFW) algorithm. Appendix~\ref{sec:affine_invariant} is about the affine invariant formulation of our algorithms, therein, we introduce some affine invariant constants and prove relevant bounds. 
  Appendix~\ref{sec:relations_gaps} presents some relationships between the primal suboptimalities and dual gaps useful for the convergence proof.
  Appendix~\ref{appendix:analysis} gives the affine invariant convergence proofs of SP-FW and SP-AFW in the strongly convex function setting introduced in Section~\ref{sec:SP-FW_strong_convexity}. 
  Appendix~\ref{sec:strong_conv_proof} gives the proof of linear convergence of SP-FW in the strongly convex set setting as defined in Section~\ref{sec:strongly_set}.
  Finally, Appendix~\ref{sec:details_on_the_experiments} provides details on the experiments.

  \section{Saddle point away-step Frank-Wolfe (SP-AFW)} %
  \label{sec:away_step_frank_wolfe}
  In this section, we describe our algorithms SP-AFW and SP-PFW with a main focus on how the away direction is chosen. We also rigorously define a \emph{drop step} and prove an upper bound on their number. In this section, we will assume that there exist two finites sets $\A$ and $\B$ such that 
  $\X = \conv(\A)$ and $\Y = \conv(\B)$.
  \paragraph{Active sets and away directions.} %
  \label{par:active_set}
   Our definition of \emph{active set} is an extension of the one provided in \citet{lacoste2015global}, we follow closely their notation and their results. Assume that we have the current expansion,
  \begin{equation}
    \xt = \sum_{\vv_x \in \mathcal S_x^{(t)}} \alpha_{\vv_x}^{(t)} \vv_x
    \quad \text{where} \quad
     \mathcal S_x^{(t)} := \left\{ \vv_x \in \A\;;\; \alpha_{\vv_x}^{(t)} >0 \right\},
  \end{equation}
  and a similar one for $\yt$.
   Then, the current iterate has a sparse representation as a convex combination of all possible pairs of atoms belonging to $S_x^{(t)}$ and $S_y^{(t)}$, i.e. 
   \begin{equation} \label{eq:ActiveSetImplicit}
     \zt = \sum_{\vv \in \mathcal S^{(t)}} \alpha_{\vv}^{(t)} \vv
     \quad \text{where} \quad
     \mathcal S^{(t)} := \left\{ \vv \in \Vertices\;;\; \alpha_{\vv}^{(t)}:= \alpha_{\vv_x}^{(t)}\alpha_{\vv_y}^{(t)} > 0\right\}.
   \end{equation}
  The set $\mathcal S^{(t)}$ is the current (implicit) \emph{active set} arising from the special product structure of the domain and it defines potentially more away directions than proposed in~\citep{lacoste2015global} for AFW, and these are the directions that we use in SP-AFW and SP-PFW.
  Namely, for every corner $\vv=(\vv_x,\vv_y)$ and $\vv' = (\vv'_x,\vv'_y)$ already picked, $\x-\vv_x$ is a feasible directions in $\X$ and $\y-\vv_y'$ is a feasible direction in $\Y$.
  Thus the combination $(\x-\vv_x,\y-\vv'_y)$ is a feasible direction even if the particular corners $\vv_x$ and $\vv_y'$ have never been picked together. We thus maintain the iterates on $\X$ and $\Y$ as independent convex combination of their respective active sets of corners (Line~\ref{algLine:activeSet} of Algorithm~\ref{alg:AFW}).
  
  Note that after $t$ iteration, the current iterate $\zt$ is $t$-sparse whereas the size of the active set $\mathcal S^{(t)}$ defined in~\eqref{eq:ActiveSetImplicit} can be of size~$t^2$. Nevertheless, because of the block formulation of the away oracle (Line~\ref{algLine:awayCorner} in Algorithm~\ref{alg:AFW}), the away direction can be found in $O(t)$ since we only need to use $S_x^{(t)}$ and $S_y^{(t)}$ separately to compute the away direction in $\mathcal S^{(t)}$.  %
  Moreover, we only need to track at most $t$ corners in $\A$ and $t$ ones in $\B$ to get this bigger active set.
  We can now define the maximal step size for an away direction.
  \paragraph{Maximal step size.} %
  \label{par:maximal_step size}
 For the standard AFW algorithm,  
 \citet{lacoste2015global} suggest to use the maximum step size $\stepmax = \alpha_{\vv^{(t)}}/(1- \alpha_{\vv^{(t)}})$ when using the away direction $\z^{(t)}-\vv^{(t)}$, to guarantee that the next iterate stays feasible. Because we have a product structure of two blocks, we actually consider more possible away directions by maintaining a separate convex combination on each block in our Algorithm~\ref{alg:AFW} (SP-AFW) and~\ref{alg:PFW} (SP-PFW). More precisely, suppose that we have $\xt = \sum_{\vv_x \in S^{(t)}_x} \alpha^{(t)}_{\vv_x} \vv_x$ and $\yt = \sum_{\vv_y \in S^{(t)}_y} \alpha^{(t)}_{\vv_y} \vv_y$, then the following maximum step size $\gamma_{\max}$ (for AFW) ensures that the iterate $\ztt$ stays feasible:
   \begin{equation} \label{eq:update}
    \ztt := \zt + \gamma_t  \dt_A 
    \quad \text{with} \quad \gamma_t \in [0,\gamma_{\max}]
    \quad \text{and} \quad
     \gamma_{\max} := \min\left\{\frac{\alpha^{(t)}_{\vt_x}}{1 -\alpha^{(t)}_{\vt_x}},\frac{\alpha^{(t)}_{\vt_y}}{1 -\alpha^{(t)}_{\vt_y}} \right\}.
    \end{equation}
    A larger $\stepsize_t$ makes one of the coefficients in the convex combination for the iterate negative, thus no more guaranteeing that the iterate stays feasible. A similar argument can be used to derive the maximal step size for the PFW direction in Algorithm~\ref{alg:PFW}.
  \paragraph{Drop steps.} %
  \label{par:drop_steps}
     A \emph{drop step} is when $\gamma_t=\stepmax$ for the away-step update~\eqref{eq:update} \citep{lacoste2015global}. In this case, at least one corner is removed from the active set. We show later in Lemma~\ref{lemme:Lin} that we can still guarantee progress for this step, i.e. $w_{t+1} < w_t$, but this progress be arbitrarily small since $\stepmax$ can be arbitrarily small. \citet{lacoste2015global} shows that the number of drop steps for AFW is at most half of the number of iterations. Because we are maintaining two independent active sets in our formulation, we can obtain more drop steps, but we can still adapt their argument to obtain that the number of drop steps for SP-AFW is at most two thirds the number of iterations (assuming that the algorithm is initialized with only one atom per active set). In the SP-AFW algorithm, either a FW step is jointly made on both blocks, or an away-step is done on both blocks. Let us call~$A_t$ the number of FW steps (which potentially adds an atom in $S^{(t)}_x$ and $S^{(t)}_y$) and $D_t^{(x)}$ (resp $D_t^{(y)}$) the number of steps that removed at least one atom from $S^{(t)}_x$  ($S^{(t)}_y$). Finally, we call $D_t$ the number of \emph{drop steps}, i.e., the number of \emph{away steps} where at least one atom from $S^{(t)}_x$ or $S^{(t)}_y$ have been removed (and thus $\stepsize_t = \stepmax$ for these). Because a step is either a FW step or an away step, we have:
     \begin{equation} \label{eq:NumberSteps}
     A_t + D_t \leq t \, .
     \end{equation}
     We also have that $D_t^{(x)} + D_t^{(y)} \geq D_t$ by definition of $D_t$. Because a FW step adds at most one atom in an active set while a drop step removes one, we have (supposing that $|S^{(0)}_x|=| S^{(0)}_y|=1$):
    \begin{equation}
      1 + A_t -  D_t^{(x)} \geq |S^{(t)}_x| \quad \text{and} \quad 1 + A_t -  D_t^{(y)} \geq |S^{(t)}_y| .
    \end{equation}
    Adding these two relations, we get:
    \begin{equation}
    2 + 2 A_t \geq |S^{(t)}_x| + |S^{(t)}_y| + D_t^{(x)} + D_t^{(y)} \geq 2 + D_t \, ,
    \end{equation}
    using the fact that each active set as at least one element. We thus obtain
    $D_t \leq 2 A_t$. Combining with~\eqref{eq:NumberSteps}, we get:
    \begin{equation} \label{eq:upper_bound_drop_steps}
    D_t \leq \frac{2}{3} t \, ,
    \end{equation}
  as claimed.
\section{Affine invariant formulation of SP-FW}
  \label{sec:affine_invariant}
  In this section, we define the affine invariant constants of a convex function $f$ and their extension to a convex-concave function $\L$. 
  These constants are important as the FW-type algorithms are affine invariant if their step size are defined using affine invariant quantities. 
  We can upper bound these constants using the non affine invariant constants defined in the main paper. 
  Hence a convergence rate with affine invariant constants will immediately imply a rate with the constant introduced in the main paper.
  \subsection{The Lipschitz constants} %
    \label{sub:the_lipschitz_constants}
    We define the Lipschitz constant~$L$ of the gradient of the function $f$ with respect to the norm~$\| \cdot\|$ by using a dual pairing of norms, i.e. $L$ is a constant such that 
    \begin{equation}
       \label{def:lip}
        \forall \x, \x' \in \X, \qquad \|\nabla f(\x) - \nabla f (\x') \|_* \leq L \| \x - \x'\|,
     \end{equation} 
     where $\|\y\|_* := \sup_{\x \in \R^d, \|\x\| \leq 1} \y^{T} \x$ is the dual norm of $\|\cdot \|$. For a convex-concave function, we also consider the partial Lipschitz constants with respect to different blocks as follows.

     For more generality, we consider the dual pairing of norms $({\|\cdot\|_\X}, {\|\cdot\|_{\X^*}})$ on $\X$, and similarly $({\|\cdot\|_\Y}, {\|\cdot\|_{\Y^*}})$ on $\Y$. We also define the norm on the product space $\X \times \Y$ as the $\ell_1$-norm on the components: $\|(\x,\y)\|_{\X \times \Y} := \|\x\|_\X + \|\y\|_\Y$. We thus have that the dual norm of $\X \times \Y$ is the $\ell_\infty$-norm of the dual norms: $\|(\x,\y)\|_{(\X \times \Y)^*} = \max ( \|\x\|_{\X^*}, \|\y\|_{\Y^*})$.
     The \emph{partial} Lipschitz constants $L_{XX},L_{YY},L_{XY} \text{ and } L_{YX}$ of the gradient of the function $\L$ with respect to these norms are the constants such that for all $\x, \,\x' \in \X$ and $\y, \,\y' \in \Y$,
     \begin{equation}
     \label{def:lipL}
     \begin{aligned}
        \hspace{-3mm} \|\nabla_x \L (\x,\y) - \nabla_x \L(\x',\y) \|_{\X^*} &\leq L_{XX} \| \x - \x'\|_\X, &\! 
        \|\nabla_y \L (\x,\y) - \nabla_y \L(\x,\y') \|_{\Y^*} &\leq L_{YY}\,\| \y - \y'\|_\Y, \\
        \hspace{-3mm} \|\nabla_x \L (\x,\y) - \nabla_x \L(\x,\y') \|_{\X^*} &\leq L_{XY} \| \y - \y'\|_\Y,  
        &\!
        \|\nabla_y \L (\x,\y) - \nabla_y \L(\x',\y) \|_{\Y^*} &\leq L_{YX} \,\| \x - \x'\|_\X.
     \end{aligned}
     \end{equation}
  Note that the cross partial Lipschitz constants $L_{XY}$ and $L_{YX}$ do not necessarily use a dual pairing as $\X$ and $\Y$ could be very different spaces. 
  On the other hand, as the possibilities in~\eqref{def:lipL} are special cases of~\eqref{def:lip} when considering the $\ell_1$-norm of this product domain, one can easily deduce that the partial Lipschitz constants can always be taken to be smaller than the full Lipschitz constant for the gradient of $\L$, i.e., we have that $L \geq \max(L_{XX}, L_{XY}, L_{YX}, L_{YY})$.
    \subsection{The curvature: an affine invariant measure of smoothness} \label{sub:curv}
    
    To prove the convergence of the Frank-Wolfe algorithm, the typical affine invariant analysis proof in the FW literature assumes that the curvature of the objective function is bounded, where the curvature is defined by~\citet{jaggi2013revisiting} for example. 
    We give below a slight generalization of this curvature notion in order to handle the convergence analysis of FW with away-steps.\footnote{\label{foot:Cfcomment}The change is to consider the more general directions $\s - \vv$ instead of just $\s - \x$, and also any feasible positive step size. 
    See also Footnote~8 in~\citesup{lacoste2013affine} for a related discussion. A different (bigger) constant was required in~\citep{lacoste2015global} for the analysis of AFW because they used a line-search.} 
    It has the same upper bound as the traditional curvature constant (see Proposition~\ref{prop:C_fbounded}).

    \paragraph{Curvature.}[Slight generalization of~\citet{jaggi2013revisiting}]  
    \label{par:curvature}
    Let $f : \X \to \R$ be a convex function, we define the curvature $C_f$ of $f$ as 
      \begin{equation} 
    \label{CurvB}
      \Cf := 
        \sup_{\substack{\scalebox{0.65}{
        $\begin{matrix}
                    \x, \s, \vv \in \X , \\ \stepsize > 0 \text{ s.t.} \\
                    \x_{\stepsize} := \x + \stepsize \dd \in \X \\ 
                    \text{with } \dd := \s - \vv 
                  \end{matrix}$
          }}}
          \frac{2}{\stepsize^2}( f(\x_{\stepsize})-f(\x)- \stepsize \prodscal{\dd}{ \nabla f(\x)}). 
          \end{equation}
    Note that only the \emph{feasible} step sizes~$\stepsize$ are considered in the definition of~$\Cf$, i.e., $\stepsize$ such that $\x_{\stepsize} \in \X$. 
    If the gradient of the objective function is Lipschitz continuous, the curvature is upper bounded. 
    \begin{proposition}[Simple generalization of Lemma~7 in~\citet{jaggi2013revisiting}]
        \label{prop:C_fbounded}
    Let $f$ be a convex and continuously differentiable function on~$\X$ with its gradient $\nabla f$ L-Lipschitz continuous w.r.t. some norm $\|.\|$ in dual pairing over the domain $\X$. Then
      \begin{equation} \label{eq:CfUpperBound}
      C_f \leq D_\X^2 L \, ,
      \end{equation}
      where $D_\X : = \sup_{\x,\x' \in \X} \|\x - \x'\|$ is the diameter of $\X$.
  \end{proposition}  
  \proof[Lemma 1.2.3 in~\citetsup{nesterov2004introductory}, \citet{jaggi2013revisiting}] 
    Let $\x, \s, \vv \in \X$, set $\dd := \s - \vv$ and $\x_{\stepsize} = \x + \gamma \dd $ for some $\stepsize > 0$ such that $\x_{\stepsize} \in \X$. Then by the fundamental theorem of calculus, 
    \begin{equation}
      f(\x_{\stepsize}) = f(\x) + \int_0^\gamma \prodscal{ \dd }{\nabla f(\x + t \dd )} dt.
    \end{equation}
    Hence, we can write 
    \begin{align}
      f(\x_{\stepsize})-f(\x)- \stepsize \prodscal{\dd}{ \nabla f(\x)}
      & = \int_0^\gamma \prodscal{ \dd }{\nabla f(\x + t \dd) -  \nabla f(\x)} dt  \notag\\
      & \leq \|\dd\| \int_0^\gamma  \| \nabla f(\x + t \dd) - \nabla f(\x)\|_*dt \notag \\
      & \leq D_\X^2 L\int_0^\gamma t dt \notag \\ 
      & \leq \frac{\gamma^2}{2} D_\X^2 L.
    \end{align}
  Thus for all $\x,\s,\vv \in \X$ and $\x_{\stepsize} = \x + \stepsize (\s - \vv)$ for $\stepsize > 0$ such that $\x_{\stepsize} \in \X$, we have
  \begin{equation}
    \frac{2}{\gamma^2}( f(\x_{\stepsize})-f(\x)- \stepsize \prodscal{\s - \vv}{ \nabla f(\x)}) \leq L D_\X^2 .
  \end{equation}
  The supremum is then upper bounded by the claimed quantity. 
  \endproof
  
  \citetsup[Appendix C.1]{osokin2016gapBCFW} illustrate well the importance of the affine invariant curvature constant for Frank-Wolfe algorithms in their paragraph titled ``Lipschitz and curvature constants". They provide a concrete example where the wrong choice of norm for a specific domain~$\X$ can make the upper bound of Proposition~\ref{prop:C_fbounded} extremely loose, and thus practically useless for an analysis.
  
  We will therefore extend the curvature constant to the convex-concave
    function $\L$ by simply defining it as the maximum of the curvatures of the functions
    belonging to the family $\left(\x' \mapsto \L(\x',\y),\y' \mapsto
      -\L(\x,\y')\right)_{\x \in \X, \y \in \Y}$ (see Section~\ref{sub:curvature_and_interior_strong_convexity_constant_for_a_convex_concave_function}). But before that, we review affine invariant analogues of the strong convexity constants that will be useful for the analysis.

  \subsection{Affine invariant measures of strong convexity} %
  \label{sub:affine_invariant_measures_of_strong_convexity}

  In this section, we review two affine invariant measures of strong convexity that were proposed by~\citetsup{lacoste2013affine}\citep{lacoste2015global} for the affine invariant linear convergence analysis of the standard Frank-Wolfe algorithm (using the ``interior strong convexity constant'') or the away-step Frank-Wolfe algorithm (using the ``geometric strong convexity constant''). We will re-use them for the affine invariant analysis of the convergence of SP-FW or SP-AFW algorithms. In a similar way as the curvature constant~$\Cf$ includes information about the constraint set~$\X$ and the Lipschitz continuity of the gradient of~$f$ together, these constants both include the information about the constraint set~$\X$ and the strong convexity of a function~$f$ together.

  \paragraph{Interior strong convexity constant.} %
  \label{par:interior_strong_convexity_constant_}
  [based on \citetsup{lacoste2013affine}]
  Let $\xc$ be a point in the relative interior of $\X$.
   The \emph{interior strong convexity constant} for~$f$ with respect to the reference point~$\xc$ is defined as
    \begin{equation}
  \label{def:int_strong}
    \mu^{\xc}_f :=  \inf_{\substack{\scalebox{0.65}{
        $\begin{matrix}
                  \x \in \X \setminus \{\xc\} \\
                  \s = \bar{\s}(\x,\xc,\X) \\
                  \gamma \in (0,1], \\
                  \z = \x + \gamma(\s-\x)
        \end{matrix}$
        }}}
       \frac{2}{\gamma^2} \left(f(\z)-f(\x)- \prodscal{\z-\x}{ \nabla f(\x)}\right).
       \end{equation}
  Here, we follow the notation of \citetsup{lacoste2013affine} and take the point $\s$ to be the point where the ray from $\x$ to the reference point $\xc$ pinches the boundary of the set $\X$,
    i.e. $\bar{\s}(\x,\xc,\X) := {\rm ray}(\x,\xc) \cap \partial \X$, where $ \partial \X$ is the boundary of the convex set~$\X$.
    
  We note that in the original definition~\citepsup{lacoste2013affine}, $\xc$ was the (unique) optimum point for a strongly convex function~$f$ over~$\X$. The optimality of~$\xc$ is actually not needed in the definition and so we generalize it here to any point~$\xc$ in the relative interior of~$\X$, as this will be useful in our convergence proof for SP-FW.

  For completeness, we include here the important lower bound from~\citepsup{lacoste2013affine} on the interior strong convexity constant in terms of the strong convexity of the function~$f$.

  \begin{proposition}[Lower bound on~$\mu^{\xc}$ from {\citetsup[Lemma~2]{lacoste2013affine}}] \label{prop:delta}
  Let $f$ be a convex differentiable function and suppose that $f$ is strongly convex w.r.t. to some arbitrary norm $\norm{\cdot}$ over the domain~$\X$ with strong-convexity constant $\mu_f > 0$. Furthermore, suppose that the reference point $\xc$ lies in the relative interior of~$\X$, i.e.,
  $\delta_c := \min_{\s \in \partial \X}
    \; ||\s-\xc|| >0$. Then the interior strong convexity constant~$\mu^{\xc}_f$~\eqref{def:int_strong} is lower bounded as follows:
    \begin{equation} \mu^{\xc}_f \geq \mu_f \delta_c^2. \end{equation}
  \end{proposition}

  \proof 
  Let $\x$ and $\z$ be defined as in~\eqref{def:int_strong}, i.e., $\z = \x + \stepsize (\s - \x)$ for some $\stepsize > 0$ and where $\s$ intersects the boundary of $\X$ with the ray going from $\x$ to $\xc$.
  By the strong convexity of~$f$, we have
    \begin{equation} \label{eq:prop7strongConvex}
    f(\z)-f(\x) - \prodscal{\z-\x}{\nabla f(\x)} \geq ||\z-\x||^2 \frac{\mu_f}{2}
    = \stepsize^2 ||\s-\x||^2 \frac{\mu_f}{2}.
    \end{equation}
  From the definition of $\s$, we have that $\xc$ lies between $\x$ and $\s$ and thus:
  $||\s-\x||
  \geq ||\s - \xc|| \geq \delta_c$. Combining with~\eqref{eq:prop7strongConvex}, we conclude
    \begin{equation} f(\z)-f(\x) - \prodscal{\z-\x}{\nabla f(\x)} \geq
   \stepsize^2 \delta_c^2 \frac{\mu_f}{2}, \end{equation}
  and therefore 
    \begin{equation} \mu^{\xc}_f \geq \delta_c^2 \mu_f. \end{equation}
  \endproof

  We now present the affine invariant constant used in the global linear convergence analysis of Frank-Wolfe variants when the convex set~$\X$ is a polytope. The \emph{geometric strong convexity constant} was originally introduced by~\citetsup{lacoste2013affine} and~\citep{lacoste2015global}. To avoid any ambiguity, we will re-use their definitions verbatim in the rest of this section, starting first with a few geometrical definitions and then presenting the affine invariant constant. In these definitions, they assume that a \emph{finite} set~$\A$ of vectors (that they call \emph{atoms}) is given such that~$\X = \conv(\A)$ (which always exists when~$\X$ is a polytope).  

  \paragraph{Directional Width.}[\citet{lacoste2015global}]
  The \emph{directional width} of a set $\A$ with respect to a direction $\r$
  is defined as $\dirW(\A,\r) := \max_{\s, \vv \in\A}
  \big\langle \frac{\r}{\normEucl{\r}}, \s - \vv \big\rangle$. The \emph{width} of~$\A$
  is the minimum directional width over all possible directions in its affine
  hull.

  \paragraph{Pyramidal Directional Width.}[\citet{lacoste2015global}] We define the \emph{pyramidal directional
  width} of a set $\A$ with respect to a direction $\r$ and a base point
  $\x \in \domain$ to be
  \begin{equation} \label{eq:TruePdirW}
  \PdirW(\A,\r, \x) := \min_{\SS \in \SS_{\x}} \dirW( \SS \cup
  \{\s(\A,\r) \} , \; \r) = \min_{\SS \in \SS_{\x}} \max_{\s \in
  \A, \vv \in \SS} \textstyle \big\langle \frac{\r}{\normEucl{\r}}, \s - \vv \big\rangle
  ,
  \end{equation}
  where $\SS_{\x} := \{ \SS \, | \, \SS \subseteq \A$ such that $\x$ is a
  proper\footnote{By \emph{proper} convex combination, we mean that all
  coefficients are non-zero in the convex combination.} convex combination of
  all the elements in $\SS\}$, and $\s(\A,\r) := \argmax_{\vv \in
  \A} \innerProdCompressed{\r}{\vv}$ is the FW atom used as a summit, when using the convention in this section that $\r := -\nabla f(\x)$.

  \paragraph{Pyramidal Width.}[\citet{lacoste2015global}]
  To define the pyramidal width of a set, we take the minimum over the cone of
  possible \emph{feasible} directions~$\r$ (in order to avoid the problem of zero width).\\
  A direction~$\r$ is \emph{feasible} for $\A$ from $\x$ if it points inwards $\conv(\A)$, 
  (i.e. $\r \in \text{cone}(\A-\x)$).\\
  We define the \emph{pyramidal width} of a set $\A$ to be the smallest pyramidal width of all its faces, i.e.
  \begin{equation} \label{eq:Pwidth}
  \PWidth(\A) := \displaystyle \min_{\substack{\Kface \in \textrm{faces}(\conv(\A)) \\
                            \x \in \Kface \\
                            \r \in \text{cone}(\Kface-\x) \setminus \{\0\}} 
                                     } \PdirW(\Kface \cap \A,\r, \x)                   .                               
  \end{equation}

\paragraph{Geometric strong convexity constant.}[\citet{lacoste2015global}] %
\label{par:strong_}
   The \emph{geometric strong convexity constant} of~$f$ (over the set of atoms $\A$ which is left implicit) is: 
    \begin{equation} \label{def:geom_strong}
    \mu^{\away}_f := 
            \underset{ \x \in \X}{\inf}
            \inf_{\substack{\scalebox{0.6}{
            $\begin{matrix}
                          \x^* \in \X  \\
                          s.t \;   \prodscal{\nabla f(\x)}{\x^* - \x} <0
             \end{matrix}$
              }}}
          \frac{2}{\gamma^\away(\x,\x^*)^2}( f(\x^*)-f(\x)- \prodscal{\x^*-\x}{ \nabla f(\x)}) 
      \end{equation}
  where $\gamma^\away(\x,\x^*) := \frac{\prodscal{-\nabla f(\x)}{\x^*-\x}}{\prodscal{-\nabla f(\x)}{\s_f(\x)-\vv_f(\x)}}$ and $\X = \conv(\A)$. The quantity $\s_f(\x)$ represents the FW corner picked when running the FW algorithm on~$f$ when at~$\x$; while $\vv_f(\x)$ represents the worst-case possible away atom that AFW could pick (and this is where the dependence on $\A$ appears). We now define these quantities more precisely. Recall that the set of possible active sets is $\SS_{\x} := \{ \SS \, | \, \SS \subseteq \A$ such that $\x$ is a
  proper convex combination of
  all the elements in $\SS\}$.
  For a given set~$\SS$, we write $\vv_{\SS}(\x) := \argmax_{\vv \in \SS }
  \left\langle \nabla f(\x), \vv \right\rangle$ for the away atom in the
  algorithm supposing that the current set of active atoms is $\SS$.
  Finally, we define $\vv_f(\x) := \hspace{-3mm}\displaystyle\argmin_{\{\vv = \vv_{\SS}(\x)
  \,|\, \SS \in \SS_{\x} \}} \textstyle \hspace{-3mm}\left\langle \nabla f(\x), \vv
  \right\rangle$ to be the worst-case away atom (that is, the atom which would
  yield the smallest away descent). An important property coming from this definition that we will use later is that for $\st$ and $\vt$ being possible FW and away atoms (respectively) appearing during the AFW algorithm (consider Algorithm~\ref{alg:AFW} ran only on~$\X$), then we have:
  \begin{equation}
    \label{eq:gammaA}
    \gap^{\PW} := \prodscal{\st - \vt}{- \nabla f (\xt)} \geq \prodscal{\s_f(\xt) - \vv_f(\xt)}{- \nabla f (\xt)}.
  \end{equation}

  The following important theorem from~\citep{lacoste2015global} lower bounds the geometric strong convexity constant of~$f$ in terms of both the strong convexity constant of $f$, as well as the pyramidal width of~$\X = \conv{(\A)}$ defined as $\PWidth(\A)$~\eqref{eq:Pwidth}.

  \begin{proposition}[Lower bound for~$\strongConvAFW$ from {\citet[Theorem 6]{lacoste2015global}}]
  \label{thm:muFdirWinterpretation2}
  Let $f$ be a convex differentiable function and suppose that $f$ is
  $\mu$-\emph{strongly convex} w.r.t. to the Euclidean norm
  $\normEucl{\cdot}$ over the domain $\domain=\conv(\A)$ with strong-convexity constant $\mu
  \geq 0$. Then
  \begin{equation} \label{eq:muDirWidthBound}
  \strongConvAFW \geq \mu \cdot \left( \PWidth(\A) \right)^2.
  \end{equation}
  \end{proposition}
  
  The pyramidal width~\eqref{eq:Pwidth} is a geometric quantity with a somewhat intricate definition. Its value is still unknown for many sets (though always strictly positive for finite sets), but \citet[Lemma~4]{lacoste2015global} give its value for the unit cube in $\R^d$ as $1/{\sqrt{d}}$.

\subsection{Curvature and interior strong convexity constant for a convex-concave function} %
\label{sub:curvature_and_interior_strong_convexity_constant_for_a_convex_concave_function}

  In this subsection, we propose simple convex-concave extensions of the definitions of the affine invariant constants defined introduced in the two previous sections.
    
  To define the convex-concave curvature, we introduce the sets $\mathcal{F}$ and $\mathcal G$ of the marginal convex functions. 
    \begin{equation}\label{eq:F_x}
       \mathcal{F} := \{\x' \mapsto \L(\x',\y)  \}_{\y \in \Y} 
       \quad \text{and} \quad
        \mathcal{G} := \{\y' \mapsto -\L(\x, \y') \}_{\x \in \X} .
    \end{equation}
  Let $\L : \M \to \R$ a convex-concave function, we define the curvature pair $(C_{\L_x},C_{\L_y})$ of $\L$ as 
    \begin{equation} \label{eq:CLx}
    (C_{\L_x},C_{\L_y}) := \left(\underset{f \in \mathcal{F}}
        {\sup}\; C_{f},\underset{g \in \mathcal{G} }
        {\sup}\; C_{g} \right).
      \end{equation}
 and the curvature of $\L$ as
  \begin{equation}\label{eq:CL}
  C_\L := \frac{C_{\L_x}+C_{\L_y}}{2}. 
  \end{equation}
 An upper bound on this quantity follows directly from the upper bound on the convex case (Lemma~7 of \citet{jaggi2013revisiting}, repeated in our Proposition~\ref{prop:C_fbounded}) :

  \begin{proposition}\label{prop:Cbounded}
    Let $\L: \M \to \R$ be a differentiable convex-concave function. If $\X$ and $\Y$ are compact and $\nabla\L$ is Lipschitz continuous, then the curvature of $\L$ is bounded by $\frac{1}{2}(L_{XX} D_\X^2+L_{YY}D_\Y^2)$, where $L_{XX}$ (resp $L_{YY}$) is the largest Lipschitz constant respect to $\x$ ($\y$) of $\x \mapsto \nabla_x \L(\x,\y)$ ($\y \mapsto \nabla_y \L(\x,\y)$).
  \end{proposition}
  \proof
    Let $f$ in $\mathcal{F}$,
    \begin{equation}
     C_f \leq Lip(\nabla f) D_\X^2 \leq L_{XX} D_\X^2.
     \end{equation}
     Similarly, let $g$ in $\mathcal{G}$,
     \begin{equation}
     C_g \leq Lip(\nabla g) D_\Y^2 \leq L_{YY} D_\Y^2.
     \end{equation}
     Consequently,
     \begin{equation}
     C_\L = \frac{1}{2}(\underset{f \in \mathcal{F}}{\sup} C_f + \underset{g \in \mathcal{G}}{\sup} C_g) \leq \frac{1}{2} (L_{XX}D_\X^2 + L_{YY} D_\Y^2).
     \end{equation} 
     Where $D_\X$ and $D_\Y$ are the respective diameter of $\X$ and $\Y$.
  \endproof
  Note that $L_{XX}$ and $L_{YY}$ are upper bounded by the global Lipschitz constant of $\nabla \L$.
  Similarly, we define various notions of strong convex-concavity in the following.
\paragraph{Uniform strong convex-concavity constant.} %
\label{par:uniform_strong_convex_concavity_constant}
    The uniform strong convex-concavity constants is defined as 
    \begin{equation}
  \label{def:strongL}
      (\mu_\X,\mu_\Y) := \left(\underset{f \in \mathcal{G}}
        {\inf} \; \mu_f,\underset{g \in \mathcal{G}}
        {\inf} \; \mu_g\right)
    \end{equation}
    where $\mu_f$ is the strong convexity constant of $f$ and $\mu_g$ the strong convexity of $g$.

  Under some assumptions this quantity is positive.

  \begin{proposition} \label{prop:strongpos}
    If the second derivative of $\L$ is continuous, $\X$ and $\Y$ are compact and if for all $f \in \mathcal{F}\cup \mathcal G, \,\mu_f >0$, then
     $\mu_\X$ and $\mu_\Y$ are positive.  
  \end{proposition}

  \proof
    Let us introduce $H_x(\x,\y):= \nabla_x^2 \L(\x,\y)$ the Hessian of the function $\x \mapsto \L(\x,\y)$. We want to show that the smallest eigenvalue is uniformly bounded on $\M$. We know that the smallest eigenvalue lower bounds $\mu_\X$,
      \begin{equation} 
        \mu_{\X} \geq \inf_{\substack{\scalebox{0.65}{
                $\begin{array}{ll}
                    \; (\x,\y) \in \M \\
                    \|\u\|_2=1
                \end{array}$
                }}} \prodscal{\u}{H_x(\x,\y) \cdot\u}.
      \end{equation}
    But $H_x(\cdot)$ is continuous (because $\nabla_x^2 \L(\cdot)$ is continuous by assumption) and then the function $(\u,\x,\y) \mapsto \prodscal{\u}{H_x(\x,\y) \cdot\u}$ is continuous. 
    Hence since $\M$ and the unit ball are compact, the infimum is a minimum which can't be 0 by assumption. Hence $\mu_\X$ is positive. 
    Doing the same thing with the smallest eigenvalue of $-\nabla_y^2 \L(\x,\y)$, we get that 
    $\mu_\Y>0$.
  \endproof
  A common  family of saddle point objectives is of the form $f(x) +x^TMy - g(y)$. In this case, we get simply that $\left( \mu_\X, \mu_\Y \right) = (\mu_f,\mu_g).$
  An equivalent definition for the uniform strong convex-concavity constant is: $\L$ is $(\mu_\X,\mu_\Y)$-uniform strongly convex-concave function if
  \begin{equation}
    \left( \x,\y \right) \mapsto\L(\x,\y) - \frac{\mu_\X}{2} \|\x\|^2 + \frac{\mu_\Y}{2} \|\y\|^2
  \end{equation}
  is convex-concave.

  The following proposition relates the distance between the saddle point and the values of the function. It is a direct consequence from the uniform strong convex-concavity definition~\eqref{def:strongL}.
  \begin{proposition}\label{prop:strongL}
    Let $\L$ be a uniformly strongly convex-concave function and $(\x^*,\y^*)$ the saddle point of $\L$. 
    Then we have for all $\x$ in $\X$ and $\y \in \Y$,
    \begin{equation} \sqrt{\L(\x,\y^*) - \L^*} \geq \|\x^*-\x\|\sqrt{\frac{\mu_\X}{2}} \quad  
     \text{and} \quad
     \sqrt{\L^* - \L(\x^*,\y)} \geq \|\y^*-\y\|\sqrt{\frac{\mu_\Y}{2}}.
    \end{equation} 
  \end{proposition}

  \proof 
    The saddle point $(\x^*,\y^*)$ is the optimal point of the two strongly convex functions $\x \mapsto \L(\x,\y^*)$ and the function $\y \mapsto -\L(\x^*,\y)$, so we can use the property of strong convexity on each function and the fact that $\mu_\X$ lower bounds the strong convexity constant of $\L(\cdot,\y^*)$ (and similarly for $\mu_\Y$ with $-\L(\x^*,\y)$) as per the definition~\eqref{def:strongL}, to get the required conclusion.
  \endproof
  Now we will introduce the uniform strong convex-concavity constants relatively to our saddle point.
\paragraph{Interior strong convex-concavity.} %
\label{par:interior_strong_convex_concavity}
   The SP-FW interior strong convex-concavity constants (with respect to the reference point $(\xc, \yc)$) are defined as: 
    \begin{equation} 
  \label{def:interior_strongL}
    \left(\mu^{\xc}_\L,\mu^{\yc}_\L \right) := \left(\inf_{f \in \mathcal F} \mu^{\xc}_f, \inf_{g\in \mathcal G} \mu^{\yc}_g \right)
    \end{equation}
    where $\mu^{\xc}_f$ is the  interior strong convexity constant of $f$ w.r.t to the point $\xc$ and $ \mu^{\yc}_g$ is the interior strong convexity constant w.r.t to the point $\yc$. The sets $\mathcal F$ and $\mathcal G$ are defined in~\eqref{eq:F_x}. We also define the smallest quantity of both (with the reference point $(\xc, \yc)$ implicit):
    \begin{equation}
    \muIntL = \min\{\mu^{\xc}_\L, \mu^{\yc}_\L\}.
    \end{equation}
    We can lower bound this constant by a quantity depending on the uniform strong convexity constant and the distance of the saddle point to the boundary.
   The propositions on the strong convex-concavity directly follow from the previous definitions and the analogous proposition on the convex case (Proposition~\ref{prop:delta})

  \begin{proposition}\label{prop:deltaL}
    Let $\L$ be a convex-concave function. If the reference point
    $(\xc,\yc)$ belongs to the relative interior of $\X \times \Y$ and if the function
    $\L$ is strongly convex-concave with a strong convex-concavity constant $\mu > 0$,
    then $\muIntL$ is lower bounded away from zero. More precisely, define $\delta_x :=\min_{\s_x\in \partial\X } \; \|\s_x-\xc\|>0$ and $\delta_y := \min_{\s_y \in \partial\Y}\|\s_y-\yc\|$. Then we have,
    \begin{equation} \mu^{\xc}_\L \geq  \mu_\X \delta_x^2 
    \qquad \text{and}\qquad
    \mu^{\yc}_\L \geq \mu_\Y \delta^2_y. \end{equation}
  \end{proposition}

  \proof 
     Using the Proposition~\ref{prop:delta} we have, 
      \begin{equation}
        \mu^{\xc}_f \geq \mu_f \cdot \delta_x^2 \geq \mu_\X\cdot \delta_x^2,
      \end{equation}
     and,
        \begin{equation}
        \mu^{\yc}_g \geq \mu_g \cdot \delta_y^2 \geq \mu_\Y\cdot \delta_y^2.
      \end{equation}
  \endproof
  When the saddle point is not in the interior of the domain, we define next a constant that takes in consideration the geometry of the sets. If the sets are polytopes, then this constant is positive.
\paragraph{Geometric strong convex-concavity.} %
\label{par:geometric_strong_convex_concavity}
   The SP-FW \emph{geometric} strong convex-concavity constants are defined analogously as the interior strong convex-concavity constants,
    \begin{equation} 
  \label{def:geom_strongL}
      \left(\mu^{\away}_{\L_x}, \mu^{\away}_{\L_y}\right) := \left( \min_{f \in \mathcal F} \mu^{\away}_f,  \min_{g \in \mathcal G} \mu^{\away}_g \right) \, ; \quad \mu^{\away}_\L:=\min \left(\mu^{\away}_{\L_x}, \mu^{\away}_{\L_y}\right)  ,
    \end{equation}
    where $\mu^{\away}_f$ is the  geometric strong convexity constant of $f \in \mathcal F$ (over $\A$) as defined in~\eqref{def:geom_strong} (and similarly $\mu^{\away}_g$ is the geometric strong convexity constant of $g \in \mathcal G$ over $\B$).

  It is straightforward to notice that the lower bound on the geometric strong convexity constant (Proposition~\ref{thm:muFdirWinterpretation2}) can be extended to the geometric strong convex-concavity constants (where $\mu_\X$ and $\mu_\Y$ are now assumed to be defined with respect to the Euclidean norm): 
  \begin{equation} \label{eq:muAwayLInequality}
  \mu^{\away}_{\L_x} \geq \mu_\X \PWidth(\A)^2 \quad \text{ and } \quad \mu^{\away}_{\L_y} \geq \mu_\Y \PWidth(\B)^2 .
  \end{equation} 

\subsection{The bilinearity coefficient}

In our proof, we need to relate 
  the gradient at the point $(\xt,\yt)$ with the one at
  the point $(\xt, \y^*)$. We can use the Lipschitz continuity of the gradient for this. We define below affine invariant quantities that can upper bound this difference.
\paragraph{Bilinearity coefficients.} %
\label{par:bilinearity_coefficients}
  let $\L$ be a strongly convex-concave function, and let $(\x^*,\y^*)$ be its unique saddle point. We define the bilinearity coefficients $(M_{XY},M_{YX})$ as, 
    \begin{equation}   \label{def:M}
    M_{XY} = \hspace{-2mm} \sup_{\substack{\scalebox{0.65}{
       $ \begin{matrix}
                 \y \in \Y \\
                 \x, \s,\vv \in \X \\
                 \dd = \s -\vv   %
               \end{matrix}$
        }}} \hspace{-3mm} \prodscal{\dd}{\frac{\nabla_x\L(\x,\y^*)-\nabla_x\L(\x,\y)}{\sqrt{\L^*-\L(\x^*,\y)}}} 
    \end{equation}
  and, 
  \begin{equation}
      M_{YX}:=  \hspace{-2mm} \sup_{\substack{\scalebox{0.65}{
       $ \begin{matrix}
                 \x \in \X\\
                 \y, \s,\vv \in \Y \\
               \dd = \s-\vv
          \end{matrix}$
        }}} \hspace{-3mm} \prodscal{\dd}{\frac{\nabla_y\L(\x,\y)-\nabla_y\L(\x^*,\y)}{\sqrt{\L(\x,\y^*)-\L^*}}}.
  \end{equation}
  We also define the global bilinearity coefficient as
  \begin{equation}\label{eq:def_bilin_coef}
    M_\L := \max\{M_{XY},M_{YX}\}.
  \end{equation}
  We can upper bound these affine invariant constants with the Lipschitz constant of the gradient, the uniform strong convex-concavity constants and the diameters of the sets.
  \begin{proposition}\label{Mbounded}
  If  $\X$ and $\Y$ are compact, $\nabla \L$ is Lipschitz continuous and $\L$ is uniformly strongly convex-concave with constants $(\mu_\X, \mu_\Y)$, then
    \begin{equation}
      M_{XY} \leq \sqrt{\frac{2}{\mu_\Y}}L_{XY} \cdot D_\X 
      \quad \text{and} \quad
      M_{YX} \leq \sqrt{\frac{2}{\mu_\X}}L_{YX} \cdot D_\Y
    \end{equation}
  where $L_{XY}$ and $L_{YX}$ are the \emph{partial} Lipschitz constants defined in Equation~\eqref{def:lipL}. The quantity $D_\X$ is the diameter of the compact set $\X$
  and $D_\Y$ is the diameter of $\Y$.
  \end{proposition}
  \proof
   \begin{align*}
   M_{XY}
   &= \sup_{\substack{\scalebox{0.65}{
          $ \begin{matrix}
                    \y \in \Y \\
                    \x, \s,\vv \in \X \\
                    \dd = \s -\vv 
                  \end{matrix}$
     }}} \;\prodscal{\dd}{\frac{\nabla_x\L(\x,\y^*)-\nabla_x\L(\x,\y)}{\sqrt{\L^*-\L(\x^*,\y)}}} \\        
    &\leq \sup_{\substack{\scalebox{0.65}{
              $ \begin{matrix}
                        \y \in \Y \\
                        \x, \s,\vv \in \X \\
                        \dd = \s -\vv 
                      \end{matrix}$
         }}} \; {\frac{ \|\dd\|_\X \, \|\nabla_x\L(\x,\y^*)-\nabla_x\L(\x,\y) \|_{\X^*}}{\sqrt{\L^*-\L(\x^*,\y)}}} \\
    & \leq \sup_{\substack{\scalebox{0.65}{
              $ \begin{matrix}
                        \y \in \Y \\
                        \s,\vv \in \X \\
                        \dd = \s -\vv 
                      \end{matrix}$
         }}} \; {\frac{ \|\dd\|_\X \,  L_{XY} \|\y^*-\y \|_\Y}{\sqrt{\L^*-\L(\x^*,\y)}}} \\
     & \leq  \sup_{\y \in \Y} D_\X L_{XY} \frac{\|\y^*-\y\|_\Y}{\sqrt{\L^*-\L(\x^*,\y)}}.
  \end{align*}
  Then using the relation between $\|\y^*-\y\|_\Y$ and $\sqrt{\L^*-\L(\x^*,\y)}$ due to strong convexity (Proposition~\ref{prop:strongL})
  \begin{equation}
  M_{XY} \leq \sqrt{\frac{2}{\mu_\Y}}L_{XY}\cdot D_\X.
  \end{equation}
  We use a similar argument for $M_{YX}$ which allows us to conclude.
  \endproof

\subsection{Relation between the primal suboptimalities} %
\label{sub:relation_between_the_primal_errors}
In this section, we are going to show that if the objective function $\L$ is uniformly strongly convex-concave, then we have a relation between $h_t$ and $w_t$. First let us introduce affine invariant constants to relate these quantities (in the context of a given saddle point $(\x^*, \y^*)$):
\begin{equation}
 \label{def:constant_relation_primal}
  P_\X := \sup_{\x \in \X} \frac{\prodscal{\nabla_x \L (\x,\hat \y(\x))}{\x-\x^*}}{\sqrt{\L(\x,\y^*)-\L(\x^*,\y^*)}}
\quad
\text{and}
\quad
  P_\Y := \sup_{\y \in \Y} \frac{\prodscal{\nabla_y \L (\hat \x (\y),\y)}{\y-\y^*}}{\sqrt{\L(\x^*,\y^*)-\L(\x^*,\y)}},
\end{equation}
where $\hat \y (\x) : = \argmax_{\y \in \Y} \L(\x,\y)$ and $\hat \x (\y) : = \argmin_{\x \in \X} \L(\x,\y)$. We also define: 
\begin{equation} 
  P_\L := \max \{P_\X, P_\Y\}.
\end{equation}
These constants can be upper bounded by easily computable constants.
\begin{proposition} For any $(\mu_\X,\mu_\Y)$-uniformly convex-concave function $\L$,
\label{prop:upper_bound_P_L}
\begin{equation}
  P_\X \leq \sqrt{\frac{2}{\mu_\X}}\sup_{\z \in \M}\| \nabla_x \L(\z)\|_{\X^*}
\quad
\text{and}
\quad
   P_\Y \leq \sqrt{\frac{2}{\mu_\Y}} \sup_{\z\in \M}\|\nabla_y \L (\z)\|_{\Y^*}.
\end{equation}
\end{proposition}
\proof 
Let us start from the definition of $P_\X$, let $\x \in \X$,
\begin{align*}
  \frac{\prodscal{\nabla_x \L (\x,\hat \y(\x))}{\x-\x^*}}{\sqrt{\L(\x,\y^*)-\L(\x^*,\y^*)}} 
  & \leq \frac{\|\x-\x^*\|_\X \cdot \sup \left(\| \nabla_x \L (\z)\|_{\X^*}\right) }{\sqrt{\L(\x,\y^*)-\L(\x^*,\y^*)}} \\
  & \leq \sqrt{\frac{ 2}{\mu_\X}} \sup_{\z \in \M}\| \nabla_x \L(\z) \|_{\X^*} \qquad \qquad\text{(by strong convexity.)} 
\end{align*}
The same way we can get 
\begin{equation}
  P_\Y \leq \sqrt{\frac{2}{\mu_\Y}} \sup_{\z\in \M}\|\nabla_y \L(\z) \|_{\Y^*}.
\end{equation}
It concludes our proof.
\endproof 
One way to compute an upper bound on the supremum of the gradient is to use any reference point $\bar{\z}$ of the set:
\begin{equation}
  \forall \bar{z} \in \M, \quad \sup_{\z \in \M}\| \nabla_x \L(\z) \|_{\X^*} 
  \leq \nabla_x \L(\bar{\z}) + L_{XX} D_\X + L_{XY} D_\Y.
\end{equation}
We recall that $L_{XX}$ is the largest (with respect to $\y$) Lipschitz constant of $\x \mapsto \nabla_x \L(\x,\y)$. Note that $L_{XX}$ is upper bounded by the global Lipschitz constant of $\nabla \L$. We can compute an upper bound on the supremum of the norm of $\nabla_y \L$ the same way.

With these above defined affine invariant constants,  we can finally relate the two primal suboptimalities as $h_t \leq \mathcal{O}(\sqrt{w_t})$.
\begin{proposition}
\label{prop:relation_primal}
For a $(\mu_\X,\mu_\Y)$-uniformly strongly convex-concave function $\L$,
\begin{equation}
  h_t \leq P_\L \sqrt{2w_t} 
  \quad 
  \text{and} \quad
  P_\L \leq \sqrt{2} \underset{\z \in \M}{\sup} \left\{\frac{\|\nabla_x \L(\z)\|_{\X^*}}{\sqrt{\mu_\X}} ,\frac{\|\nabla_y \L(\z)\|_{\Y^*}}{\sqrt{\mu_\Y}}   \right\}.
\end{equation}
\end{proposition}

\proof We will first work on $h_t^{(x)}$:
  \begin{align*}
   h_t^{(x)} &=\L(\xt,\ytm) -\L^* \\
    &\leq \L(\xt,\ytm) -\L(\x^*,\ytm)\\
    &\leq \prodscal{\xt - \x^*}{\nabla_x \L(\xt,\ytm} 
      \qquad\qquad \text{(by convexity)}\\
    & \leq P_\X \sqrt{\wt^{(x)}} \qquad (\text{def of } P_\X \: \eqref{def:constant_relation_primal}).
  \end{align*}
We can do the same thing for $h_t^{(y)}$ and $w_t^{(y)}$, thus 
\begin{equation}
h_t \leq P_\L \left(\sqrt{\wt^{(x)}} + \sqrt{\wt^{(y)}} \right) \leq P_\L \sqrt{2 w_t} ,    
\end{equation}
where the last inequality uses $\sqrt{a} + \sqrt{b} \leq \sqrt{2 (a+b)}$.
Finally, the inequality on $P_\L$ is from Proposition~\ref{prop:upper_bound_P_L}.
\endproof

\section{Relations between primal suboptimalities and dual gaps}\label{sec:relations_gaps}
\subsection{Primal suboptimalities}
  Recall that we introduced $\xtm := \argmin_{\x \in \X} \L(\x, \yt)$ 
  and similarly $\ytm := \argmax_{\y \in \Y} \L(\xt, \y).$
  Then the primal suboptimality is the positive quantity 
    \begin{equation}\label{def:primal}
    h_t := \L(\xt,\ytm)- \L(\xtm,\yt). 
    \end{equation}
  To get a convergence rate, one has to upper bound the primal suboptimality defined in~\eqref{def:primal}, but it is hard to work with the moving quantities $\xtm$ and $\ytm$ in the analysis. This is why we use in our analysis a different merit function that uses the (fixed) saddle point $(\x^*,\y^*)$ of $\L$ in its definition. We recall its definition below.

   \paragraph{Second primal suboptimality.} %
   \label{par:second_primal_gap}
  We define the second primal suboptimality for $\L$ of the iterate $(\xt, \yt)$ with respect to the saddle point $(\x^*,\y^*)$ as the positive quantity: 
    \begin{equation} 
  \label{def:primal2}
    \wt := \L(\xt,\y^*)- \L(\x^*,\yt). 
    \end{equation}
  It follows from $\L(\xt,\ytm) \geq \L(\xt,\y^*)$ and $\L(\x^*,\yt) \geq \L(\xtm,\yt)$ that $w_t \leq h_t$. 
  Furthermore, under the assumption of uniform strong convex-concavity, we proved in Proposition~\ref{prop:relation_primal} that the square root of $w_t$ upper bounds $h_t$ up to a constant.

\subsection{Gap inequalities} \label{app:gapInequalities}

  In this section, we will prove the crucial inequalities relating suboptimalities and the gap function.
  Let's recall the definition of $\st$ and $\vt$:
    \begin{equation}\st:=  \argmin_{\s \in \M}
         \prodscal{\s}{
       \rt}
    \quad \text{  and } \quad 
    \vt:=  \argmax_{\vv \in \SS^{t}_x\times \SS^{t}_y}
         \prodscal{\vv}{
        \rt}
         \end{equation}
  where $(\rt)^\top := ( (\rt_x)^\top, (\rt_y)^\top ) :=  \left(
        \nabla_x \L(\xt,\yt),
         -\nabla_y \L(\xt,\yt)
         \right)$. 
  Also, the following various gaps are defined as 
    \begin{equation}\label{eq:gaps} 
    \gap^{\FW} :=\prodscal{\dt_{\FW}}{-\rt} 
    , \qquad
    \gap^{\PW} :=\prodscal{\dt_{\PW}}{-\rt}
    \quad \text{and} \quad
    \gap := \prodscal{\dt}{-\rt}
    \end{equation}
  where $\dt_{\FW} = \st - \zt$ and $\dt_{\PW} = \st - \vt$. The direction $\dt$ is the direction chosen by the algorithm at step $t$: it is always $\dt_{FW}$ for SP-FW, and can be either $\dt_{FW}$ or $\dt_{\away} := \zt - \vt$ for SP-AFW. Even if the definitions of these gaps are different, the formalism for the analysis of the convergence of both algorithms is going to be fairly similar.
  It is straightforward to notice that $\gap^{\PW} \geq \gap$ and one can show that the current gap $\gap$ is lower bounded by half of $\gap^{\PW}$:
  \begin{lemma}\label{lemma:gap_FW_PW}
     For the SP-AFW algorithm, the current gap $\gap$ can be bounded as follows:
     \begin{equation}
      \frac{1}{2} \gap^{\PW} \leq \gap \leq \gap^{\PW}
     \end{equation}
   \end{lemma} 
   \proof
    First let's show the RHS of the inequality, 
    \begin{equation}
      \gap^{\PW} := \prodscal{\dt_{\PW}}{-\rt} = \prodscal{\dt_{\away}}{-\rt} + \prodscal{\dt_{\FW}}{-\rt} \geq \prodscal{\dt}{-\rt}
    \end{equation}
    because both $\prodscal{\dt_{\away}}{-\rt}\geq 0$ and $ \prodscal{\dt_{\FW}}{-\rt} \geq 0$ from their definition.
    For the LHS inequality, we use the fact that $\gap = \max \left\{ \prodscal{\dt_{\away}}{-\rt} , \prodscal{\dt_{\FW}}{-\rt} \right\}$ for SP-AFW and thus:
    \begin{equation}
       \gap^{\PW} = \prodscal{\dt_{\away}}{-\rt} + \prodscal{\dt_{\FW}}{-\rt} \leq 2 \gap .
    \end{equation}
   \endproof
   In the following, we will assume that we are in one of the two following cases: 
   \leqnomode
  \begin{equation}
    \label{case:1}
    \text{The saddle point of } \L \text{ belongs to the relative interior of } \M.
    \tag{I}
  \end{equation}
  \begin{equation}
    \label{case:2}
    \X \text{ and } \Y \text{ are polytopes}, \quad \text{i.e.} \; \exists \A, \B \text{ finite s.t}\;\; \X = \conv(\A), \; \Y = \conv(\B).
    \tag{P}
  \end{equation}
  \reqnomode

  Then either $\muIntL>0$ (case~\ref{enu:situation1}) or $\mu^{\away}_\L>0$ (case~\ref{enu:situation2}).
  Let's write the gap function as the sum of two smaller gap functions:
    \begin{equation} \label{eq:gtGeneralized}
      \gap 
      = \quad \underbrace{\prodscal{ \dt_{(x)}}{-\rt_x}}_{=:\gap^{(x)}} \quad
      + \quad
    \underbrace{\prodscal{ \dt_{(y)}}{ -\rt_y}}_{=:\gap^{(y)}} 
    \end{equation}
  Because of the convex-concavity of $\L$, this scalar product bounds
  the differences between the value of $\L$ at the point
  $(\xt,\yt)$ and the value of $\L$ at another point. Hence this
  gap function upper-bounds $h_t$ and $w_t$ defined in~\eqref{def:primal} and~\eqref{def:primal2}. More concretely, we
  have the following lemma.

  \begin{lemma} \label{lemme:g}
    For all $t$ in $\N$, $\x\in \X$ and $ \y \in \Y$
      \begin{equation} \label{eq:encadre}
      \gap^{\PW} \geq \gap^{\FW} \geq  \L(\xt,\y)-\L(\x, \yt),
      \end{equation}
    and, furthermore,
      \begin{equation} \label{eq:gtInequality_with_ht}
      \gap \geq  h_t  \geq \wt.
      \end{equation}
  \end{lemma}
  \proof 
    First let's show the LHS of \eqref{eq:encadre}, 
    \begin{equation}
      \gap^{\PW} = \prodscal{\dt_{\PW}}{-\rt} = \prodscal{\dt_{\away}}{-\rt} + \prodscal{\dt_{\FW}}{-\rt} \geq \prodscal{\dt_{\FW}}{-\rt} = \gap^{\FW}
    \end{equation}
    because one can easily derive that $\prodscal{\dt_{\away}}{-\rt}\geq 0$ from the definition of the away direction $\dt_{\away}$.
    It follows from convexity of $\x \mapsto \L(\x, \yt)$ that for all $\x$ in $\X$, 
        \begin{align} \label{eq:gapx}
           (\gap^{\FW})_{x} :=\prodscal{ (\dt_{\FW})_x}{ -\gnx}
              & \geq \prodscal{ \x - \xt}{-\gnx}\\
              &\geq \L(\xt,\yt) - \L(\x, \yt). 
        \end{align}
    A similar inequality emerges through the convexity of $\y \mapsto -\L(\xt,\y)$,
      \begin{equation}
        (\gap^{\FW})_{y} := \prodscal{ (\dt_{\FW})_y}{\gny} \geq \L(\xt, \y) - \L(\xt,\yt),
      \end{equation}
    which gives us 
      \begin{equation} 
      \gap^{\FW} \geq \L(\xt, \yt) - \L(\x, \yt) + \L(\xt, \y) - \L(\xt , \yt) ,
      \end{equation}
     which shows~\eqref{eq:encadre}.
     By using $\x = \xtm$ and $\y = \ytm$ in~\eqref{eq:encadre}, we get $\gap^{\FW} \geq h_t$. We also know that $g_t = \max(\gap^\away, \gap^\FW) \geq \gap^{\FW}$ for SP-AFW. So combining with $h_t \geq w_t$ that we already knew, we get~\eqref{eq:gtInequality_with_ht}.
  \endproof
  Next, we recall two lemmas, one from~\citesup{lacoste2013affine} and the other one from \citep{lacoste2015global}. These lemmas upper bound the primal suboptimality with the square of the gap times a constant depending on the geometric (or the interior) strong convexity constant.

  \begin{lemma}[\citet{lacoste2015global}, \citetsup{lacoste2013affine}] \label{lemme:lingap}
   If $f$ is strongly convex, then for any $\xt \in \X$, 
    \begin{equation}
     f(\xt)-f(\xc) \leq \frac{\left( \gap^{\FW} \right)^2}{2 \mu^{\xc}_f}
     \quad
     \text{if } \xc \in \text{ interior of } \X
     \quad \text{\citepsup{lacoste2013affine}} 
     \end{equation}
  and 
     \begin{equation}
     f(\xt)-f^* \leq \frac{\left( \gap^{\PW} \right)^2}{2 \mu^{\away}_f}
     \quad
     \text{if } \X = \conv(\A) 
     \quad
     \text{\citep{lacoste2015global}}
    \end{equation}
  where $\gap^{\FW}= \prodscal{\xt-\st}{\nabla f(\xt)}$, $\gap^{\PW}= \prodscal{\vt-\st}{\nabla f(\xt)}$ and $f^* = \min_{x \in \X} f$.
  \end{lemma}
  Notice once again that in this lemma~\emph{we do not need $\xc$ to be optimal}.
  \proof 
  Let $\xt \neq \xc$. Using the definition of interior strong convexity~\eqref{def:int_strong} and choosing $\gamma$ such that $\xc = \xt + \gamma \left( \bar \s\left(\xt,\xc\right) - \xt \right)$, we get
    \begin{align*}
    f(\xc)-f(\xt) &\geq \gamma\prodscal{\bar\s(\xc,\xt)-\xt}{\nabla f(\xt)} + \gamma^2 \frac{ \mu^{\xc}_f}{2}\\
    & \geq - \gamma\gap^{\FW}+ \gamma^2 \frac{ \mu^{\xc}_f }{2}  \\
    & \geq  \frac{ - \left( \gap^{\FW} \right)^2 }{2\mu^{\xc}_f}.
    \end{align*}
  The last line of this derivation is obtained through the inequality: $-a^2 + 2ab - b^2 \leq 0$. 
  If $\xt = \xc$ the inequality is just the positivity of the gap.

  For the second statement, we will use the definition of the geometric strong convexity constant (Equation~\eqref{def:geom_strong}) at the point $\x = \xt$ and $\x^* \in \argmin_{\x \in \X} f(x)$. 
  Recall that $\gamma^\away(\x,\x^*) = \frac{\prodscal{-\nabla f(\x)}{\x^*-\x}}{\prodscal{-\nabla f(\x)}{\s_f(\x)-\vv_f(\x)}}$.
    \begin{align*}
    f(\x^*)-f(\xt) &\geq \prodscal{\x^*-\xt}{\nabla f(\xt)} + \frac{ \mu^{\away}_f}{2}\gamma^\away(\xt, \x^*)^2 \\
    & = -\gamma^\away(\xt, \x^*) \prodscal{\s_f(\xt) - \vv_f(\xt)}{-\nabla f(\xt)}  + \frac{ \mu^{\away}_f}{2} \gamma^\away(\xt, \x^*)^2  \\
    & \geq - \gamma^\away(\xt, \x^*)\gap^{\PW}+  \frac{ \mu^{\away}_f }{2}\gamma^\away(\xt, \x^*)^2 \qquad \text{(Equation~\eqref{eq:gammaA})}  \\
    & \geq  \frac{ - \left( \gap^{\PW} \right)^2 }{2\mu^{\away}_f}.
    \end{align*}
  \endproof
  Lemma~\ref{lemme:lingap} is useful to understand the following lemma and its proof which is just an extension to
  the convex-concave case.

  \begin{lemma}[Quadratic gap upper bound on second suboptimality for~\eqref{case:1} or~\eqref{case:2}]\label{lemme:gap}
   If $\L$ is a strongly convex-concave function, then for any $(\xt,\yt) \in \M$,
    \begin{equation}
     \wt  \leq \frac{(\gap^{\FW})^2}{2 \muIntL } 
     \quad 
      \text{for \eqref{case:1}}
     \quad \text{and} \quad
     \wt \leq h_t \leq \frac{(\gap^{\PW})^2}{2 \mu^{\away}_\L}
     \quad
     \text{for \eqref{case:2}} %
    \end{equation}
  where the gaps are defined in \eqref{eq:gaps}, $\muIntL := \min \{ \mu^{\x^*}_\L, \mu^{\y^*}_\L\}$ (i.e. using the reference points $(\xc,\yc) := (\x^*,\y^*)$ in the definition~\eqref{def:interior_strongL})
  and $\mu^{\away}_\L$ is the geometric strong convex-concavity of $\L$ over $\Vertices$, as defined in~\eqref{def:geom_strongL}.
  \end{lemma}

  \proof For \eqref{case:1}:

   Let the function $f$ on $\X$ be defined by $f(\x) = \L(\x',\yt)$, and the function $g$ on $\Y$ be $g(\y') = -\L(\xt, \y')$. Then using the Lemma~\ref{lemme:lingap} on the function $f$ with the reference point $\x^*$, and on $g$ with reference point $\y^*$, we get
    \begin{align*}
      \L(\xt, \yt) - \L(\x^*,\yt)  & \leq  \frac{ \prodscal{\st_x - \xt }{-\nabla_x \L(\xt,\yt)}^2 }{2 \mu_f^{\x^*}}  \\
      \L(\xt, \y^*) - \L(\xt,\yt)  & \leq   \frac{ \prodscal{\st_y - \yt }{\nabla_y \L(\xt,\yt)}^2 }{2 \mu^{\y^*}_g} .
    \end{align*}
    As $\muIntL$ is smaller than both $\mu_f^{\x^*}$ and $\mu_g^{\y^*}$ by the definition~\eqref{def:interior_strongL}, we can use it in the denominator of the above two inequalities.
    As we saw from Section~\ref{app:gapInequalities} in~\eqref{eq:gtGeneralized}, the gap can be split as sum of the gap of the block $\X$ and the gap of the block $\Y$, i.e. $\gap^{\FW} =  \prodscal{\st_x- \xt}{-\nabla_x \L(\xt,\yt)} +  \prodscal{\st_y -\yt}{\nabla_y \L(\xt,\yt)} $.
    Then, using the inequality: $a^2+b^2 \leq (a+b)^2$ for $(a  ,b \geq
    0)$, we obtain
      \begin{equation} \label{eq:proof_w_FWgap}
      w_t \leq \frac{ (\gap^{\FW})^2}{2 \muIntL}.
      \end{equation}

      For \eqref{case:2}:

    Using the Lemma~\ref{lemme:lingap} for case~\eqref{enu:situation2} on the same functions $f$ and $g$ defined above, we get
      \begin{align*}
        \L(\xt, \yt) - \L(\xtm,\yt)  & \leq  \frac{  \prodscal{\st_x - \vt_x }{\nabla_x \L(\xt,\yt)}^2 }{2 \mu^{\away}_{f}}  \\
        \L(\xt, \ytm) - \L(\xt,\yt)  & \leq   \frac{  \prodscal{\st_y -\vt_y }{-\nabla_y \L(\xt,\yt)}^2 }{2 \mu^{\away}_{g}}.
      \end{align*}
    Using a similar argument as the one to get~\eqref{eq:proof_w_FWgap}, using that $\mu^{\away}_\L$ is smaller than both $\mu^{\away}_{f}$ and $\mu^{\away}_{g}$, and referring to the separation of the gap~\eqref{eq:gtGeneralized}, we get
      \begin{equation} 
       h_t \leq \frac{ (\gap^{\PW})^2}{2 \mu^{\away}_\L}.
      \end{equation}

  \endproof

\section{Convergence analysis}
\label{appendix:analysis}
In this section, we are going to show two important lemmas. The first one shows that under some assumptions we can get a Frank-Wolfe-style induction scheme relating the second suboptimality of the potential update $w_{\gamma}$, the current value of the second suboptimality $w_t$, the gap~$\gap$ and any step size $\gamma \in [0,\stepmax]$. The second lemma will relate the gap and the square root of $w_t$; this relation enables us to get a rate on the gap after getting a rate on $w_t$.

\subsection{First lemmas}
  The first lemma in this section is inspired from the standard FW progress lemma, such as Lemma~C.2 in~\citep{lacoste2013block}, though it requires a non-trivial change due to the compensation phenomenon for~$\L$ mentioned in the main text in~\eqref{eq:oscilation}.
  In the following, we define the possible updated iterate $\z_\stepsize$ for $\stepsize \in [0, \stepmax]$:
  \begin{equation}
  \z_\gamma := (\x_\gamma, \y_\gamma):= \zt + \gamma \dt, \quad \text{ where } \quad \dt \text{ is the direction of the step.}
  \end{equation}
  For a FW step $\dt = \dt_{\FW}:= \st - \zt$ and for an away step $\dt = \dt_{\away} := \zt - \vt$. We also define the corresponding new suboptimality for $\z_\gamma$:
  \begin{equation}
    w_\gamma := \L(\x_\gamma, \y^*) - \L(\x^*,\y_\gamma).
  \end{equation}
  \begin{lemma}[Suboptimality progress for SP-FW and SP-AFW] \label{lemme:ineg}
  Let $\L$ be strongly convex-concave, 

    If we are in case~\eqref{case:1} and $\dt = \dt_{\FW}$ is a FW direction, we have for any $\stepsize \in [0,1]$:
      \begin{equation}  \label{eq:wInequalityCaseI}
      w_\gamma 
       \leq \wt  
      - \CondNumb^{\FW}\gamma \gap^{\FW}
      + \gamma^2 {C_\L},
      \quad \text{where} \quad 
      \CondNumb^{\FW} :=1- \frac{M_\L }{\sqrt{\muIntL}}.
      \end{equation}
    If we are in case~\eqref{case:2} and $\dt$ is defined from a step of SP-AFW (Algorithm~\ref{alg:AFW}), we have for any $\stepsize \in [0, \stepmax]$:
      \begin{equation} \label{eq:wInequalityCaseP} 
      w_\gamma 
       \leq \wt  
      - \CondNumb^{\PW} \gamma \gap^{\PW}
      + \gamma^2 {C_\L},
      \quad \text{where} \quad
      \CondNumb^{\PW} := \frac{1}{2}- \frac{M_\L }{\sqrt{\mu^{\away}_\L}}.
      \end{equation}
  \end{lemma}

  \proof 
    The beginning of the argument works with any direction $\dt$. 
    Recall that $\wt^{(x)} := \L(\xt, \y^*)- \L^*$ and $\wt^{(y)} := -\L(\x^*,\yt) + \L^*$.
    Now writing $\x_\gamma = \xt + \gamma \dt_x$ and using the definition of the curvature $C_f$~\eqref{CurvB} for the function $\x \mapsto f(\x) := \L(\x,\y^*)$, we get
      \begin{align}
      w_\gamma^{(x)} &:= \L(\x_\gamma,\y^*) - \L^* \notag\\
        &\leq \L(\xt, \y^*) - \L^* + \gamma \prodscal{\dt_x}{\nabla_x\L(\xt, \y^*)} + \gamma^2 \frac{C_\L}{2},
      \end{align} 
    since $C_\L \leq C_f$ by definition~\eqref{eq:CL}. 
    Recall that the gap function $\gap$ can be decomposed by~\eqref{eq:gtGeneralized} into two smaller gap functions $\gap^{(x)} := \prodscal{\dt_x}{-\rt_x}$ and $\gap^{(y)} := \prodscal{\dt_y}{-\rt_y}$. We define $\epsilon_t := \prodscal{\dt_x}{\nabla_x\L(\xt, \y^*)- \gnx}$ to be the sequence representing the error between the gradient used for the minimization and the gradient at the point $(\xt, \y^*)$. Then,
       \begin{equation}\label{eq:un}
          w_\gamma^{(x)}   \leq   \wt^{(x)} - \gamma \gap^{(x)} +\gamma \epsilon_t + \gamma^2 \frac{C_\L}{2}.
      \end{equation}
    Now, as $M_{XY}$ is finite (under Lipschitz gradient assumption), we can use the definition of the bilinearity constant~\eqref{def:M} to get
      \begin{equation}\label{eq:grad}
      |\epsilon_t| = \left|\prodscal{\dt_x}{\nabla_x\L(\xt, \y^*)- \gnx}\right| \leq \sqrt{\wt^{(y)}} M_{XY}.
       \end{equation}
    Combining equations \eqref{eq:un} and \eqref{eq:grad} we finally obtain
      \begin{equation}
        w_\gamma^{(x)} \leq \wt^{(x)} - \gamma \gap^{(x)} + \gamma M_{XY} \sqrt{\wt^{(y)}}
          + \gamma^2\frac{ C_\L}{2}.
      \end{equation}
     We can get an analogous inequality for $w_\gamma^{(y)}$,
       \begin{equation}
        w_\gamma^{(y)} \leq \wt^{(y)} - \gamma \gap^{(y)} + \gamma M_{YX} \sqrt{ \wt^{(x)}}
          + \gamma^2 \frac{C_\L}{2}.
      \end{equation}
    Then adding $w_\gamma^{(x)}$ and $w_\gamma^{(y)}$ and using $\sqrt{a} + \sqrt{b} \leq \sqrt{2(a+b)}$ (coming from the concavity of $\sqrt{\cdot}$), we get
      \begin{equation} \label{eq:wgammaInt}
      w_\gamma\leq  \wt - \gamma \gap +\gamma M_{\L} \sqrt{2\wt}+ 2\gamma^2 \frac{C_\L}{2}.
      \end{equation}
    We stress that the above inequality~\eqref{eq:wgammaInt} is valid for any direction $\dt$, using $g_t := \prodscal{\dt}{-\rt}$, and for any feasible step size $\stepsize$ such that $\z_\stepsize \in \X \times \Y$ (the last condition was used in the definition of $C_f$; see also footnote~\ref{foot:Cfcomment} for more information).
      
  To finish the argument, we now use the specific property of the direction $\dt$ and use the crucial Lemma~\ref{lemme:gap} that relates~$w_t$ with the square of the appropriate gap.
  
  For the case~\eqref{case:1} of interior saddle point, we consider $\dt = \dt_{\FW}$ and thus $g_t = \gap^{\FW}$. Then combining Lemma~\ref{lemme:gap} (using the interior strong convexity constant) with~\eqref{eq:wgammaInt}, we get
      \begin{equation}\label{eq:int_strong_usual_gap}w_\gamma\leq  \wt - \gamma \left( 1 - \frac{M_\L}{\sqrt{\muIntL}}\right)\gap^{\FW}+ \gamma^2 C_\L.
      \end{equation}
      
    For the case~\eqref{case:2} of polytope domains, we consider $\dt$ as defined by the SP-AFW algorithm. We thus have $\gap \geq \frac{1}{2} \gap^{\PW}$ by Lemma~\ref{lemma:gap_FW_PW}. Then combining Lemma~\ref{lemme:gap} (using the geometric strong convexity constant) with~\eqref{eq:wgammaInt}, we get
      \begin{equation}\label{eq:geom_strong_pfw_gap}
      w_\gamma\leq  \wt - \gamma \left( \frac{1}{2} - \frac{M_\L}{\sqrt{\mu^{\away}_\L}}\right)\gap^{\PW}+ \gamma^2 C_\L ,
      \end{equation}
      which finishes the proof. Still in the case~\eqref{case:2}, we also present an inequality in terms of the direction gap $\gap$ (which yields a better constant that will be important for the sublinear convergence proof in Theorem~\ref{thm:conv_sublin}) by using instead the inequality $\gap^{\PW} \geq \gap$ (Lemma~\ref{lemma:gap_FW_PW}) with Lemma~\ref{lemme:gap} and~\eqref{eq:wgammaInt}:
      \begin{equation}\label{eq:geom_strong_gt_gap}
      w_\gamma\leq  \wt - \gamma \left( 1 - \frac{M_\L}{\sqrt{\mu^{\away}_\L}}\right)\gap + \gamma^2 C_\L.
      \end{equation}
  \endproof

The above lemma uses a specific update direction~$\dt$ to get a potential new suboptimality~$w_\stepsize$. By using the property that $w_\stepsize \geq 0$ always, we can actually derive an upper bound on the gap in terms of $w_t$ \emph{irrespective of any algorithm} (i.e. this relationship holds for any possible feasible point $(\xt,\yt)$). More precisely, for SP-FW algorithm (case~\eqref{enu:situation1}) the only thing we need to set is a feasible point $\zt$ but for the SP-AFW algorithm (case~\eqref{enu:situation2}) we also need an active set expansion for $\zt$ for which the maximum away step size is larger than $\frac{\CondNumb\gap}{2C_\L}$ (which can potentially not be the active set calculated by an algorithm). 
This is stated in the following theorem, which is a saddle point generalization of the gap upper bound given in Theorem~2 of~\citep{lacoste2015global}. 

\begin{theorem}[Bounding the gap with the second suboptimality] \label{thm:subopt_gap}
  If $\L$ is strongly convex-concave and has a finite curvature constant then
  \begin{itemize}
    \item case~\eqref{enu:situation1}: For any $\zt \in \M$,
      \begin{equation}    
        \gap^{\FW} \leq \frac{2}{\CondNumb^{\FW}} \max \left\{\sqrt{C_\L w_t},w_t \right\}.
      \end{equation}
      Since $\zt$ is fixed, this statement is algorithm free.
    \item case~\eqref{enu:situation2}:  For any $\zt \in \M$, if there exists  an active set expansion for $\zt$ for which $\stepmax = 1$ or $\stepmax \geq \frac{\CondNumb^{\PW}\gap^{\PW}}{2C_\L}$ (see~\eqref{eq:maxstep} for the definition of $\stepmax$) then,
      \begin{equation}    
        \gap^{\PW} \leq \frac{2}{\CondNumb^{\PW}} \max \left\{\sqrt{C_\L w_t},w_t \right\}.
     \end{equation}
  \end{itemize}
  Both statement are algorithm free but $\gap^{\PW}$ depends on a chosen expansion of $\zt$:
   \begin{equation}\label{eq:expansion}
    \zt = \sum_{(\vv_x,\vv_y) \in \mathcal S^{(t)}} \alpha_{\vv_x}^{(t)}\alpha_{\vv_y}^{(t)} (\vv_x,\vv_y)
    \quad \text{where} \quad
     \mathcal S^{(t)} := \left\{ (\vv_x,\vv_y) \in \A\times\B\;;\; \alpha_{\vv_x}^{(t)}\alpha_{\vv_y}^{(t)} >0 \right\},
  \end{equation}
    because,
   \begin{equation}
   \begin{aligned}
   &\gap^{\PW} := \innerProdCompressed{-\rt}{\dt_\FW + \dt_\away} \\[1mm] 
    &\text{where} \quad \dt_\away := \zt - \argmax_{\vv \in \Coreset_x \times \Coreset_y} \innerProdCompressed{\rt}{\vv}, \\ 
    &\text{and} \quad \dt_\FW := \argmin_{\vv \in \A \times \B} \innerProdCompressed{\rt}{\vv} - \zt.
    \end{aligned} 
  \end{equation}
  The maximum step size associated with the active set expansion described in Equation~\eqref{eq:expansion} is 
  \begin{equation} \label{eq:maxstep}
     \gamma_{\max} := \left\{ \begin{array}
     {ll}
     1 \quad \text{if} \quad \innerProd{-\rt}{\dt_\away} \leq \innerProd{-\rt}{\dt_\FW},&\\
     \min\left\{\frac{\alpha^{(t)}_{\vt_x}}{1 -\alpha^{(t)}_{\vt_x}},\frac{\alpha^{(t)}_{\vt_y}}{1 -\alpha^{(t)}_{\vt_y}} \right\} \quad\text{otherwise}.&
     \end{array}\right.
  \end{equation}
  \end{theorem}
  \proof In this proof, we let $(\gap, \CondNumb)$ to stand respectively for $(\gap^{\FW},\CondNumb^{\FW})$ for case~\eqref{case:1} or $(\gap^{\PW},\CondNumb^{\PW})$ for case~\eqref{case:2}.
    We will start from the inequalities~\eqref{eq:wInequalityCaseI} and~\eqref{eq:wInequalityCaseP} in Lemma~\ref{lemme:ineg}. Equation~\eqref{eq:wInequalityCaseI} is valid considering a FW direction $\dt_\FW$, Equation~\eqref{eq:wInequalityCaseP} is valid if we consider the direction that would have be set by the SP-AFW algorithm if it was run at point $\zt$ with the active set expansion described in the theorem statement.
    Since $w_\gamma \geq 0$,  for both cases become :
    \begin{equation}
     0 \leq \wt - \gamma \CondNumb\gap+ \gamma^2 C_\L,
    \end{equation}
    then we can put the gap on the LHS,
    \begin{equation} \label{eq:int_subopt_gap}
    \stepsize\CondNumb\gap - \stepsize^2 C_\L \leq \wt.
    \end{equation}
    This inequality is valid for any $\stepsize \in [0,\stepmax]$.
    In order to get the tightest bound between the gap and the suboptimality, we will maximize the LHS. It can be maximized with $\bar{\gamma} := \frac{\CondNumb\gap}{2C_\L}= \stepsize_t$.
    Now we have two cases: 

    If $\bar{\gamma} \in [0,\stepmax]$, then we get: $\: \CondNumb\gap \leq 2\sqrt{C_\L \, w_t}.$

    And if $\stepmax = 1$ and $\bar{\gamma} = \frac{\CondNumb\gap}{2C_\L} > 1$, then setting $\stepsize=1$ we get: $\: \CondNumb\gap \leq 2 w_t$. By taking the maximum between the two options, we get the theorem statement.
  \endproof
  The previous theorem guarantees that the gap gets small when $w_t$ gets small (only for the non-drop steps if in situation~\eqref{case:2}).
  As we use the gap as a stopping criterion in the algorithm, this is a useful theorem to provide an upper bound on the number of iterations needed to get a certificate of suboptimality.
  
  The following corollary provides a better bound on $h_t$ than the inequality $h_t \leq \textrm{cst} \sqrt{w_t}$ previously shown \eqref{eq:relate_h_w} when the situation is \eqref{enu:situation2}. It will be useful later to get a better rate of convergence for $h_t$ under hypothesis~\eqref{case:2}. 
  \begin{corollary}[Tighter bound on $h_t$ for non-drop steps in situation~\eqref{case:2}]
  \label{cor:h_w}
  Suppose that $\L$ is strongly convex-concave and has a finite curvature constant, and that the domain is a product of polytopes (i.e. we are in situation~\eqref{case:2}). Let $\zt \in \M$  be given. If there exists an active set expansion for $\zt$ for which the maximum step size is larger than $\frac{\CondNumb\gap}{2C_\L}$ (see Theorem~\eqref{thm:subopt_gap} for more details) ,
    \begin{equation}
            h_t \leq \frac{2 \max \{C_\L,w_t\}}{\mu^{\away}_\L (\CondNumb^\PW)^2} w_t ,
    \end{equation}
    where $\CondNumb^\PW$ is defined in~\eqref{eq:wInequalityCaseP}.
  \end{corollary}
  \proof
  By Lemma~\ref{lemme:gap} in situation~\eqref{case:2}, we have:
  \begin{equation*}
  h_t \leq \frac{(\gap^{\PW})^2}{2 \mu^{\away}_\L} \leq \frac{2\max \{C_\L,w_t\}}{\mu^{\away}_\L (\CondNumb^\PW)^2} w_t.
  \end{equation*}
  The last inequality is obtained by applying the upper bound on the gap given in the previous Theorem~\ref{thm:subopt_gap}.
  \endproof

\subsection{Proof of Theorem~\ref{thm:conv}}
  \label{subsec:proof_main_thm}
  In this section, we will prove that under some conditions on the constant defined in subsections~\ref{sub:curv} and~\ref{sub:affine_invariant_measures_of_strong_convexity}, the suboptimalities $\wt$
  vanish linearly with the adaptive step size $\gamma_t = \min \left\{\stepmax,\frac{\CondNumb}{2C_\L}\gap\right\}$ or sublinearly with the universal step size $\gamma_t = \min \left\{\stepmax,\frac{2}{2+k(t)}\right\}$.
  \begin{lemma}[Geometric decrease of second suboptimality] 
    \label{lemme:Lin} Let $\L$ be a 
    strongly convex-concave function with a smoothness constant $C_\L$, a positive interior strong convex-concavity constant $\muIntL$~\eqref{def:interior_strongL} or a positive geometric strong convex-concavity $\mu^{\away}_\L$~\eqref{def:geom_strongL}. Let us also define the rate multipliers~$\CondNumb$ as
  \begin{equation}
    \CondNumb^{\FW} :=1- \frac{M_\L }{\sqrt{\muIntL}} \quad\text{and}\quad \CondNumb^{\PW} := \frac{1}{2}- \frac{M_\L }{\sqrt{\mu^{\away}_\L}} \quad(\text{see Equation~\eqref{eq:def_bilin_coef} for the definition of } M_\L).
  \end{equation}
    Let the tuple $(\gap,\CondNumb,\mu_\L)$ refers to
  either $(\gap^{\FW},\CondNumb^{\FW},\muIntL)$ for case~\eqref{case:1} where the algorithm is SP-FW , or $(\gap^{\PW},\CondNumb^{\PW},\mu^{\away}_\L)$ for case~\eqref{case:2} %
  where the algorithm is SP-AFW,

        If $\CondNumb > 0$, then  at each non-drop step (when $\gamma_t < \gamma_{\max}$ or $ \gamma_{\max} \geq 1$ ),
        the suboptimality $w_t$ of the algorithm with
        step size $\stepsize_t = \min(\stepmax, \frac{\CondNumb}{2C_\L}\gap)$ decreases
        geometrically as
          \begin{equation} 
             \wtt \leq (1- \rho_\L) \wt
          \end{equation}
        where $\rho_\L := \frac{\CondNumb^2}{2} \frac{\mu_\L}{C_\L}$.
        Moreover, for case \eqref{enu:situation1} there is no drop step and for case~\eqref{enu:situation2} the number of drop step (when $\gamma_t = \gamma_{\max}$) is upper bounded by two third of the number of iteration (see Section~\ref{par:drop_steps}, Equation~\eqref{eq:upper_bound_drop_steps}), while when we have a drop step, we still have:
        \begin{equation}
         \wtt \leq \wt.
        \end{equation}
  \end{lemma}

  \proof 
  The bulk of the proof is of a similar form for both SP-FW and SP-AFW, and so in the following, we let $(\gap, \mu_\L, \CondNumb)$ to stand respectively for $(\gap^{\FW},\muIntL,\CondNumb^{\FW})$ for SP-FW (case~\eqref{case:1}) or $(\gap^{\PW},\mu^{\away}_\L,\CondNumb^{\PW})$ for SP-AFW (case~\eqref{case:2}). As $\stepsize_t \leq \stepmax$, we can apply the important Lemma~\ref{lemme:ineg} with $\stepsize = \stepsize_t$ (the actual step size that was taken in the algorithm) to get:
  \begin{equation} \label{eq:wtIneq}
  w_{t+1} = w_{\stepsize_t}  
         \leq \wt  
        - \CondNumb \stepsize_t \gap
        + \gamma_t^2 {C_\L} . 
  \end{equation}
  We note in passing that the adaptive step size rule $\stepsize_t = \min(\stepmax, \frac{\CondNumb}{2 C_\L}\gap)$ was specifically chosen to minimize the RHS of~\eqref{eq:wtIneq} among the feasible step sizes. 
  
  If $\frac\CondNumb{2C_\L}\gap\leq \stepmax$, then we have $\gamma_t = \frac\CondNumb{2C_\L}\gap$ and so~\eqref{eq:wtIneq} becomes:
      \begin{equation}
        \wtt \leq \wt - \frac{\CondNumb^2}{ 2C_\L}  (\gap)^2+ \frac{\CondNumb^2}{ 4C_\L}(\gap)^2 = \wt - \frac{\CondNumb^2}{4C_\L} (\gap)^2.
      \end{equation}
    Applying the fact that the square of the appropriate gap upper bounds $w_t$ (Lemma~\ref{lemme:gap} with a similar form for both cases~\eqref{case:1} and~\eqref{case:2}), we directly obtain the claimed geometric decrease
      \begin{equation}
        \wtt \leq \wt \left(1- \frac{\CondNumb^2}{2} \frac{\mu_\L}{C_\L} \right).
      \end{equation}
    If $ \frac\CondNumb{2 C_\L}\gap> \stepmax$, then we have $\gamma_t = \stepmax$ and so~\eqref{eq:wtIneq} becomes:
    \begin{align}
      \wtt 
      &\leq  \wt - \CondNumb \stepmax \gap + \stepmax^2 C_\L & \notag\\
      & \leq \wt - \CondNumb \stepmax \gap + \frac\CondNumb{2} \stepmax \gap &( \text{using}\quad  C_\L < \frac{\CondNumb}{2 \stepmax}\gap) \\
      & \leq \wt - \frac\CondNumb{2} \stepmax \gap &\notag\\
       & \leq \wt \left(1-\frac{\CondNumb}{2}\stepmax\right).& (w_t \leq g_t \text{ by Lemma \ref{lemme:g}})
    \end{align}
    If $\stepmax \geq 1$ (either we are taking a FW step or an away step with a big step size), then the geometric rate is at least $(1-\frac{\CondNumb}{2})$, which is a better rate than $\rho_\L$ since $\CondNumb^2 \leq \CondNumb$ as $\CondNumb \leq 1$, and one can show that $\frac{\mu_\L}{C_\L} \leq 1$ always (see Remark~7 in Appendix D of~\citep{lacoste2015global} for case~\eqref{enu:situation2} and use a similar argument for case~\eqref{enu:situation1}). Thus $\rho_\L$ is valid both when $\stepsize_t < \stepmax$ or $\stepmax \geq 1$, as claimed in the theorem.
    
    When $\stepsize_t = \stepmax < 1$, we cannot guarantee sufficient progress as $\stepmax$ could be arbitrarily small (this can only happen for an away step as $\stepmax = 1$ for a FW step). These are the problematic \emph{drop steps}, but as explained in Appendix~\ref{par:drop_steps} with Equation~\eqref{eq:upper_bound_drop_steps}, they cannot happen too often for SP-AFW.
    
    Finally, to show that the suboptimality cannot increase during a drop step ($\gamma_t= \gamma_{\max}$), we point out that 
    the function $\gamma \mapsto \wt - \gamma \CondNumb^{\PW}\gap^{\PW}+ \gamma^2 C_\L$ is a convex function that is minimized by $\bar{\gamma} = \frac{\CondNumb^{\PW}}{2C_\L} \gap^{\PW}$ and so is decreasing on $[0, \bar{\stepsize}]$. When $\stepsize_t = \stepmax$, we have that $\stepmax \leq \bar{\gamma}$, and thus the value for $\gamma = \gamma_{\max}$ is lower than the value for $\gamma = 0$, i.e.
      \begin{equation}
            \wtt \leq \wt - \stepmax \CondNumb^{\PW}\gap^{\PW}+ \gamma_{\max}^2 C_\L \leq \wt.
      \end{equation}
  \endproof
    The previous lemma (Lemma~\ref{lemme:Lin}), the fact that the gap upper bounds the suboptimality (Lemma~\ref{lemme:g}) and the primal suboptimalities analysis lead us directly
    to the following theorem. This theorem is the affine invariant formulation with adaptive step size of Theorem~\ref{thm:conv}.

  \begin{theorem}\label{thm:conv_affine_invariant}
  Let $\L$ be a 
    strongly convex-concave function with a finite smoothness constant $C_\L$, a positive interior strong convex-concavity constant $\muIntL$~\eqref{def:interior_strongL} or a positive geometric strong convex-concavity $\mu^{\away}_\L$~\eqref{def:geom_strongL}. Let us also define the rate multipliers~$\CondNumb$ as
  \begin{equation}\label{eq:thm24_state_multiplier}
    \CondNumb^{\FW} :=1- \frac{M_\L }{\sqrt{\muIntL}} \quad\text{and}\quad \CondNumb^{\PW} := \frac{1}{2}- \frac{M_\L }{\sqrt{\mu^{\away}_\L}} \quad(\text{see Equation~\eqref{eq:def_bilin_coef} for the definition of } M_\L).
  \end{equation}
    Let the tuple $(\gap,\CondNumb,\mu_\L)$ refers to
  either $(\gap^{\FW},\CondNumb^{\FW},\muIntL)$ for case~\eqref{case:1} where the algorithm is SP-FW , or $(\gap^{\PW},\CondNumb^{\PW},\mu^{\away}_\L)$ for case~\eqref{case:2} where the algorithm is SP-AFW,

  If $\CondNumb>0$, then the suboptimality $h_t$ of the iterates of the algorithm with
  step size $\stepsize_t = \min(\stepmax, \frac{\CondNumb}{2C_\L}\gap)$ decreases
  geometrically\footnote{For a non-drop step one can use Corollary~\ref{cor:h_w} to get the better rate on $h_t$ losing the square root but with a potentially worse constant $ h_t \leq \frac{2\max \{C_\L,w_0(1-\rho_\L)^{k(t)}\}}{\mu^{\away}_\L (\CondNumb^\PW)^2} w_0(1-\rho_\L)^{k(t)}$.} as,
        \begin{equation} 
          h_t \leq P_\L \sqrt{2w_0}(1- \rho_\L)^{k(t)/2}
        \end{equation}
  where $\rho_\L := \frac{\CondNumb^2\mu_\L }{2C_\L}$ and $k(t)$ is the number of non-drop step after $t$ steps. For SP-FW, $k(t)=t$ and for SP-AFW, $k(t)\geq t/3$. Moreover we can also upper bound the minimum gap observed, for all $T \in \N$
  \begin{equation}
    \min_{t \leq T} \gap \leq \frac{2\max\{\sqrt{C_\L},\sqrt{w_0} (1-\rho_\L)^{k(T)/2}\}}{\CondNumb{}} \sqrt{w_0} \left( 1- \rho_\L \right)^{k(T)/2}.
  \end{equation}
  The Theorem~\ref{thm:conv} statement can be deduced from this theorem using the lower and upper bounds on the affine invariant constant of this statement. More precisely, one can upper bound $C_\L$, $M_\L$, $P_\L$ respectively with Propositions~\ref{prop:Cbounded}, \ref{Mbounded} and \ref{prop:upper_bound_P_L} and lower bound $\muIntL$ and $\mu^{\away}_\L$ respectively with Proposition~\ref{prop:deltaL} and Equation~\eqref{eq:muAwayLInequality}.\footnote{Note that only the definition of $\delta_\mu$ in~Theorem~\ref{thm:conv} for case~\eqref{case:2} requires to use the Euclidean norm (because inequality~\eqref{eq:muAwayLInequality} with the pyramidal width only holds for the Euclidean norm). On the other hand, any norm could be used for (separately) bounding $C_\L$, $M_\L$ and $P_\L$.} 
  If we apply these bounds to the rate multipliers in~\eqref{eq:thm24_state_multiplier}, it gives the smaller rate multipliers $\CondNumb$ stated in Theorem~\ref{thm:conv}. 
  \end{theorem}
  \proof
  We uses the Lemma~\ref{lemme:Lin} giving a geometric scheme, with a straightforward recurrence we prove that,
    \begin{equation}
      w_t \leq w_0 (1-\rho_\L)^{k(t)},
    \end{equation}
    where $k(t)$ is the number of non-drop step steps. This number is equal to $t$ for the SP-FW algorithm and it is lower bounded by $t/3$ for the SP-AFW algorithm (see Section~\ref{par:drop_steps} Equation~\eqref{eq:upper_bound_drop_steps}).
    Then by using Proposition \eqref{prop:relation_primal}  relating $h_t$ and the square root of $w_t$ we get the first statement of the theorem,
      \begin{equation}
       h_t \leq   P_\L \sqrt{2w_0}(1- \rho_\L)^{k(t)/2}.
      \end{equation}
    To prove the second statement of the theorem we just use Theorem~\ref{thm:subopt_gap} for the last \emph{non-drop step} after $T$ iterations (let us assume it was at step $t_0$),
    \begin{align}
     g_{t_0} &\leq   \frac{2 \max\{\sqrt{C_\L},\sqrt{w_0} (1-\rho_\L)^{k(t_0)/2}\}}{\CondNumb{}} \sqrt{w_0}(1- \rho_\L)^{k(t_0)/2} &\\
          &=\frac{2\max\{\sqrt{C_\L},\sqrt{w_0} (1-\rho_\L)^{k(T)/2}\}}{\CondNumb{}} \sqrt{w_0}(1- \rho_\L)^{k(T)/2}. & (\text{because } k(t_0)= k(T))
    \end{align}
    The minimum of the gaps observed is smaller than the gap at time $t_0$ then,
    \begin{equation}
      \min_{t\leq T}\gap \leq g_{t_0} \leq \frac{2\max\{\sqrt{C_\L},\sqrt{w_0} (1-\rho_\L)^{k(T)/2}\}}{\CondNumb{}} \sqrt{w_0}(1- \rho_\L)^{k(T)/2}.
    \end{equation}
  \endproof
  The affine invariant formulation with the universal step size $\stepsize_t = \min \left\{\stepmax,\frac{2}{2+k(t)}\right\}$ of Theorem~\ref{thm:conv} also follows from Lemma~\ref{lemme:Lin} by re-using standard FW proof patterns.
   \begin{theorem}\label{thm:conv_sublin}
   Let $\L$ be a 
    strongly convex-concave function with a finite smoothness constant $C_\L$, a positive interior strong convex-concavity constant $\muIntL$~\eqref{def:interior_strongL} or a positive geometric strong convex-concavity $\mu^{\away}_\L$~\eqref{def:geom_strongL}. Let us also define the rate multipliers~$\CondNumb$ as
  \begin{equation}
    \CondNumb^{\FW} :=1- \frac{M_\L }{\sqrt{\muIntL}} \quad\text{and}\quad \tilde \CondNumb^{\PW} := 1- \frac{
        M_\L}{\sqrt{\mu^{\away}_\L}} \quad(\text{see Equation~\eqref{eq:def_bilin_coef} for the definition of } M_\L).
  \end{equation} 
    Let $\CondNumb$ refers to
  either $\CondNumb^{\FW}$ for case~\eqref{case:1} where the algorithm is SP-FW, or $\tilde\CondNumb^{\PW}$ for case~\eqref{case:2} where the algorithm is SP-AFW,

  If $\CondNumb>\frac{1}{2}$, then the suboptimality $w_t$ of the iterates of the algorithm with universal step size $\stepsize_t = \min\left\{\stepmax,\frac{2}{2+k(t)}\right\}$ (see Equation~\eqref{eq:update} for more details about $\stepmax$) has the following decreasing upper bound:
        \begin{equation}\label{eq:sub_rate_w}
          w_t \leq \frac{C}{{2+k(t)}}
        \end{equation}
      where $C = 2 \max\left(w_0,\frac{2 C_\L}{2 \CondNumb-1}\right)$ and $k(t)$ is the number of non-drop step after $t$ steps. For SP-FW, $k(t)=t$ and for SP-AFW, $k(t)\geq t/3$. 
      Moreover we can also upper bound the minimum FW gap observed for $T\geq 1$,
  \begin{equation}
    \min_{t \leq T} \gap^\FW \leq \frac{5C}{\CondNumb{}(k(T)+1)}.
  \end{equation}
    Note that in this theorem the constant $\tilde \CondNumb^{\PW}$ is slightly different from the constant $\CondNumb^{\PW}$ in Theorem~\ref{thm:conv_affine_invariant}.
  \end{theorem}
  \proof
    We can put both the recurrence~\eqref{eq:int_strong_usual_gap} for the SP-FW algorithm and the recurrence~\eqref{eq:geom_strong_gt_gap} for the SP-AFW algorithm (from the proof of Lemma~\ref{lemme:ineg}) in the following form by using our unified notation introduced in the theorem statement:
          \begin{equation}\label{eq:proof_sublinear2}
         \wtt 
       \leq \wt  
      - \stepsize_t \CondNumb \gap
      + \stepsize_t^2 {C_\L}. 
      \end{equation}
    Note that the gap $\gap$ is the one defined in Equation~\eqref{eq:gaps} and depends on the algorithm.
    Let $(\CondNumb)$ to stand respectively for $(\CondNumb^{\FW})$ for SP-FW (case~\eqref{case:1}) or $(\tilde \CondNumb^{\PW})$ for SP-AFW (case~\eqref{case:2}). With this notation, the inequality $\gap \geq w_t$ leads to,
      \begin{equation}\label{eq:def_f}
          \wtt 
       \leq \wt \left( 1- \CondNumb\gamma_t \right)   
      + \gamma_t^2 {C_\L}. 
      \end{equation}
      Our goal is to show by induction that 
      \begin{equation}\label{eq:induction_sublin}\tag{$\star$}
      w_t \leq \frac{C}{2+k(t)} \quad \text{ where } C := 2\max \left( w_0, \frac{2C_\L}{2 \CondNumb -1} \right).
      \end{equation}
      Let us first define the convex function $f_t :\gamma \mapsto  \wt \left( 1- \CondNumb\gamma \right)
      + \gamma^2 {C_\L}$. 
      We will show that under~\eqref{eq:induction_sublin}, the function $f_t$ has the following property:
      \begin{equation}
        f_t\left(\frac{2}{2+k(t)}\right) \leq \frac{C}{3+k(t)} .
      \end{equation}
      This property is due to a simple inequality on integers; let $k = k(t)$, from the crucial induction assumption, we get:
      \begin{equation}
        f_t\left(\frac{2}{2+k}\right) 
       = 
          w_{t} \frac{2+k - 2 \CondNumb}{2+k}+ \frac{4}{(2+k)^2} {C_\L} 
       \leq 
          \frac{C}{3+k} \left[ \frac{(3+k)(k+1 - (2 \CondNumb - 1) +\frac{4C_\L}{C})}{(2+k)^2} \right], 
      \end{equation}
      but $(2 \CondNumb - 1) \geq \frac{4C_\L}{C}$ and $(3+k)(1+k)<(2+k)(2+k)$ for any $k$, thus 
      \begin{equation}\label{eq:sublin_f}
       f_t\left(\frac{2}{2+k}\right)  \leq \frac{C}{3+k}.
      \end{equation}
      Equation~\eqref{eq:sublin_f} is crucial for the inductive step of our recurrence.
      \begin{itemize}
        \item 
      Hypothesis~\eqref{eq:induction_sublin} is true for $t=0$ because $k(0)=0$. 
        \item 
      Now let us assume that \eqref{eq:induction_sublin} is true for a $t \in \N$.
      We set the stepsize $\stepsize_t := \min \left\{\stepmax,\frac{2}{2+k(t)}\right\}$.

      If $k(t+1)=k(t)+1$, it means that $\stepsize_t = \frac{2}{2+k(t)}$ and then by~\eqref{eq:def_f} and~\eqref{eq:sublin_f},
      \begin{equation} \label{eq:proof_w_sublin}
         w_{t+1} 
         \leq 
          f_t\left(\frac{2}{2+k(t)}\right) 
         \leq \frac{C}{3+k(t)}
         =  \frac{C}{2+k(t+1)}.
      \end{equation}

      If $k(t+1)= k(t)$, then it means that $0 \leq \stepsize_t < \frac{2}{2+k(t)}$. Hence, the convexity of the function $f_t$ leads us to the inequality
      \begin{align} \notag 
       w_{t+1} \leq f_t(\stepsize_t) 
       &\leq \max \left\{ f_t(0),f_t\left(\frac{2}{2+k(t)}\right) \right\} \\
       & = \max \left\{ w_t,f_t\left(\frac{2}{2+k(t)}\right) \right\}\label{eq:ineg_wt_wtt} \\
       &\leq \max \left\{ \frac{C}{2+k(t)},\frac{C}{3+k(t)} \right\} \\
       &\leq \frac{C}{2+k(t)}.
       \end{align} 
      where we used~\eqref{eq:sublin_f} and the induction hypothesis~\eqref{eq:induction_sublin} to get the penultimate inequality~\eqref{eq:ineg_wt_wtt}. Since we assumed that $k(t+1)= k(t)$, we get
      \begin{equation}
        w_{t+1} \leq  \frac{C}{2+k(t+1)}, 
      \end{equation}
      completing the induction proof for~\eqref{eq:sub_rate_w}.
      \end{itemize}

      In case \eqref{case:1}, $k(t) = t$ and in case \eqref{case:2}, $k(t) \geq t/3$ (see Equation~\eqref{eq:upper_bound_drop_steps}), leading us to the first statement of our theorem.

The proof of the second statement is inspired by the proof of Theorem C.3 from~\citep{lacoste2013block}.

With the same notation as the proof of Lemma~\ref{lemme:ineg}, we start from Equation~\eqref{eq:proof_sublinear2} where we isolated the gap $\gap$ to get the crucial inequality
  \begin{equation}
   \gap
         \leq \frac{\wt  -w_{t+1}}{\CondNumb \stepsize_t }
        + \gamma_t \frac{C_\L}{\CondNumb}. 
  \end{equation}
  Since the gap $\gap$ is the one depending on the algorithm  defined by $\gap := \innerProd{-\rt}{\dt}$, we have $\gap = \gap^\FW$ for SP-FW and $\gap = \max \left( \gap^\FW, \gap^\away \right)  \geq \gap^\FW$ for SP-AFW. Thus,
     \begin{equation} \label{eq:Ineg_g_wt_wtt}
    \gap^\FW \leq \gap
         \leq \frac{\wt  -w_{t+1}}{\CondNumb \stepsize_t }
        + \gamma_t \frac{C_\L}{\CondNumb}. 
  \end{equation}
  In the following in order not to be too heavy with notation we will work with de FW gap and note $\gap$ for $\gap^\FW$. 

  The proof idea is to take a convex combination of the inequality~\eqref{eq:Ineg_g_wt_wtt} to obtain a new upper-bound on a convex combination of the gaps computed from step $0$ to step $T$. Let us introduce the convex combination weight $\rho_t := \frac{\stepsize_t\cdot  k(t)( k(t) +2) }{S_T}$ where $k(t)$ is the number of non-drop steps after $t$ steps and $S_T$ is the normalization factor. Let us  also call $N_T := \{ t \leq T \, | \, t \text{ is a non-drop step}\}$.  Taking the convex combination of~\eqref{eq:Ineg_g_wt_wtt}, we get
  \begin{equation} \label{eq:Ineg_mean_g_wt_wtt}
        \sum_{t=0}^T \rho_t \gap \leq 
        \sum_{t=0}^T \rho_t \frac{\wt  -w_{t+1}}{\CondNumb \stepsize_t }+ 
        \sum_{t=0}^T \rho_t \gamma_t \frac{C_\L}{\CondNumb}. 
  \end{equation}
  By regrouping the terms and ignoring the negative term, we get
  \begin{equation} \label{eq:sum_gt_wt_rhot}
        \sum_{t=0}^T \rho_t \gap \leq 
        \frac{w_0 \rho_0}{\CondNumb \stepsize_0} +
        \frac{1}{\CondNumb}\sum_{t=0}^{T-1} w_{t+1}\left( \frac{\rho_{t+1}}{ \stepsize_{t+1} }   -\frac{\rho_t}{ \stepsize_t}\right)+ 
        \sum_{t=0}^{T} \rho_t \gamma_t \frac{C_\L}{\CondNumb}. 
  \end{equation}
  By definition $\frac{\rho_t}{\stepsize_t} := \frac{k(t)(k(t)+2)}{S_T}$ and notice that $\rho_0 = 0$. We now consider two possibilities: if $\gamma_t$ is a drop step, then $k(t+1) = k(t)$ and so
  \begin{equation}
   \frac{\rho_{t+1}}{ \stepsize_{t+1}} - \frac{\rho_{t}}{ \stepsize_{t}} = 0.
  \end{equation}
  If $\stepsize_t$ is a non-drop step, then $k(t+1) = k(t) + 1$ and thus we have
  \begin{align}
     \frac{\rho_{t+1}}{ \stepsize_{t+1}} - \frac{\rho_{t}}{ \stepsize_{t}} 
     &=   \frac{(k(t)+1)(k(t)+3)}{S_T} -  \frac{k(t)(k(t)+2)}{S_T} \\
     &= \frac{2k(t)+3}{S_T}. \label{eq:rho/gamma} 
  \end{align}
  As $\stepsize_t \leq \frac{2}{k(t)+2}$, we also have $\rho_t \stepsize_t \leq \frac{4 k(t)}{S_T (k(t)+2))} $.
  The normalization factor $S_T$ to define a convex combination is equal to 
  \begin{equation}
    S_T := \sum_{u=0}^T \stepsize_u \cdot k(u)( k(u) +2)  
    \geq  \sum_{\small\substack{u=0 \\ u \in N_{T} }}^T  \frac{2}{2 +k(u)} \cdot k(u)( k(u) +2) 
    = \sum_{k = 0}^{k(T)} 2 k
    = k(T)(K(T)+1).  
  \end{equation}
  Plugging this and~\eqref{eq:rho/gamma} in the inequality~\eqref{eq:sum_gt_wt_rhot} with the rate~\eqref{eq:induction_sublin} shown by induction gives us,
  \begin{align}
       \sum_{t=0}^T \rho_t \gap 
        &\leq 
        0 +
        \frac{1}{\CondNumb} \sum_{\small\substack{t=0 \\ t \in N_{T} }}^{T-1} \frac{C}{2+k(t+1)}\frac{2k(t)+3}{k(T)(k(T)+1)} + 
         \sum_{t=0}^T \frac{4k(t)}{k(T)(k(T)+1)(k(t)+2)}\frac{C_\L}{\CondNumb} \\
        & \leq \frac{2C}{\CondNumb{}} \frac{1}{k(T)(k(T)+1)} \Big( \sum_{\small\substack{t=0 \\ t \in N_{T} }}^{T-1}  1 + \frac{2C_\L}{C} \sum_{t=1 }^T  1 \Big) \\
        & \leq \frac{2C}{\CondNumb (k(T)+1)} (1 + \frac{\CondNumb}{2}\frac{T}{k(T)}) \leq \frac{5C}{\CondNumb (k(T)+1)},
  \end{align}
  by using $k(T) \geq T/3$.
  Finally, the minimum of the gaps is always smaller than any convex combination, so we can conclude that (for $T\geq1$):
  \begin{equation}
    \min_{0 \leq t \leq T} \gap \leq \frac{5C}{\CondNumb{}(k(T)+1)}.
  \end{equation}
  \endproof

\section{Strongly convex sets}
  \label{sec:strong_conv_proof}
  In this section, we are going to prove that the function $\s(\cdot)$ is Lipschitz continuous when the sets $\X$ and $\Y$ are strongly convex and when the norm of the two gradient components are uniformly lower bounded. We will also give the details of the convergence rate proof for the strongly convex sets situation.
  Our proof uses similar arguments as~\citet[Theorem~3.4 and~3.6]{dunn1979rates}.
  \begin{reptheorem}{thm:s_lip}  
   Let $\X$ and $\Y$ be $\beta$-strongly convex sets. 
   If $\min(\|\nabla_{\!x} L(\z)\|_{\X^*}, \|\nabla_{\!y} L(\z)\|_{\Y^*}) \geq\delta>0$ for all $\z \in \M$, then the oracle function $\z \mapsto \s(\z) := \arg\min_{\s \in \M}\innerProd{\s}{F(\z)}$ is well defined and is $\frac{4L}{\delta \beta }$-Lipschitz continuous (using the norm $\|(\x,\y)\|_{\X \times \Y} := \|\x\|_\X + \|\y\|_\Y$), where $F(\z) := \left(
     \nabla_x \L(\z),
     -\nabla_y \L(\z)
     \right)$.
  \end{reptheorem}
 \proof First note that since the sets are strongly convex, the minimum is reached at a unique point. Then, we introduce the following lemma which can be used to show that each component of the gradient is Lipschitz continuous irrespective of the other set.
    \begin{lemma} \label{lemma:FxLip}
    Let $F_x:\X\times \A \to \R^d$ be a $L$-Lipschitz continuous function (i.e. $\|F_x(\z)-F_x(\z')\|_{\X^*} \leq L \|\z - \z'\|_{\X \times \A}$) and $\X$ a $\beta$-strongly convex set. If $\forall \z \in \X\times \A, \; \|F_x(\z)\|_{\X^*}\geq \delta > 0$, then $\s_x:\z \mapsto \argmin_{\s \in \X} \prodscal{\s}{F_x(\z)}$ is $ \frac{2L}{\delta \beta }$-Lipschitz continuous.
    \end{lemma}
    \proof
    Let $\z, \, \z' \in \X\times \A$ and let $\bar \x = \frac{\s_x(\z)+ \s_x(\z')}{2}$, then
    \begin{align} 
      \prodscal{\s_x(\z)-\s_x(\z')}{-F_x(\z)} 
        & = 2\prodscal{\s_x(\z)-\bar \x}{-F_x(\z)} \notag \\
        & \geq  2\prodscal{\x-\bar \x}{-F_x(\z)}. \quad \forall \x \in \X \quad \text{(by definition of $s_x$)} \label{eq:Sinequality}
    \end{align}
    Now~\eqref{eq:Sinequality} holds for any $\x \in B_\beta\left(\tfrac{1}{2},\s_x(\z),\s_x(\z') \right)$ as this set is included in $\X$ by $\beta$-strong convexity of~$\X$.
    Then since $\bar \x$ is the center of $B_\beta\left(\tfrac{1}{2},\s_x(\z),\s_x(\z') \right)$, we can choose a $\x$ in this ball such that $\x - \bar \x$ is in the direction which achieves the dual norm of $-F_x(\z)$.\footnote{For the Euclidean norm, we choose $\x - \bar \x$ proportional to $-F_x$; but for general norms, it could be a different direction.} More specifically, we have that:
    $$
    \|-F_x(\z)\|_{\X^*} = \sup_{\| \vv \|_{\X} \leq 1} \prodscal{-F_x(\z)}{\vv}.
    $$
    As we are in finite dimensions, this supremum is achieved by some vector $\vv$. So choose $\x := \bar{\x} + \frac{\vv}{\|\vv\|_\X} \frac{\beta}{8} \|\s_x(\z)-\s_x(\z')\|_\X^2 \in B_\beta\left(\tfrac{1}{2},\s_x(\z),\s_x(\z') \right)$ and plug it in~\eqref{eq:Sinequality}:
    
    \begin{align}
      \prodscal{\s_x(\z)-\s_x(\z')}{-F_x(\z)} 
       & \geq \frac{\beta\|\s_x(\z)-\s_x(\z')\|_\X^2}{4 \|\vv\|_\X} \prodscal{\vv}{-F_x(\z)} \\
        & =  \frac{\beta\|\s_x(\z)-\s_x(\z')\|_\X^2}{4 \|\vv\|_\X}  \|F_x(\z)\|_{\X^*} \\
        & \geq \frac{\beta}{4} \|\s_x(\z)-\s_x(\z')\|_\X^2 \|F_x(\z)\|_{\X^*} . \label{eq:s_lip1}
    \end{align}
    Switching $\z$ and $\z'$ and using a similar argument, we get,
    \begin{equation} \label{eq:s_lip2}
         \prodscal{\s_x(\z')-\s_x(\z)}{-F_x(\z')} \geq \frac{\beta}{4} \|\s_x(\z)-\s_x(\z')\|_\X^2 \|F_x(\z')\|_{\X^*}  .
    \end{equation}
    Hence summing \eqref{eq:s_lip1} and \eqref{eq:s_lip2},
    \begin{align}\label{eq:s_lip_final}
       \frac{\beta}{4} (\|F_x(\z)\|_{X^*} +\|F_x(\z')\|_{\X^*}) \, \|\s_x(\z)-\s_x(\z')\|_\X^2 
       &\leq 
       \prodscal{\s_x(\z)-\s_x(\z')}{F_x(\z')-F_x(\z)} \notag\\
       &\leq
       \|\s_x(\z)-\s_x(\z')\|_\X \|F_x(\z')-F_x(\z)\|_{\X^*} \\
       & \leq L \|\s_x(\z)-\s_x(\z')\|_\X \|\z'-\z\|_{\X \times \Y} \;\;\;\text{(Lip. cty. of  $F_x$)}\notag
    \end{align}
    and finally
    \begin{equation}
      \|\s_x(\z)-\s_x(\z')\|_\X \leq \frac{4L}{\beta \left( \|F_x(\z)\|_{\X^*} +\|F_x(\z')\|_{\X^*} \right) } \|\z-\z'\|_{\X \times \Y} \leq \frac{2L}{\delta \beta } \|\z-\z'\|_{\X \times \Y}.
    \end{equation}
    \endproof
    To prove our theorem, we will notice that for the saddle point setup,
    the oracle function $\s(\cdot) := \argmin_{\s \in \M} \prodscal{\s}{F(\cdot)}$ can be decomposed as $\s(\cdot)=(\s_x(\cdot),\s_y(\cdot))$
    where $\s_x(\cdot) := \argmin_{\s \in \X} \prodscal{\s}{F_x(\cdot)}$ 
    and   $\s_y(\cdot) := \argmin_{\s \in \Y} \prodscal{\s}{F_y(\cdot)}$. 
    Then applying our lemma, the function $\s_x(\cdot)$ is Lipschitz continuous. 
    The same way $\s_y(\cdot)$ is Lipschitz continuous. Then, for all $\z,\z'$ in $\M$
    \begin{equation}
      \|\s(\z)-\s(\z')\|_{\X \times \Y} = \|\s_x(\z)-\s_x(\z')\|_\X + \|\s_y(\z)-\s_y(\z')\|_\Y \leq \frac{4L}{\delta\beta} \|\z-\z'\|_{\X \times \Y},
    \end{equation}
    which gives the definition of the Lipschitz continuity of our function and proves the theorem.
    \endproof
    In this theorem,  we introduced the function $F$. This function is monotone in the following sense:
  \begin{equation}\label{eq:nondecreasing}
        \forall \z, \z' \quad \prodscal{\z - \z'}{F(\z)-F(\z')} \geq 0. 
  \end{equation}
  Actually this property follows directly from the convexity of $\L(\cdot,\y)$ and the concavity of $\L(\x,\cdot)$.
  We can also prove that when the sets $\X$ and $\Y$ are strongly convex and when the gradient is uniformly lower bounded, we can relate the gap and the distance between $\zt$ and $\st$. 
   \begin{lemma}\label{lemma:lower_gap}
   If $\X$ is a $\beta$-strongly convex set and if $\|\nabla f\|_{\X^*}$ is uniformly lower bounded by $\delta$ on $\X$, then
  \begin{equation}
    \max_{\s \in \X} \prodscal{\s - \x}{- \nabla f(\x)}\geq \frac{\beta}{4}  \delta \|\s(\x) - \x\|^2,
  \end{equation}
  where $\s(\x) = \argmax_{\s \in \X}\prodscal{\s-\x}{-\nabla f (\x)}$.
  \end{lemma}
  \proof Let $\x$ and $\s(\x)$ be in $\X$. We have $B_\beta \left( \frac{1}{2}, \s(\x), \x \right) \subset \X$ by $\beta$-strong convexity. So as in the proof of Lemma~\ref{lemma:FxLip}, let $\vv$ be the vector such that $\|\vv\|_\X \leq 1$ and $\prodscal{-\nabla f(\x)}{\vv} = \|\nabla f(\x)\|_{\X^*}$. Let
  \begin{equation}
    \bar \s :=\frac{\s(\x)+ \x}{2} + \frac{\beta}{8} \| \s(\x) - \x\|^2 \frac{\vv}{\|\vv\|_\X} \in \X.
  \end{equation}
  Then 
  \begin{align}
    \prodscal{\s(\x) - \x}{- \nabla f(\x)} 
    & \geq \prodscal{ \bar \s - \x}{- \nabla f(\x)} \notag \\
    &= \frac{1}{2}\prodscal{\s(\x) - \x}{- \nabla f(\x)} +  \frac{\beta}{8}  \|\s(\x) - \x\|^2 \frac{\|\nabla f (\x)\|_{\X^*}}{\|\vv\|_\X} \notag \\
    &\geq \frac{1}{2}\prodscal{\s(\x) - \x}{- \nabla f(\x)} +  \frac{\beta}{8} \delta \|\s(\x) - \x\|^2 
  \end{align}
  which leads us to the desired result.
  \endproof
  From this lemma, under the assumption that $\min(\|\nabla_{\!x} L(\z)\|_{\X^*}, \|\nabla_{\!y} L(\z)\|_{\Y^*}) \geq \delta$ $\forall \z \in \X \times \Y$, it directly follows that 
  \begin{align}
    \gap^\FW = \gap^{(x)} + \gap^{(y)} &\geq \frac{\beta}{4} \delta \left( \|\st_x - \xt\|_\X^2 + \| \st_y - \yt\|_\Y^2\right) \notag \\ 
    &\geq \frac{\beta}{8} \delta \left( \|\st_x - \xt\|_\X + \| \st_y - \yt\|_\Y\right)^2 = \frac{\beta}{8} \delta  \|\st - \zt\|_{\X \times \Y}^2. \label{eq:gap_lower_dist}
  \end{align}
  Now we recall the convergence theorem for strongly convex sets from the main text, Theorem~\ref{thm:conv_strong}:
  \begin{reptheorem}{thm:conv_strong}
       Let $\L$ be a convex-concave function and $\X$ and $\Y$ two compact $\beta$-strongly convex sets. 
       Assume that the gradient of  $\L$ is $L$-Lipschitz continuous and that there exists $\delta>0$ such that $\min(\|\nabla_{\!x} L(\z)\|_{\X^*}, \|\nabla_{\!y} L(\z)\|_{\Y^*}) \geq \delta \;\, \forall \z \in \M$. Set $C_\delta := 2L + \frac{8L^2}{\beta \delta}$. 
       Then the gap $\gap^{\FW}$~\eqref{eq:gap} of the SP-FW algorithm with step size $\gamma_t = \tfrac{\gap^{\FW}}{\|\st-\zt\|^2 C_\delta}$
       converges linearly as 
     \begin{equation}
       \gap^{\FW} \leq g_0 \left( 1- \rho \right)^{t}    
     \end{equation}
    where $\rho :=  \frac{\beta \delta}{16 C_\delta}$. The initial gap $g_0$ is cheaply computed during the first step of the SP-FW algorithm. Alternatively, one can use the following upper bound to get uniform guarantees:
    \begin{equation}  
     g_0 \leq \sup_{\z\in \M} \|\nabla_x \L(\z)\|_{\X^*} D_\X + \sup_{\z\in \M} \|\nabla_y \L(\z)\|_{\Y^*} D_\Y.
     \end{equation} 
    \end{reptheorem}

     \proof
      We compute the following relation on the gap:
      \begin{align}
      g_{t+1} 
      & = \prodscal{\ztt - \stt}{F(\ztt)} \notag\\
      & = \prodscal{\zt - \stt}{F(\ztt)} + \stepsize_t\prodscal{\st - \zt}{F(\ztt)}\notag\\
      & = \prodscal{\zt - \stt}{F(\zt)} +  \prodscal{\zt - \stt}{F(\ztt)-F(\zt)} \notag\\
      &\quad + \stepsize_t\prodscal{\st - \zt}{F(\zt)} + \stepsize_t \prodscal{\st-\zt}{F(\ztt)-F(\zt)} \notag\\
      & \leq \prodscal{\zt - \stt}{F(\zt)} +  \prodscal{\zt - \stt}{F(\ztt)-F(\zt)} \notag\\
      &\quad + \stepsize_t\prodscal{\st - \zt}{F(\zt)} + \stepsize_t^2 \|\st-\zt\|^2 L \label{al:lip gradient}
      \end{align}
      where in the last line we used the fact that the function $F(\cdot)$ is Lipschitz continuous. Then using that $\prodscal{\zt - \stt}{F(\zt)} \leq \prodscal{\zt - \st}{F(\zt)}$ (by definition of $\st$), we get
      \begin{align}
      g_{t+1} 
      & \leq \gap(1- \stepsize_t) + \prodscal{\zt - \stt}{F(\ztt)-F(\zt)} +\stepsize_t^2 \|\st-\zt\|^2 L  \notag\\
      & \leq \gap(1- \stepsize_t) + \prodscal{\st - \stt}{F(\ztt)-F(\zt)} +\stepsize_t^2 \|\st-\zt\|^2 L. \label{eq:183}
      \end{align}
      The last line uses the fact that $F$ is monotone by convexity (Equation~\eqref{eq:nondecreasing}). Finally, using once again the Lipschitz continuity of $F$ and the one of $\s(\cdot)$ (by Theorem~\ref{thm:s_lip}), we get
      \begin{align}
      \prodscal{\st- \stt}{F(\ztt) - F(\zt)} 
      & \leq \| \st - \stt\| L \| \ztt - \zt\| \notag\\
      & \leq \frac{4L^2}{\beta \delta} \|\ztt - \zt \|^2 \qquad\quad (\text{Lipschitz continuity of } \s)\notag\\
      & = \frac{4L^2}{\beta \delta} \stepsize_t^2 \|\st - \zt \|^2. \label{eq:184} 
      \end{align}
      Combining~\eqref{eq:184} with~\eqref{eq:183}, we get
      \begin{equation}
        \label{eq:ineg_strongly convex_sets}
      g_{t+1} \leq \gap (1- \stepsize_t) + \stepsize_t^2 \|\st - \zt \|^2 \frac{C_\delta}{2} \qquad \text{where} \qquad C_\delta := 2L  + \frac{8L^2}{\beta \delta}.
      \end{equation}
      Thus by setting the step size $\stepsize_t = \frac{\gap}{\|\st - \zt\|^2 C_\delta}$,
       we get 
     \begin{equation}
        g_{t+1} \leq \gap - \frac{\gap}{2 C_\delta} \left( \frac{\gap}{\|\st-\zt\|^2} \right) \leq \gap \left( 1- \frac{\beta \delta}{16 C_\delta} \right), 
     \end{equation}
     using the fact that the gap is lower bounded by a constant times the square of the distance between $\st$ and $\zt$ (Equation~\eqref{eq:gap_lower_dist}).
    \endproof
    Note that the bound in this theorem is not affine invariant because of the presence of Lipschitz constants and strong convexity constants of the sets. The algorithm is not affine invariant either because the step size rule depends on these constants as well as on $\|\st-\zt\|$. Deriving an affine invariant step size choice and convergence analysis is still an interesting open problem in this setting.

\section{Details on the experiments} %
\label{sec:details_on_the_experiments}

  \paragraph{Graphical Games.} %
  \label{par:graphical_games}
  The payoff matrix $M$ that we use encodes the following simple
model of competition between universities with their respective benefits: 
\begin{enumerate}
  \item University~1 (respectively University~2) has benefit $b_i^{(1)}$ ($b_i^{(2)}$) to get student $i$. 
  \item Student $i$ ranks the possible
roommates with a permutation $\sigma_i \in \mathcal S_p$. Let $\sigma_i(j)$ represents the rank of $j$ for $i$ (first in the list is the preferred one). 
  \item They go to the university that matched them with their preferred roommate, in case of equality the student chooses randomly.
  \item Supposing that $\x$ encodes the roommate assignment proposed by University~1 (and $\y$ for University~2), then the expectation of the benefit of University~1 is $\x^\top M \y$,
with the following definition for the payoff matrix $M$ indexed by pairs
of matched students. For the pairs $(i,j)$ with $i < j$ and $(k,l)$ with $k < l$ with elements in $1, \ldots, s$, we have:
\begin{enumerate}
  \item $M_{ij,il} = \left\{ 
            \begin{array}{lll}
            b_i^{(1)} & \text{if} \; \sigma_i(j) < \sigma_i(l) \quad \text{\emph{i.e. student $i$ preferred $j$ over $l$}}\\
            -b_i^{(2)} & \text{if} \; \sigma_i(j) > \sigma_i(l) \\
            \frac{b_i^{(1)} - b_i^{(2)}}{2} & \text{otherwise}  \quad (\text{in that case }j = l).
            \end{array} \right.$
  \item $M_{ij,kj} = M_{ji,jk}$ %
  \item $M_{ij,ki} = M_{ij,ik}$
  \item $M_{ij,jl} = M_{ij,lj}$
  \item $M_{ij,kl} = 0 \quad \text{otherwise}$
\end{enumerate}
\end{enumerate}
Note that we need to do unipartite matching here (and not bipartite
matching) since we have to match students together and not students with dorms. 

For our experiments, in order to get a realistic payoff matrix, we set $\mu_i \sim \mathcal U[0,1]$ the \emph{true} value of student $i$. Then we set $b_i^{(U)}\sim \mathcal N(\mu_i, 0.1)$ the value of the student $i$ \emph{observed} by University~$U$. To solve the perfect matching problem, we used Blossom V by \citetsup{kolmogorov2009blossom}. 
  \paragraph{Sparse structured SVM.} %
  \label{par:Structured SVM}
  We give here more details on the derivations of the objective function for the structured SVM problem.
  We first recall the structured prediction setup with the same notation from~\citep{lacoste2013block}. In structured prediction, the goal is to predict a structured object $\y \in \Y(\x)$ (such as a sequence of tags) for a given input $\x \in \X$. For the structured SVM approach, a structured feature map $\bm{\phi} : \X \times \Y \rightarrow \R^d$ encodes the relevant information for input / output pairs, and a linear classifier with parameter $\bm{w}$ is defined by $h_{\bm{w}} = \argmax_{\y \in \Y(\x)} \prodscal{\bm{w}}{\bm{\phi}(\x, \y)}$. We are also given a task-dependent structured error $L(\y', \y)$ that gives the loss of predicting $\y$ when the ground truth is $\y'$. Given a labeled training set $\{(\x^{(i)}, \y^{(i)})\}_{i=1}^n$, the standard $\ell_2$-regularized structured SVM objective in its non-smooth formulation for learning as given for example in Equation~(3) from \citep{lacoste2013block} is:
  \begin{equation}
    \min_{ \bm{w} \in \R^d} \frac{\lambda}{2}\|\bm{w}\|_2^2 + \frac{1}{n} \sum_i \tilde H_i(\bm{w})
  \end{equation}
  where $\tilde H_i(\bm{w}) := \max_{\y \in \mathcal Y_i} \; L_i(\y) - \prodscal{\bm{w}}{\bm{\psi}_i(\y)}$ is the structured hinge loss, and the following notational shorthands were defined: $\Y_i := \Y(\x^{(i)})$, $L_i(\y) := L(\y^{(i)}, \y)$ and $\bm{\psi}_i(\y) := \bm{\phi}(\x^{(i)},\y^{(i)}) - \bm{\phi}(\x^{(i)},\y)$.
  
  In our setting, we consider a sparsity inducing $\ell_1$-regularization instead. Moreover, we use the (equivalent) constrained formulation instead of the penalized one, in order to get a problem over a polytope. We thus get the following challenging problem:
  \begin{equation}
    \min_{\|\bm{w}\|_1  \leq R} \frac{1}{n} \sum_i \tilde H_i(\bm{w}).
  \end{equation}
  To handle any type of structured output space $\Y$, we use the following generic encoding.
  Enumerating the elements of $\Y_i$, we can represent the $j^{th}$ element of $\Y_i$ as $(\overbrace{0,\ldots,0}^{j-1},1,0,\ldots,0 ) \in \R^{|\Y_i|}$. Let $M_i$ have $\big(\bm{\psi}_i(\y)\big)_{\y \in \mathcal Y_i}$ as columns and let $\bm{L}_i$ be a vector of length $|\Y_i|$ with $L_i(\y)$ as its entries. 
  The functions $\tilde H_i(\bm{w})$ can then be rewritten as the maximization of linear functions in $\y$: $\tilde H_i(\bm{w}) = \max_{\y \in \mathcal Y_i} \; \bm{L}_i^\top \y - \bm{w}^\top M_i \y$.
  As the maximization of linear functions over a polytope is always obtained at one of its vertex, we can equivalently define the maximization over the convex hull of $\Y_i$, which is the probability simplex in $\R^{|\Y_i|}$ that we denote $\Delta(|\Y_i|)$: 
  \begin{equation}
    \max_{\y_i \in \mathcal Y_i} \; \bm{L}_i^\top \y_i - \bm{w}^\top M_i \y_i =  \max_{\bm{\alpha}_i \in \Delta(|\mathcal Y_i|)} \!\!\bm{L}_i^\top \bm{\alpha}_i - \bm{w}^\top M_i \bm{\alpha}_i
  \end{equation}
  Thus our equivalent objective is 
  \begin{equation}
   \min_{\|\bm{w}\|_1  \leq R} \frac{1}{n} \sum_i \Big( \max_{\y_i \in \mathcal Y_i} \; \bm{L}_i^\top \y_i - \bm{w}^\top M_i \y_i \Big) 
    = 
    \min_{\|\bm{w}\|_1 \leq R} \frac{1}{n} \sum_i \Big( \max_{\bm{\alpha}_i \in \Delta(|\mathcal Y_i|)} \!\!\bm{L}_i^\top \bm{\alpha}_i - \bm{w}^\top M_i \bm{\alpha}_i \Big) ,
  \end{equation}
  which is the bilinear saddle point formulation given in the main text in~\eqref{eq:structuredSVM}.

\bigskip
\bigskip
\bibliographystylesup{abbrvnat}
\bibliographysup{bib}

\end{document}